\documentclass[12pt]{amsart}
\usepackage[margin=3cm]{geometry}                % See geometry.pdf to learn the layout options. There are lots.
\usepackage{graphicx}
\usepackage{amssymb}
\usepackage{epstopdf}
\usepackage{amsthm}
\usepackage{enumerate}
\usepackage{hieroglf}
\usepackage{bbding}
\usepackage{xcolor}
\usepackage{framed}
\usepackage{booktabs}
\usepackage{csquotes}
\usepackage{hyperref}
\usepackage[noabbrev]{cleveref}
\usepackage{float}
\usepackage{fancyhdr}
\usepackage{tabularx}

\usepackage[textwidth=2cm]{todonotes}

\definecolor{skyblue}{rgb}{0.85,0.85,1}
\usepackage[normalem]{ulem} % Prevents \emph from being underlined
\usepackage[final]{showkeys}

\numberwithin{equation}{section}

\newcommand{\e}{{\rm e}}
\newcommand{\cL}{\mathcal{L}}
\newcommand{\ve}{\varepsilon}
\newcommand{\wh}{\widehat}

\DeclareMathOperator{\sech}{sech}
\newcommand{\Hbb}{\mathbb{H}}
\newcommand{\R}{\mathbb{R}}
\newcommand{\C}{\mathbb{C}}
\newcommand{\imag}[1]{\textrm{Im}\left( #1\right)}

\definecolor{internationalkleinblue}{rgb}{0.0, 0.18, 0.65}

\newcommand{\black}{\textcolor[rgb]{0.0,0.0,0.0}}

\newtheorem*{remark}{Remark}
\title[]{Instability of solutions in a degenerate reaction diffusion equation\\ \today}
\author{R. Marangell}
\address{R. Marangell \\ School of Mathematics and Statistics, The University of Sydney.}

\author{J.J. Wylie}
\address{J.J. Wylie \\ Department of Mathematics, City University of Hong Kong.}

\author{B.H. Bradshaw-Hajek}
\address{B.H. Bradshaw-Hajek\\School of Mathematical Sciences, Adelaide University.}

\begin{document}
\pagestyle{fancy}
\footskip=30pt
                                           
\maketitle

\begin{abstract}
We study the spectral stability of travelling and stationary front and pulse solutions in a class of degenerate reaction–diffusion systems. We characterise the essential spectrum of the linearised operator in full generality and identify conditions under which it lies entirely in the left half-plane. For a number of special cases we obtain analytical results, including explicit Evans functions and complete spectral descriptions for certain stationary waves. In regimes where analytical methods are not available, we compute the point spectrum numerically using a Riccati–Evans function approach. Our results show that stable travelling fronts can occur, while travelling pulses are generically unstable.
\end{abstract}

\tableofcontents

\section{Introduction}

In this paper we will consider a system of reaction-diffusion equations of the form
\begin{eqnarray}\label{eq:sys1}
\begin{aligned}
u_t&=~~u_{xx}+g(u,v)\\
v_t&=Dv_{xx}-g(u,v),
\end{aligned}
\end{eqnarray}
Equations of this type are degenerate in the sense that the reaction terms are linearly dependent and so the only requirement for homogeneous steady states is a single equation $g(u,v)=0$ for two quantities. Such models arise very naturally for systems that have the property that the reaction terms simply move populations from one category to another without either creating or destroying the local total population. That is, any local increase in $u$ due to the reaction terms must necessarily be accompanied by a subsequent local decrease in $v$.
 
This type of degeneracy is found in a number of important applications. Perhaps the best-known example is the SIS model of infection \cite{Edel2005} in which the total population is assumed to be one of two classes: infected or susceptible. The SIS model neglects births and deaths, and moves individuals between the susceptible and infected classes via infection and recovery processes. It therefore locally conserves the total population. In fact, equations of this type are extremely widespread, and important examples can be found across a diverse range of applications including biology \cite{Kawa997,Sat2001,Mori2008,Volp2009,HuangWM2011}, 
chemistry \cite{Saul1984,Merk1989,Bill1991,Merk1993}, and physics \cite{Huang2007,Huang2008,Wylie2009}.

Another application in which degenerate reaction terms of a slightly more general type naturally arise is in the Tuckwell-Miura (TM) equations that model the ionic transport believed to be responsible for a neurological phenomenon that is associated with migraine \cite{Tuck1978}. The TM model tracks the intra-cellular and extra-cellular ion concentrations of potassium and calcium ions and therefore contains four equations. Each of the ion types can be transported from the intra-cellular and extra-cellular regions via various cellular processes that neither create nor destroy ions. Therefore, the TM model contains two similar degeneracies to \eqref{eq:sys1}. Moreover, the diffusion coefficient for the intra-cellular ions is assumed to be zero since intra-cellar ions are constrained by the boundaries of the cell. Consequently, intra-cellular ions can only be transported spatially by first being transported into the extra-cellular space and then diffusing there. So, the TM model has the additional property that $D=0$. In fact, there are numerous models involving ion transport that share similar properties \cite{Yao2011,HuangMY2011,Yao2018}. Another simple example in which $D=0$ is the limit of the SIS model for which infected individuals are so badly affected by the illness that they are assumed to have no mobility \cite{Edel2005}. The case of $D=0$ was also considered in the context of bacterial pattern formation by Satnoianu et al. \cite{Sat2001}
 
Despite the physical importance and prevalence of systems similar to \eqref{eq:sys1}, it is extremely surprising that there has been relatively little previous work and much of the pioneering work has been restricted to relatively simple reaction terms. However, Wylie and Miura \cite{Wylie2006} considered a general type of nonlinear reaction term and determined conditions for the existence of travelling fronts. They also examined the types of localised perturbation to a stable rest state that can trigger a traveling front. Bradshaw-Hajek and Wylie \cite{BHW19} derived explicit conditions on the reaction term to guarantee the existence of travelling and stationary pulse and front solutions. This allowed them to obtain numerous examples of reaction terms for which analytical pulse and front solutions can be found. This had important consequences for the TM model which was previously believed to be the minimal model of this degenerate form that could exhibit travelling pulse solutions. However, although travelling pulse solutions were shown to exist, the stability of such solutions has not been studied and therefore the minimality of the TM model is still an open question. One of the main aims of this paper is to address this question. Bradshaw-Hajek and Wylie \cite{BHW19} also showed that the case $D=0$ is singular and that infinite families of stationary piecewise constant solutions can occur. This is in direct contrast with the $D\ne 0$ case for which all stationary solutions must be continuous.

\black{In this paper we will consider the spectral stability of travelling and stationary front and pulse solutions of \eqref{eq:sys1}, by determining the spectrum of the operator linearised around such solutions, which we will denote by $\sigma(\cL)$. 
For most of our examples of interest, the essential spectrum will be in the left half plane, typically touching the origin with a (generic) quadratic tangency, and thus will not contribute to instability. What will determine the stability in most cases is the {\em point spectrum}.  Owing to our assumptions, we have that the point spectrum will consist of a discrete set of eigenvalues of finite (algebraic) multiplicity \cite{KP13}. In the cases where we cannot make analytical progress, we will numerically compute the point spectrum using the {\em Evans function} or, more likely, the {\em Riccati Evans function}.}

The term Evans function was coined by Jones in \cite{Jon84} to describe the Wronskian-like determinant used to locate eigenvalues via a shooting argument. The idea is for values of the spectral parameter away from the essential spectrum of $\cL$, we evolve a spanning set of the stable subspaces of the equivalent, far-field linearised systems  and then take the determinant of the matrix with such evolved subspaces evaluated at a common point $z_0$. When this determinant vanishes, we have a common element which decays at both ends, and so we have an eigenvalue.

In this paper we will consider the spectral stability of travelling and stationary front and pulse solutions of \eqref{eq:sys1}. For $D=0$, the only non-trivial {\it stationary} solutions are composed of piecewise constants. We show that the eigenvalues associated with these solutions 
can be obtained analytically. 
For $D=0$ the stability of {\it travelling waves} must be determined numerically and we use an efficient, robust and accurate method to do so. We use this to determine the stability of all known analytical solutions and will hence show that linearly stable {\it travelling fronts} can exist. However, all of the {\it travelling pulses} that we examine prove to be linearly unstable. This shows that even though the TM model is not the minimal model for the existence of travelling pulses, it may indeed be the minimal model for the existence of {\it stable} travelling pulses. For the case $D=1$ we find {\it stationary fronts} for which we obtain an analytical expression for the Evans function and the entire spectrum. Hence we prove that such solutions can be spectrally, linearly and nonlinearly stable. However, we do find {\it stable stationary pulses}. For $D\ne0,1$, all stability results must be determined numerically. In this case we show that the most unstable eigenvalues can be complex. We again use our numerical scheme to test the stability of all known analytical solutions.

\section{Formulation of the problem}

Since we will be interested in travelling waves, we define $z=x-ct,~\tau=t$ as the travelling wave coordinates, where $c$ is the constant wave speed. Equations \eqref{eq:sys1} are then recast in terms of these travelling wave coordinates to obtain
\begin{eqnarray}\label{eq:TWsys1}
\begin{aligned}
u_{\tau}-cu_z&=u_{zz}+g(u,v)\,,\\
u_{\tau}-cv_z&=Dv_{zz}-g(u,v)\,.
\end{aligned}
\end{eqnarray}
Steady travelling waves which maintain their shape over time satisfy $u_{\tau}\equiv v_{\tau}\equiv0$. We will denote such solutions with hats  so that
\begin{eqnarray}\label{eq:TWsys2}
\begin{aligned}
\hat u''+c\hat u'+g(\hat u,\hat v)&=0\,,\\
D\hat v''+c\hat v'-g(\hat u,\hat v)&=0\,,
\end{aligned}
\end{eqnarray}
where $'$ means $\frac{d}{dz}$. 

We define a {\it front} to be a solution for which the constant state far ahead of the wave $(\hat u_+,\hat v_+)=\lim_{z\to\infty}(\hat u(z),\hat v(z))$ is different from the constant state far behind the wave $(\hat u_-,\hat v_-)=\lim_{z\to-\infty}(\hat u(z),\hat v(z))$. We consider a {\it pulse} to be a solution for which the constant states far ahead and far behind the waves are the same, that is, $(\hat u_+,\hat v_+)=(\hat u_-,\hat v_-)$. (Note this includes, so-called multi-pulses.) We define a wave to be {\it stationary} (pulses or fronts) if $c=0$ and {\it nonstationary} if $c\ne0$.

In order to determine the stability of solutions to equations of type \eqref{eq:TWsys2}, we consider the spectrum of the relevant linearised operators. Starting with (appropriate) solutions $\hat{u}(z)$ and $\hat{v}(z)$, we consider solutions to equations \eqref{eq:sys1} of the form $u(x,t) = \hat{u}(z) + \ve \e^{\lambda t} p(z)$ and $v(x,t) = \hat{v}(z) + \ve\e^{\lambda t} q(z)$. Substituting into equations \eqref{eq:sys1}, and keeping the leading (first order in $\ve$) terms we are led to the linear system of non-autonomous ordinary differential equations,
\begin{equation}\label{eq:eigenvalue1}
\begin{aligned}
\lambda p & = p'' + c p ' + g_u(\hat{u}(z),\hat{v}(z))p + g_v(\hat{u}(z),\hat{v}(z))q \\ 
\lambda q & = D q'' + c q' - g_u(\hat{u}(z),\hat{v}(z))p - g_v(\hat{u}(z),\hat{v}(z))q.
\end{aligned}
\end{equation}
{From here onward, for notational brevity, we will not explicitly write the dependence of $g_u$ and $g_v$ on  $\hat{u}(z)$ and $\hat{v}(z)$ when it is clear. }

\subsection{Essential spectra}
{
In order to study the essential spectra of the operator, we need to consider the behaviour of solutions of the system \eqref{eq:eigenvalue1} at $z\to\pm\infty$.
Since this is a linear fourth-order system, there will be four linearly independent solutions.
In order to determine the behaviour of the solutions as $z\to\pm\infty$ we take the limit of \eqref{eq:eigenvalue1} and obtain a fourth-order system with constant coefficients. We can thus determine the behaviour by examining the characteristic equation obtained by considering solutions of the form $e^{\nu^\pm z}$ where the $\pm$ superscript indicates the limits as $z\to\pm\infty$. This gives
\begin{align}\label{Characteristic}
D (\nu^{\pm})^4 
+ c(D+1)(\nu^{\pm})^3 
&+ \left[c^2-(\hat g_v^\pm -D \hat g_u^\pm)-\lambda(D+1)\right] (\nu^{\pm})^2\\ \nonumber
&- c\left[\hat g_v^\pm- \hat g_u^\pm+2\lambda \right] \nu^{\pm}
+\lambda\left[ \hat g_v^\pm- \hat g_u^\pm+\lambda\right] =0,
\end{align}
where
\begin{equation}
\hat g_u^\pm=\lim_{z\to\pm\infty} g_u\left(\hat u(z), \hat v(z) \right)
\qquad\text{and}\qquad
\hat g_v^\pm=\lim_{z\to\pm\infty} g_v\left(\hat u(z), \hat v(z) \right).
\end{equation}}

{
Since we require that the eigenfunctions are bounded at $z\to\pm\infty$ we must set the coefficients of any of the four linearly independent solutions to zero for which either $\Re(\nu^{+})>0$ or $\Re(\nu^-)<0$.
This means that each value of $\nu$ with $\Re(\nu^{+})>0$ or $\Re(\nu^-)<0$ provides a homogeneous constraint on 
 the coefficients of the four linearly independent solutions. The mismatch (including potentially none) between the number of $\nu^+$ with $\Re({\nu^+})>0$ versus the number of $\nu^-$ with $\Re({\nu^{-}})>0$ is called {\em the Fredholm index}. In the cases considered here, we have four constraints and for $\nu$ with large enough real part the only solution will be to set all of the coefficients to be zero. Having four constraints for four coefficients corresponds to the Fredholm index being zero. 
 However, we may be able to find isolated values of $\lambda$ for which non-trivial solutions exist. Such solutions correspond to the point spectra of the operator. 
If we have more or fewer than four constraints (i.e. the Fredholm index is non-zero) then the operator is not invertible and such values of $\lambda$ are part of the essential spectrum.
The {\em Fredholm borders} are typically defined to be the part of the essential spectrum for which at least one value of $\nu$ is purely imaginary. }

In order to determine the essential spectrum we need  to determine the number of roots of \eqref{Characteristic} with $\Re(\nu)>0$ as $\lambda$ varies across the whole complex plane.
In order to do this we first note that \eqref{Characteristic} is invariant under the symmetry $c\to-c$ and $\nu\to-\nu$. This symmetry is associated with reflection in $x$. There is also another symmetry $D\to 1/D$, $\nu\to D^{1/2}\nu$, $c\to D^{-1/2}c$, $\hat g_u^\pm\to-\hat g_v^\pm$ and $\hat g_v^\pm\to-\hat g_u^\pm$. This symmetry is associated with scaling $x$ to make the diffusion in the $v$ equation equal to unity and relabelling the two species $u$ and $v$.
Given these two symmetries, we need only to consider the behaviour for $c\ge0$ and $0\le D\le1$.

We begin by considering the Fredholm borders -- those values $\nu$ for which one of the roots of \eqref{Characteristic} is purely imaginary.
Setting $\nu=ik$ where $k$ is real into \eqref{Characteristic} and grouping terms we obtain
\begin{equation}
\left(\lambda-ikc\right)^2+\left[(D+1)k^2+\hat g_v^\pm-\hat g_u^\pm \right]\left(\lambda-ikc\right)
+k^2\left[Dk^2 +\hat g_v^\pm-D\hat g_u^\pm\right]=0,
\label{FredholmBorders}
\end{equation}
which is a quadratic in $(\lambda-ikc)$ with real coefficients.
Stability of the Fredholm borders requires that $\Re(\lambda)<0$ for all values of $k$.
This implies that $(D+1)k^2+\hat g_v^\pm-\hat g_u^\pm>0$ and $Dk^2 +\hat g_v^\pm-D\hat g_u^\pm>0$ for all $k$. Hence, this part of the essential spectrum will be stable if $\hat g_v^\pm>\hat g_u^\pm$ and $\hat g_v^\pm>D\hat g_u^\pm$. On the one hand, if $\hat g_v^\pm<\hat g_u^\pm$, then the uniform solution corresponding to the far field behaviour of $u$ and $v$ is unstable to a direct instability with zero wavenumber being the most unstable mode. On the other hand, if $\hat g_v^\pm<D\hat g_u^\pm$, then the uniform solution is unstable to a Turing instability whose most unstable mode occurs at a non-zero wavenumber.

We can readily solve \eqref{FredholmBorders} to obtain an expression for $\lambda$ given by
 \begin{equation}
\lambda=\frac{2ikc-(D+1)k^2-\hat g_v^\pm+\hat g_u^\pm \pm \sqrt{\Delta}}{2},\label{FredholmBorder}
\end{equation}
where 
\begin{align}
\Delta&=\left[(D+1)k^2+\hat g_v^\pm-\hat g_u^\pm \right]^2-4k^2\left[Dk^2 +\hat g_v^\pm-D\hat g_u^\pm\right]\\ 
         &=\left[(D-1)k^2+\hat g_v^\pm+\hat g_u^\pm\right]^2-4\hat g_v^\pm\hat g_u^\pm.\nonumber
\end{align}
We immediately see that if $D=1$ or if $\hat g_v^\pm\hat g_u^\pm\le0$ then $\Delta\ge0$.
Furthermore, noting that we only need to consider $0\le D\le1$ because of symmetry and observing that $\Delta$ can be written in the form $\Delta=\left[(D-1)k^2\right]^2+2(\hat g_v^\pm+\hat g_u^\pm)\left[(D-1)k^2\right]
+ (\hat g_v^\pm-\hat g_u^\pm)^2$ we see that if $\hat g_u^\pm<0$ and $\hat g_v^\pm<0$ then $\Delta\ge0$. Hence, $\Delta$ can only be negative for $0\le D\le1$ if both $\hat g_u^\pm>0$ and $\hat g_v^\pm>0$.

\begin{figure}
	\centering
	\renewcommand{\arraystretch}{0.1} % it makes the space between two rows
	\begin{tabularx}{\textwidth}{c c}%{c{6cm}*4{Z}}  % C{6cm} this means the gap between columns  4{Z} means it creates column               
		\includegraphics[scale=0.35]{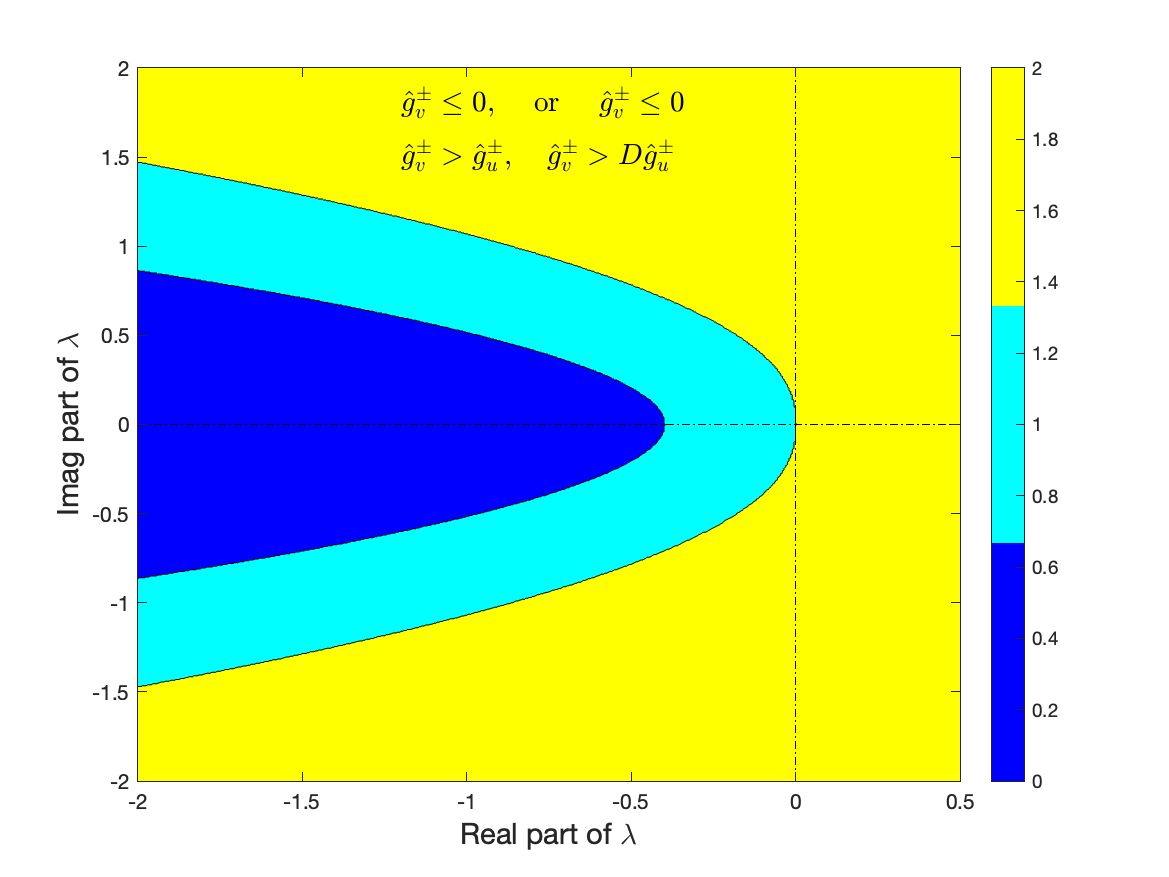}&\includegraphics[scale=0.35]{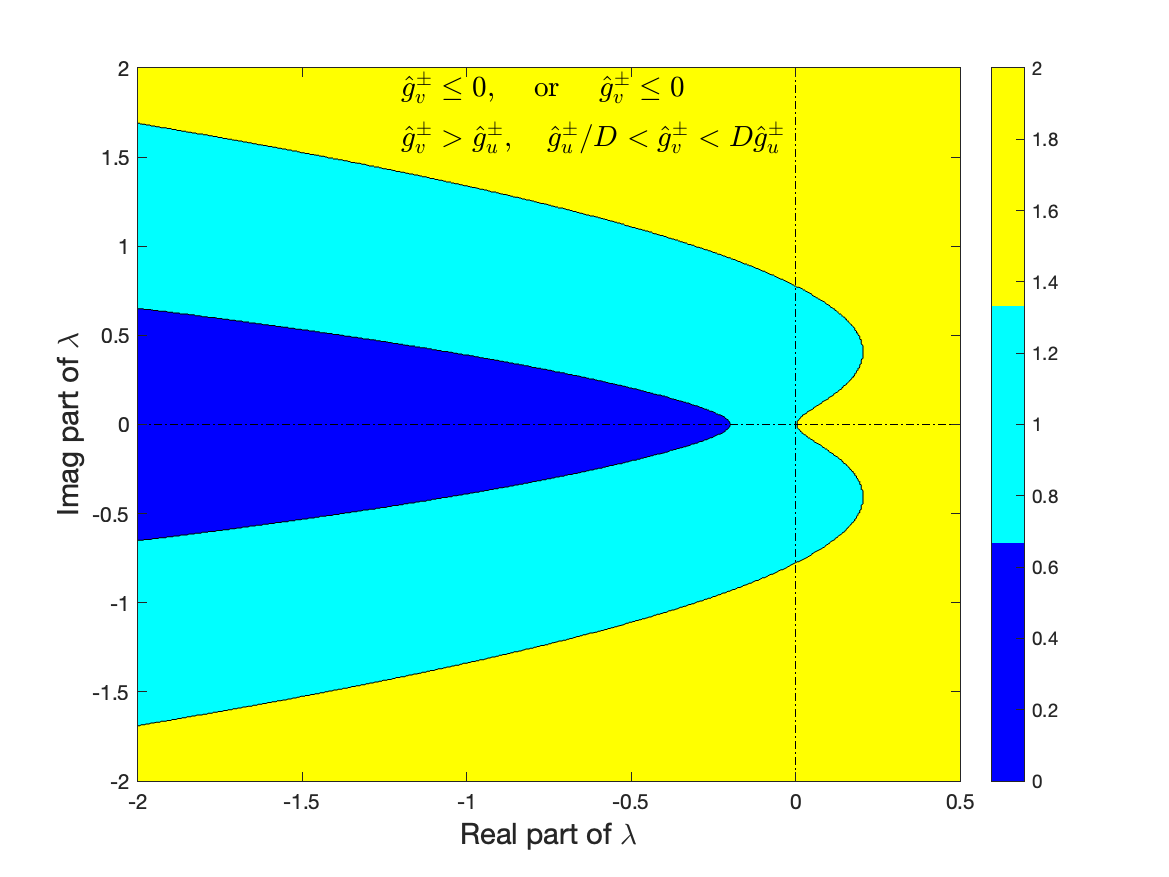}\\
        (a)                               &  (b) \\
		\includegraphics[scale=0.35]{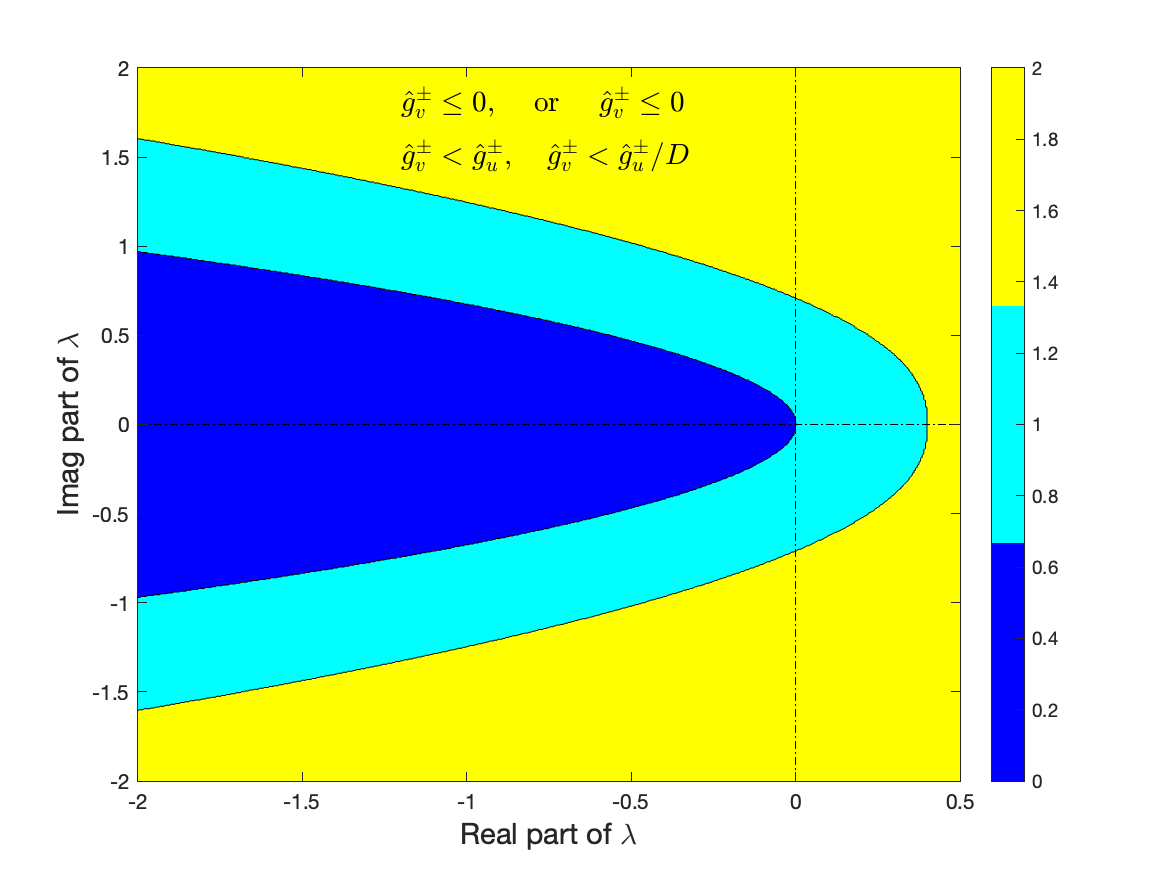}&\includegraphics[scale=0.35]{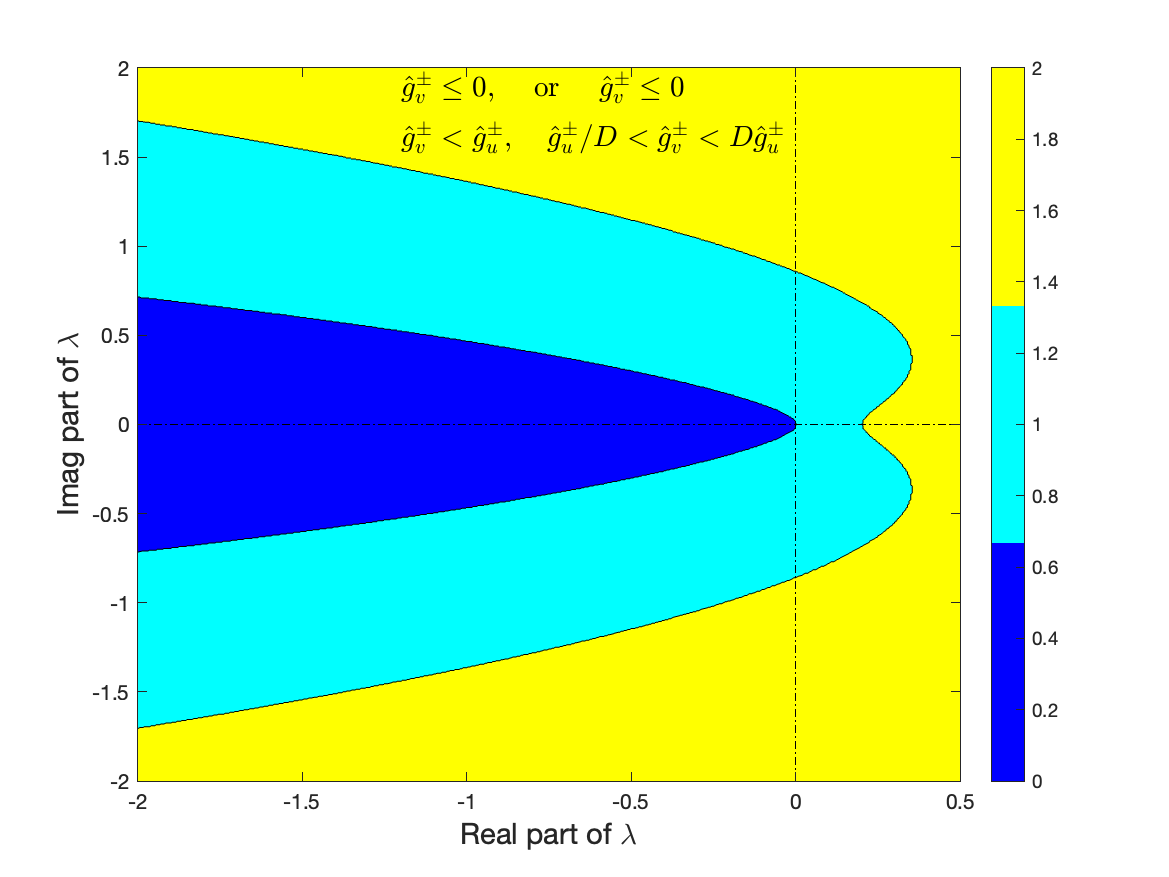}\\
		(c)                               &  (d) \\
		\includegraphics[scale=0.35]{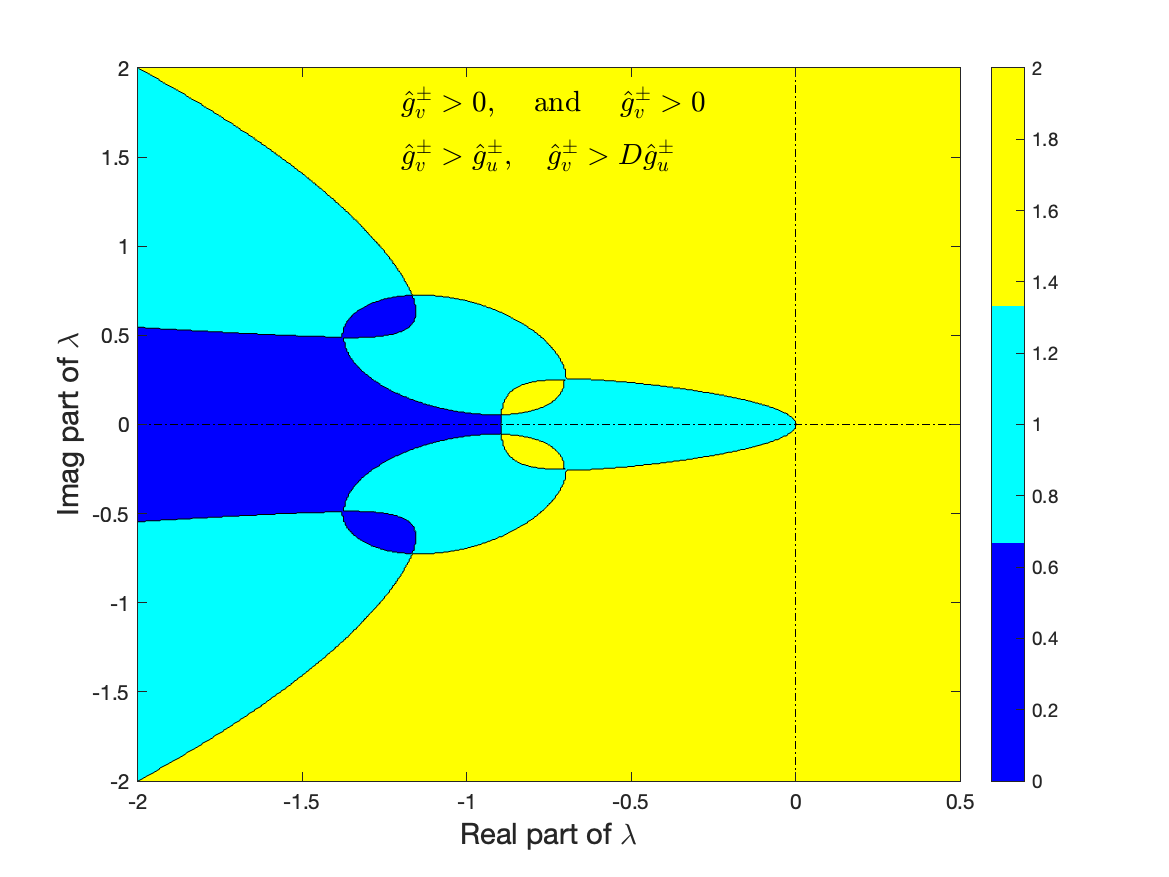}&\includegraphics[scale=0.35]{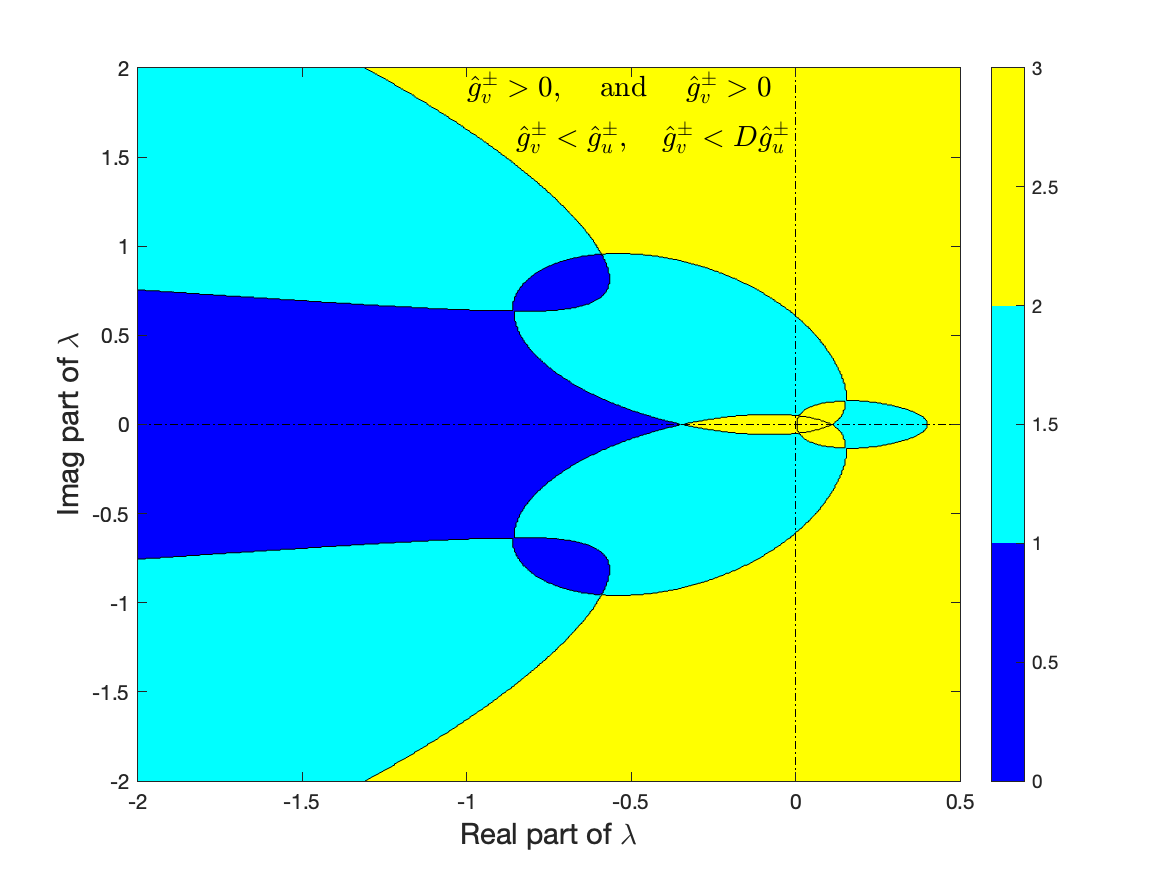}\\	
		(e)                               &  (f) \\
	\end{tabularx}
	\caption{The six types of typical Fredholm borders and the number of roots with positive real part for $\nu$ are shown for various values of $\hat g_u^\pm$ and $\hat g_v^\pm$. The results are shown for $c>0$ and $0<D<1$. Results for $c<0$ and $D>1$ can be obtained using symmetry considerations.} 
	\label{fig:FredholmBoundaries}
\end{figure}

If $\Delta\ge0$, then 
 \begin{equation}
\Im(\lambda)=kc\quad\text{and}\quad\Re(\lambda)=\frac{-(D+1)k^2-\hat g_v^\pm+\hat g_u^\pm \pm \sqrt{\Delta}}{2}.
\end{equation}
In this case, if $c=0$, then the Fredholm borders will be purely real and given by $\lambda\in\left(-\infty,\max\{0,\hat g_u^\pm- \hat g_v^\pm\}\right]$. This means that all \black{temporal} modes in the essential spectrum will be non-oscillatory for stationary waves for $\hat g_u^\pm\le0$ or $\hat g_v^\pm\le0$. 
For $c\ne0$ there is a one-to-one correspondence between the wavenumber $k$ and $\Im(\lambda)$. Therefore, each of the two roots given in \eqref{FredholmBorder} will have the property that the
$\Re(\lambda)$ must be a single-valued function of $\Im(\lambda)$. This can be observed in Figure~\ref{fig:FredholmBoundaries}a--d.
One can readily see that $\Re(\lambda)$ is symmetric in $k$ and so there must be a stationary point at $\Im(\lambda)=0$. However, there can be other values of $k$ for which $\Re(\lambda)$ has a local extremum. After some straightforward algebra we obtain an expression for the possible values of $k$ 
for which stationary points can occur
\begin{equation}
k=0\quad\text{or}\quad k^2=\frac{\hat g_v^\pm+\hat g_u^\pm}{1-D}\pm\frac{1+D}{1-D}\sqrt{\frac{\hat g_v^\pm\hat g_u^\pm}{D}}.
\end{equation}
Hence, there will be no stationary points with $k\ne0$ for $\hat g_v^\pm\hat g_u^\pm<0$.
For $0<D<1$, $\hat g_v^\pm<0$, and $\hat g_u^\pm<0$ the expression for $k^2$ can be factorised as
\begin{equation}
k^2=-\frac{\left[\sqrt{D}\sqrt{-\hat g_v^\pm}\pm\sqrt{-\hat g_u^\pm} \right]
\left[\sqrt{-\hat g_v^\pm}\pm\sqrt{D}\sqrt{-\hat g_u^\pm} \right]}{(1-D)\sqrt{D}}.
\end{equation}

One of the values of $k^2$ is necessarily negative, but the other value will be positive if $D\hat g_u^\pm>\hat g_v^\pm>\hat g_u^\pm/D$. The condition $D\hat g_u^\pm>\hat g_v^\pm$ means that the continuous spectrum is unstable to a Turing instability and one can readily show that this stationary point corresponds to a local maximum of $\Re(\lambda)$ that is positive (see Figure~\ref{fig:FredholmBoundaries}b and d).

In the above paragraph we have shown that $\Delta\ge0$ if $\hat g_v^\pm\hat g_u^\pm\le0$ or if both $\hat g_v^\pm<0$ and $\hat g_u^\pm<0$. However, it is possible that $\Delta$ may be negative if $\hat g_v^\pm>0$ and $\hat g_u^\pm>0$ . In this case we can factor 
$\Delta=\left[(\sqrt{\hat g_v^\pm}+\sqrt{\hat g_u^\pm})^2-(1-D)k^2\right]
              \left[(\sqrt{\hat g_v^\pm}-\sqrt{\hat g_u^\pm})^2-(1-D)k^2\right]$.
So, $\Delta<0$ if 
\begin{equation}
\frac{|\sqrt{\hat g_v^\pm}-\sqrt{\hat g_u^\pm}|}{\sqrt{1-D}}<|k|<\frac{\sqrt{\hat g_v^\pm}+\sqrt{\hat g_u^\pm}}{\sqrt{1-D}},
\label{kRangeDeltaNegative}
\end{equation}
that represents two symmetric ranges of $k$ that are bounded away from $k=0$.
In these two ranges of $k$ we get 
\begin{equation}
\Im(\lambda)=\frac{2kc\pm\sqrt{-\Delta}}{2}\quad\text{and}\quad\Re(\lambda)=\frac{-(D+1)k^2-\hat g_v^\pm+\hat g_u^\pm}{2}.
\end{equation}
If $c=0$, we obtain two parts of the essential spectrum. The first part comes from the range of $k$ for which $\Delta\ge0$ and is given by 
$\lambda\in\left(-\infty,\max\{0,\hat g_u^\pm- \hat g_v^\pm\}\right]$. This is the same as the case in which 
$\hat g_u^\pm\le0$ or $\hat g_v^\pm\le0$. 
The second part comes from the range of $k$ for which $\Delta<0$ and gives an ellipse in the complex plane given by
\begin{equation}
\left[ \Re(\lambda) -g_v^\pm-D\hat g_u^\pm \right]^2+(D+1)^2\Im(\lambda)^2=g_u^\pm\hat g_v^\pm.
\label{c0ellipse}
\end{equation}
This means that if $\hat g_u^\pm>0$ and $\hat g_v^\pm>0$, then the essential spectrum contains temporal oscillatory modes for stationary waves that do not occur if $\hat g_u^\pm\le0$ or $\hat g_v^\pm\le0$.

For $c\ne0$ and $\hat g_v^\pm>\hat g_u^\pm>0$ we plot the typical behaviour of the Fredholm borders in Figure~\ref{fig:FredholmBoundaries}e. The Fredholm borders visually appear as five intersecting curves. There is a closed curve that passes through the origin and $(-\hat g_v^\pm+\hat g_u^\pm,0)$. This curve corresponds to 
$|k|<|\sqrt{\hat g_v^\pm}-\sqrt{\hat g_u^\pm}|/\sqrt{1-D}$ for which $\Delta>0$. There are two symmetric closed curves that intersect the first closed curve perpendicularly and correspond to the range of $k$ given in \eqref{kRangeDeltaNegative} for which $\Delta<0$. Finally, there are two symmetric open curves in the left part of the complex plane that correspond to values of $|k|>(\sqrt{\hat g_v^\pm}+\sqrt{\hat g_u^\pm})/\sqrt{1-D}$ for which $\Delta>0$. 
The case for $c\ne0$ and $\hat g_u^\pm>\hat g_v^\pm>0$ is shown in Figure~\ref{fig:FredholmBoundaries}. The basic structure of the Fredholm borders is similar but the right-most closed curve extends into the right half of the complex plane.

Having determined the structure of the Fredholm borders in the different cases, we now turn our attention to the number of roots of \eqref{Characteristic} with $\Re(\nu)>0$ for different values of 
$\lambda$. Clearly, as $\lambda$ varies, number of roots of \eqref{Characteristic} with $\Re(\nu)>0$ can only change at values of $\nu$ that are purely imaginary. This corresponds to values of $\lambda$ that are part of the Fredholm borders that we analysed above. The Fredholm borders will divide the complex plane into different regions and we simply need to determine the number of roots with $\Re(\nu)>0$ in each of the regions.
To achieve this, we begin by considering large values of $|\lambda |$.

\black{On the one hand, if $D>0$} the asymptotic form of the four roots of \eqref{Characteristic} is independent of  
$\hat g_u^\pm$ and $\hat g_v^\pm$ and is given by
\begin{equation}
\nu=\pm\lambda^{1/2}-\frac{c}{2}+O(\lambda^{-1/2})
\qquad\text{and}\qquad
\nu=\pm\frac{\lambda^{1/2}}{D^{1/2}}-\frac{c}{2D}+O(\lambda^{-1/2}).
\end{equation}
So, if $\arg(\lambda)\ne\pi$, then we will always get two solutions for $\nu$ with positive real part and two 
with negative real part. \black{Alternatively,} if $\arg(\lambda)=\pi$, then we will get four solutions with positive real part and zero solutions with negative real part if $c<0$, and four solutions with negative real part and zero solutions with positive real part if $c>0$.

\black{On the other hand, if $D=0$ and $c\ne0$, \eqref{Characteristic} has only three roots and their asymptotic form is given by
\begin{equation}
\nu=\pm\lambda^{1/2}-\frac{c}{2}+O(\lambda^{-1/2})
\qquad\text{and}\qquad
\nu=\frac{\lambda}{c}+\frac{\hat g_v^\pm}{c}+O(\lambda^{-1}).
\end{equation}
So, if $c>0$ and $\Re(\lambda)>0$ we will have two solutions for $\nu$ with positive real part and one 
with negative real part. If $c>0$ and $\Re({\lambda})<0$ with $\arg(\lambda)\ne\pi$, then we have two solutions for $\nu$ with negative real part and one 
with positive real part. If $c>0$ and $\arg(\lambda)\ne\pi$ all three solutions will have negative real part. The case of $c<0$ follows from the symmetry $c\to-c$, $\nu\to -\nu$.
In the case $D=0$ and $c=0$, \eqref{Characteristic} has only two roots and their asymptotic form is given by $\nu=\pm\lambda^{1/2}+O(\lambda^{-1/2})$.}

A particularly simple case is for $D=1$. Then \eqref{Characteristic} can be factored and we obtain
\begin{equation}
\nu^2+c\nu-\lambda=0\quad\text{or}\quad\nu^2+c\nu-\lambda-(\hat g_v^\pm- \hat g_u^\pm)=0.
\end{equation}
 The Fredholm borders in the case $D=1$ are given by two parabolas in the complex plane 
$\lambda=-k^2+ikc$ and $\lambda=-k^2+ikc-(\hat g_v^\pm- \hat g_u^\pm)$.
In light of the behaviour for large $|\lambda |$ we see that if $c>0$, the regions in the complex plane to the right, in between, and to the left of the two parabolas will have two, one, and zero values of $\nu$ with positive real part, respectively.
For $c<0$, the symmetry mentioned above means that we will have two, three, and four values of $\nu$ with positive real part. For $c=0$, the parabolas collapse onto the part of the real axis for which 
$\lambda\in\left(-\infty,\max\{0,\hat g_u^\pm- \hat g_v^\pm\}\right]$. On this line the roots for $\nu$ will be purely imaginary and in the rest of the complex plane there are two values of $\nu$ with positive real part.

For $D<1$ and $\hat g_v^\pm\le0$ or $\hat g_u^\pm\le0$ (see Figures~\ref{fig:FredholmBoundaries}a--d) the two curves divide the complex plane into three regions in which the right side of the complex plane has two positive roots, the region between the two curves has one positive root and the region to the left has no positive roots.
 For $D<1$ and $\hat g_v^\pm>0$ and $\hat g_u^\pm>0$ (see Figures~\ref{fig:FredholmBoundaries}e--f) the situation is slightly more complicated. The region to the right of all of the curves has two positive roots. Regions enclosed by only one of the curves have one positive root. Regions enclosed by more than one of the closed curves have two positive roots. Regions enclosed by one of the closed curves and the open curves have no positive roots.
 
\begin{figure}
	\centering
	\renewcommand{\arraystretch}{0.1} % it makes the space between two rows
	%\begin{tabularx}{\textwidth} {C{6cm}*4{Z}}  % C{6cm} this means the gap between columns  4{Z} means it creates column               
		%\setlength{\tabcolsep}{0pt}
		\includegraphics[scale=0.4]{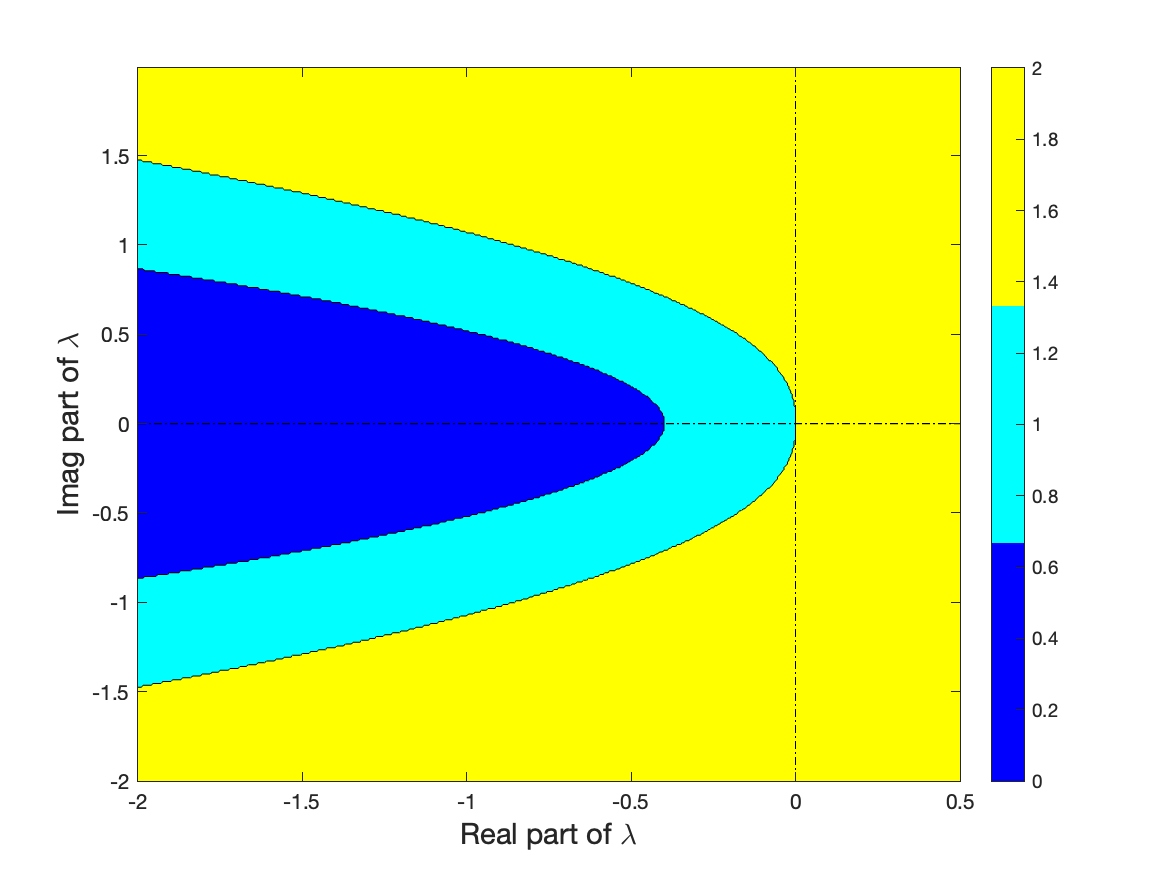}\\
		\includegraphics[scale=0.4]{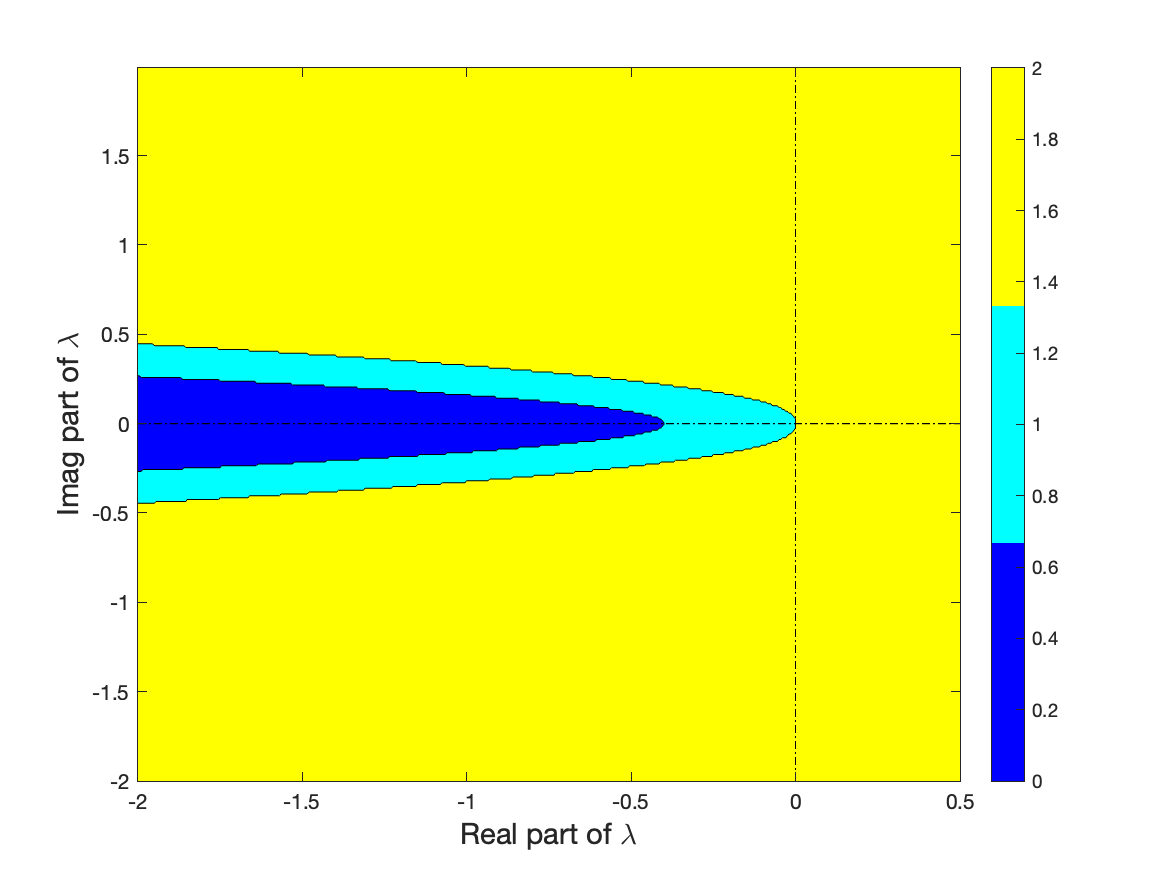}\\
		\includegraphics[scale=0.4]{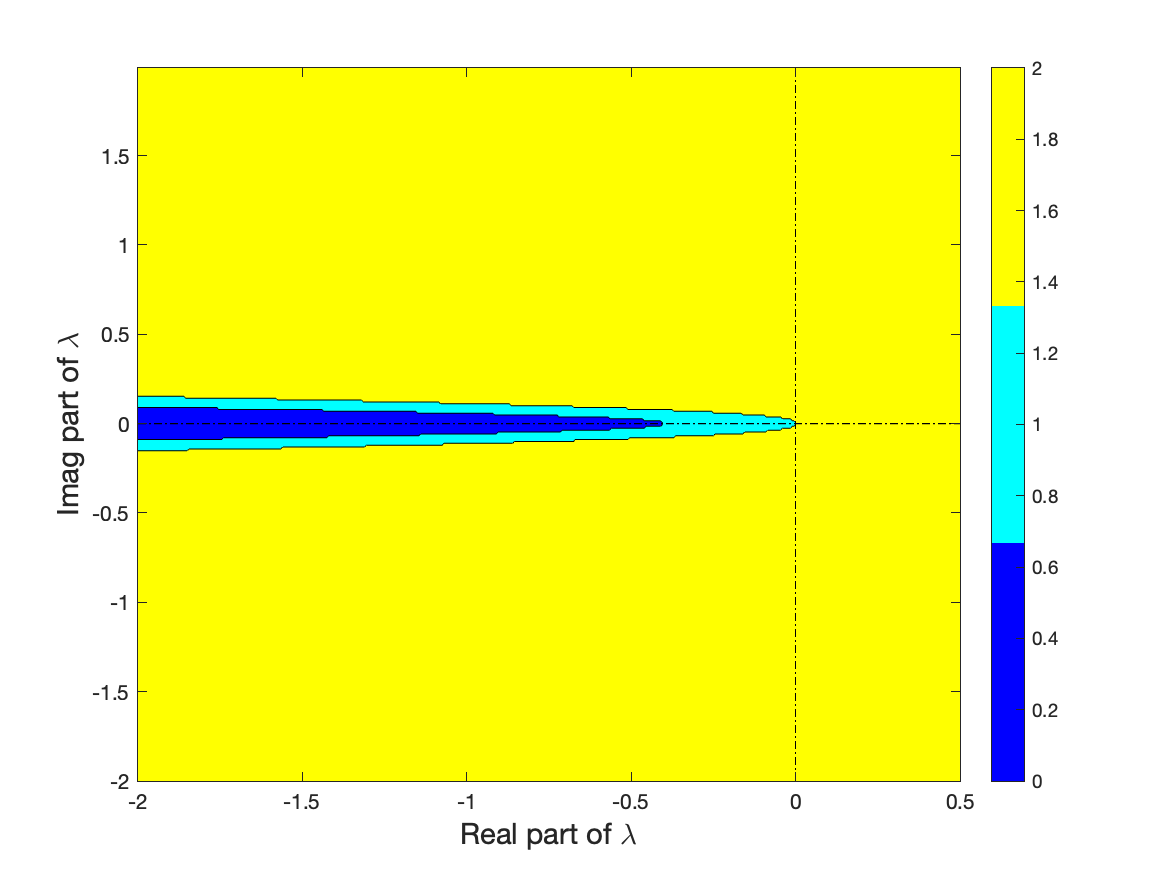}\\
		
	%\end{tabularx}
	\caption{The behaviour as $c\to0$ of the Fredholm boundaries and the number of roots with positive real part for $\nu$ are shown for a typical case in which $\hat g_u^\pm\le0$ or $\hat g_v^\pm\le0$. 
	The value of $c$ decreases from the upper subfigure to the lower subfigure. As $c\to0$ the Fredholm borders collapse onto a portion of the real axis.} 
	\label{fig:FredholmBoundariesCZeroReal}
\end{figure}

\begin{figure}
	\begin{center}
	\renewcommand{\arraystretch}{0.1} % it makes the space between two rows
	%\begin{tabularx}{\textwidth} {C{6cm}*4{Z}}  % C{6cm} this means the gap between columns  4{Z} means it creates column               
		%\setlength{\tabcolsep}{0pt}
		\includegraphics[scale=0.35]{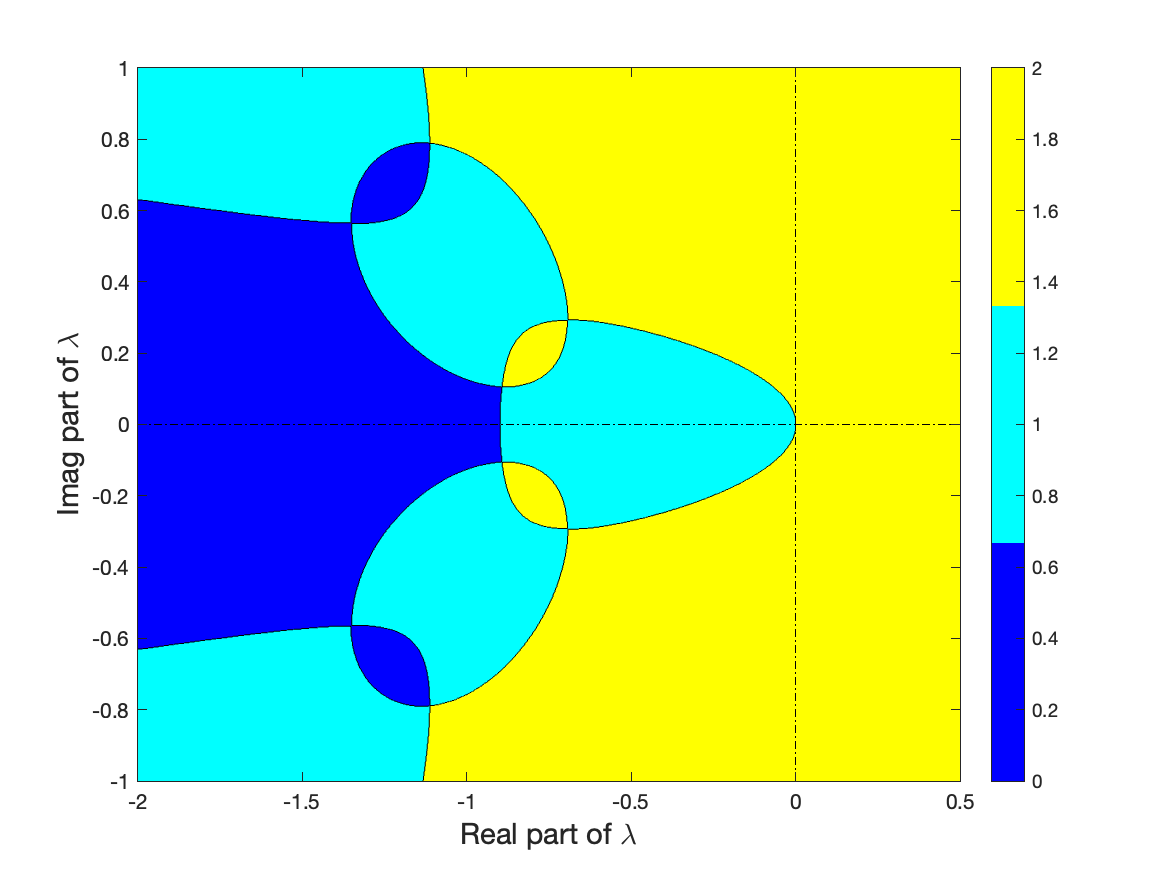}\\
		\includegraphics[scale=0.35]{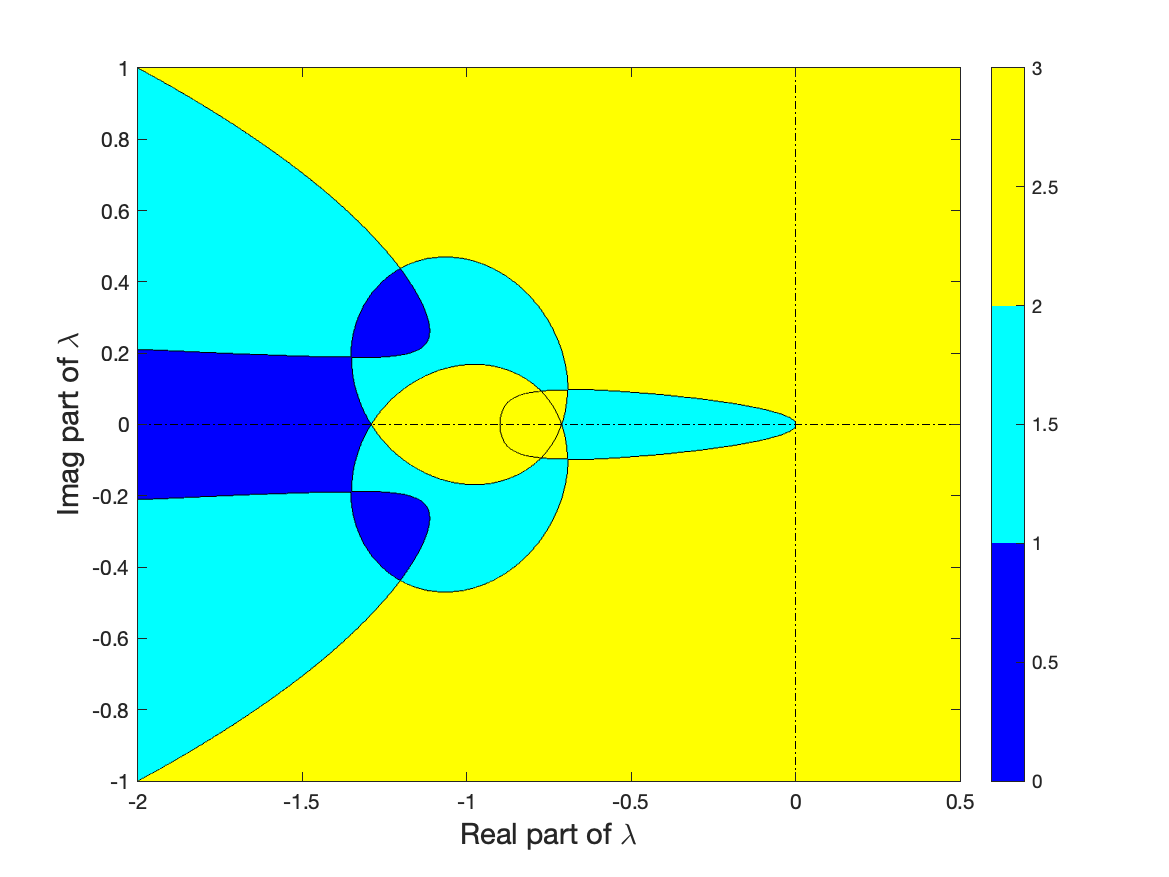}\\
		\includegraphics[scale=0.35]{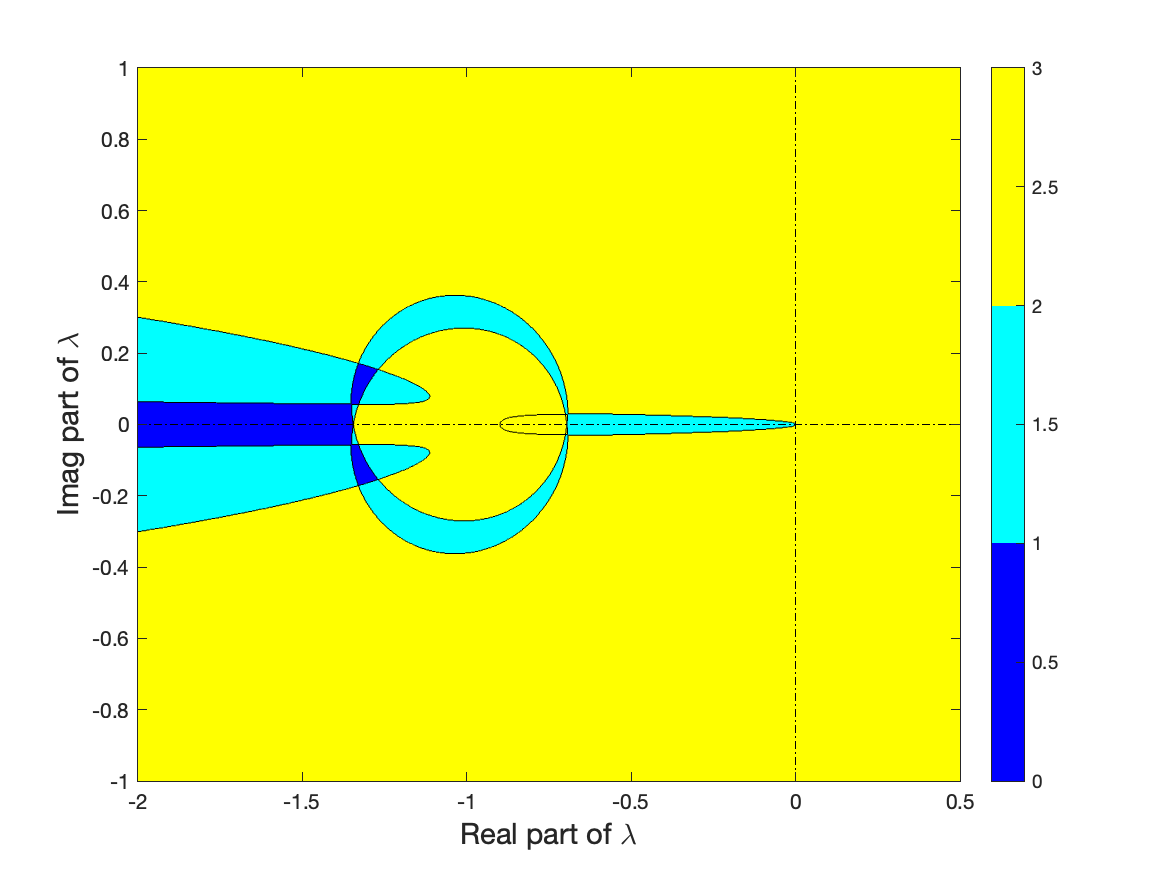}\\
		\includegraphics[scale=0.35]{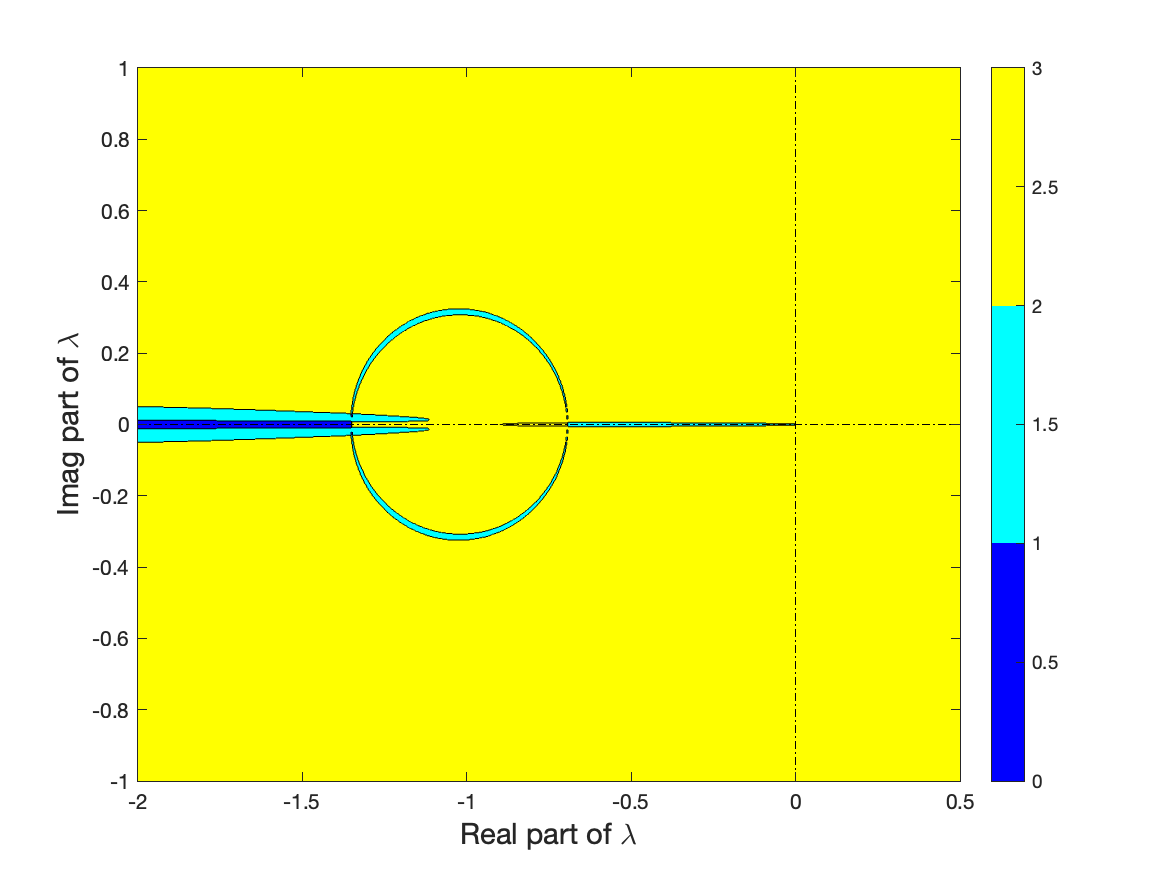}\\
	%\end{tabularx}
	\end{center}
	\caption{The behaviour as $c\to0$ of the Fredholm boundaries and the number of roots with positive real part for $\nu$ are shown for a typical case in which $\hat g_u^\pm>0$ and $\hat g_v^\pm\>0$. 
	The value of $c$ decreases from the upper subfigure to the lower subfigure. As $c\to0$ the Fredholm borders collapse onto a portion of the real axis and an ellipse.} 
	\label{fig:FredholmBoundariesCzeroComplex}
\end{figure}

Above we showed that for $c=0$, the Fredholm borders are particularly simple and always contain the part of the real axis $\lambda\in\left(-\infty,\max\{0,\hat g_u^\pm- \hat g_v^\pm\}\right]$. However, if $\hat g_v^\pm>0$ and $\hat g_u^\pm>0$, the Fredholm border also contains the ellipse \eqref{c0ellipse}.
We now consider how the number of positive roots for $\nu$ behaves in the limit as $c\to0$. 
If $\hat g_v^\pm\le0$ or $\hat g_u^\pm\le0$ the behaviour as $c\to0$ is straightforward and is shown in Figure~\ref{fig:FredholmBoundariesCZeroReal}. As $c\to0$ the two curves simply collapse onto $\lambda\in\left(-\infty,\max\{0,\hat g_u^\pm- \hat g_v^\pm\}\right]$ and so there will be two positive roots for the rest of the complex plane.
For $\hat g_v^\pm>0$ and $\hat g_u^\pm>0$ the situation is more complicated and is shown in Figure~\ref{fig:FredholmBoundariesCzeroComplex}.
As $c\to0$ the rightmost closed curve and the two leftmost open curves collapse onto the real axis while the other two closed curves both collapse onto the ellipse \eqref{c0ellipse}. The net result of this is that all of the complex plane except for the continuous spectrum will have two positive roots.

In order to determine if a given value of $\lambda$ is of part of the essential spectrum we need to consider its Fredholm index. This is determined by the difference between the number of roots of \eqref{Characteristic}, $\nu$ with positive real part for $z\to\infty$ and $z\to-\infty$. For pulses, the two end states are identical and so the away from the Fredholm borders the Fredholm index is necessarily zero. This implies that the essential spectrum for pulses consists only of these curves. This is a generic feature of all reaction-diffusion systems and is not particular for the degenerate-type systems that we are considering here.

\begin{figure}
	\centering
	\includegraphics[scale=0.5]{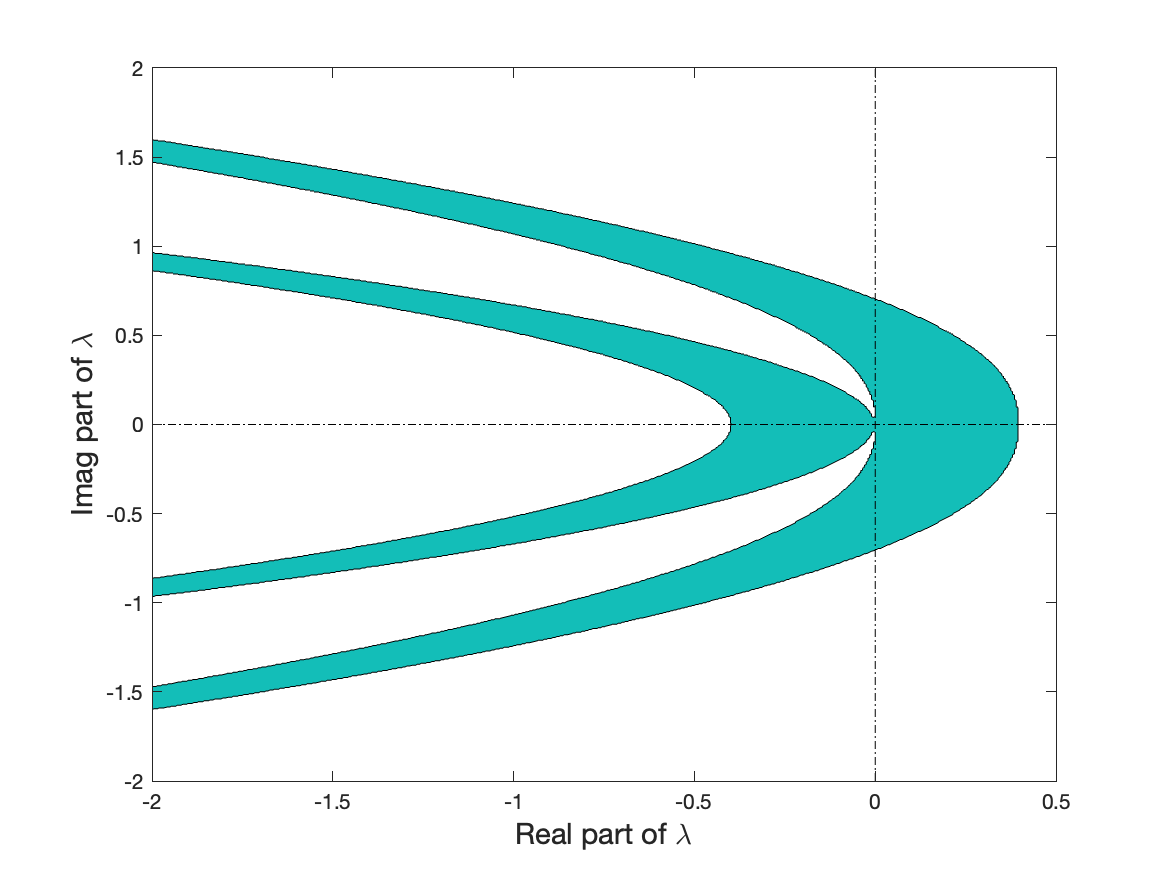}
	\caption{The essential spectrum for a travelling front solution whose asymptotic states far behind and far in front of the wave correspond to Figures~\ref{fig:FredholmBoundaries}a and c .} 
	\label{fig:Essential}
\end{figure}

We now consider the structure of the essential spectrum for travelling fronts. For $c=0$, we have shown that there will be there will be two roots with positive real part for $\nu$ for all values of $\lambda$ away from the Fredholm borders. This result is independent of the form of $g(u,v)$.
The consequence of this is that the Fredholm index is necessarily zero and the essential spectrum for stationary pulses consists only of the Fredholm borders.
However, for travelling front solutions the behaviour is more complicated and if $\hat g_v^+\ne\hat g_v^-$ or $\hat g_u^+\ne\hat g_u^+$ there will generically be portions of the complex plane in which the Fredholm index will be non-zero and so the essential spectrum will not only consist of the essential spectrum. This is illustrated in Figure~\ref{fig:Essential} where we consider a case in which a travelling wave connects two asymptotic states. The asymptotic states far behind and far in front of the travelling wave are taken to be Figures~\ref{fig:FredholmBoundaries}a and c. The plot shows the regions of the complex plane for which the Fredholm index is not zero and is hence contained in the essential spectrum.

\section{Analytical results}

A number of special cases arise in system \eqref{eq:TWsys2}. Looking at the second equation, we see that if $D=0$, the second derivative no longer appears, reducing the system from fourth order to third order. Consequently, the case when the $v$ does not diffuse must be treated separately.
The system is doubly singular if $D=0$ and the wavespeed, $c=0$, so that the problem reduces to a second order alegebro-differential equation.

In order to investigate the stability of solutions to models of type \eqref{eq:TWsys2} we now have three cases to consider: the case of nonstationary or stationary waves (fronts and pulses) when $v$ is diffusing, that is, $D\ne0$; the case of nonstationary waves when $v$ is not diffusing, that is $D=0,~c\ne0$; and the case of stationary waves when $v$ is not diffusing, that is $D=c=0$. There is a fourth, less obvious case that must be carefully considered. As explained in \cite{BHW19}, the case of equal diffusivities, $D=1$, can also be reduced to a second order problem (see below for more detail). Example solutions in these four cases were outlined in \cite{BHW19}. Straightforward examination of the linearised equations \eqref{eq:eigenvalue1} reveals that the stability problem can be separated into the same four cases.

For the cases $D=1$ and $D=c=0$ we present some analytical results below. The final two cases ($D=0,~c\ne0$ and $D\ne0,1$) must be analysed numerically.

\subsection{Equal diffusivity, $D=1$}

By adding equations \eqref{eq:TWsys1} with $D=1$, we find that the quantity $w=u+v$ must obey
\begin{equation*}
w_{\tau}=w_{zz}+cw_z,
\end{equation*}
while the second equation becomes

\begin{equation}\label{eq:fourform}
v_{\tau} = v_{zz} + c v_z -  g(w-v,v).
\end{equation}

Equation \eqref{eq:fourform} is impossible to solve in general, unless we have additional information about $w$. Since we are interested in travelling waves, set $w_{\tau} =0$. 
In this case, without loss of generality, we assume that one of the far-field rest states is (u,v)=(0,0) and follow \cite{WM06,BHW19} to integrate the first equation and apply the boundary conditions to obtain
\begin{equation}\label{eq:integrated}
w_z + c w = 0,
\end{equation}
whence we get that $w = 0$ and so $u = - v$. The second equation then becomes
$$
v_{zz} + c v_{z} - g(-v,v) = 0.
$$
Several example reaction terms, $g(-v,v)$, and solutions were given in \cite{WM06,BHW19}.  
In particular, if we have a function $g$ that produces a pair of solutions $\wh{u},\wh{v} = - \wh{u}$ the linearised equations \eqref{eq:eigenvalue1} can be written as
\begin{equation}\label{eq:linearisedcombo}
\begin{aligned}
\lambda n &  = n_{zz} + c n_z\,, \\ 
\lambda q & = q_{zz} + c q_z - g_u(\wh{u},\wh{v}) (n-q) - g_v(\wh{u},\wh{v})q \,,
\end{aligned}
\end{equation}
where $n = p + q$. The key observation here is that if we choose $\lambda$ away from the line $-k^2 + ick$ (in particular if $\lambda$ is in the right half plane), then in order for $n \in L^2(\R)$ we must have $n = 0$. The second equation in \cref{eq:linearisedcombo} then reduces to the linearisation of a nonlinear equation. Letting $G(v) := g(-v,v)$, if $n =0$, the second equation in  \eqref{eq:linearisedcombo} is then 
$$
\lambda q=q_{zz} + c q_z - G'(\wh{v})q
$$
which is the eigenvalue problem for the linearisation of the nonlinear planar ODE
\begin{equation}\label{eq:master}
v_{zz} + c v_z - G(v) = 0.
\end{equation}
We can then use phase plane analysis to describe qualitative results about solutions as well as their stability. For example, homoclinic orbits are only possible if $c = 0$. To see this we multiply \eqref{eq:master} by $v_z$ and integrate with respect to $z$. Setting  $H(v) = \int_0^v G(v) \, dv$, we have 
\begin{equation}
\begin{aligned}
&\int \frac{\partial}{\partial z} \frac{v_z^2}{2} \, dz + \frac{c}{2} \int v^2 \, dz - \int \frac{\partial}{\partial z} H(v) \, dz \\
&\qquad\qquad = \frac{v_z^2}{2} \bigg{|}_{z=-\infty}^{+\infty} + \frac{c}{2}\int v^2 \, dz - H(v)\bigg{|}_{z=-\infty}^{+\infty} = 0
\end{aligned}
\end{equation}
For homoclinic orbits of \eqref{eq:master}, the boundary conditions cause the integrated terms vanish and we have 
\begin{equation}
\frac{c}{2} \int_{\R} v^2 \, dz = 0
\label{D1PulseStationary}
\end{equation}
whence we see that $c \equiv 0$. That is, no {\it nonstationary} pulses exist in the case of equal diffusivities.

To analyse stability, we know that $q = \wh{v}'(z)$ will solve
\begin{equation}
q_{zz} + c q_z - G'(\wh{v})q = 0.
\label{DvSatisfiesLambda0}
\end{equation}
So, if we have $\wh{v} \in H^1(\R)$ for instance, then $\lambda = 0$ will be an eigenvalue of the linearised problem with eigenfunction $ \wh{v}'(z)$. Thus, by Sturm-Liouville theory \cite{KP13}, the number of zeros of $\wh{v}'(z)$ indicates the number of eigenvalues greater than or equal to $0$. The consequence of this is that stationary pulses will be unstable, while fronts (stationary or otherwise) will be stable if and only if they are monotone \cite{KP13}. In certain special cases, we can even compute the eigenvalues and we provide two examples below. 

\subsubsection*{Example 1: Stationary pulse.}\label{sec:D1pulse}

If $g(u,v)=-u(4-6v)$ in equations \eqref{eq:sys1}, \cite{BHW19} found solutions
\begin{align}\label{eq:Eg1sol}
    \wh{u}=-D\sech^2(z + \beta)=-\sech^2(z + \beta), \qquad \wh{v}=\sech^2(z + \beta)
\end{align}
for a fixed shift $\beta$ (below we set $\beta=0$, but the analysis is the same for $\beta\ne0$). 
Note that here we will set $D=1$ but these solutions are also valid for $D\ne1$ as we shall see later. The second equation in \eqref{eq:linearisedcombo} becomes
\begin{equation}\label{eq:statpulseD1}
\lambda q = q_{zz} + 12 \sech^2(z) q - 4q\,.
\end{equation}
In the far field this equation becomes $\lambda q=q_{zz}-4q$ since $\sech^2z\to0$ as $z\to\pm\infty$. The continuous spectrum from the $n$ and $q$ equations respectively is $(-\infty,-4]\cup(-\infty,0]=(-\infty,0]$. There is no point spectrum coming from the first equation, and indeed, again, any element of the point spectrum of $\cL$ must have $n \equiv 0$. The point spectrum of equation \eqref{eq:statpulseD1} can be obtained using ladder operator techniques and is given by $\lambda = -3, 0,  5$ (it is the spectrum of the $\sech^2(z)$ profile shifted by 4; see for example \cite{CDB09}), and since $\lambda = 0$ and $\lambda = -3$ will be in the continuous spectrum, we have a single eigenvalue of the linearised operator at $\lambda = 5$. Since we have one eigenvalue in the right half plane, this stationary pulse is unstable. The Evans function in this case can be computed explicitly as the Jost solutions \cite{yafaev1992mathematical} are known in each of the cases and we find\footnote{Traditionally $D(\lambda)$ is used to denote the Evans function, but as our diffusion parameter is called $D$, we use this alternate notation.}
$$
F(\lambda) = 4(\lambda-5)\lambda^{3/2}(\lambda+3)\sqrt{\lambda + 4}\,.
$$
We note that, as far as we are aware, there are very few examples for which explicit, exact Evans functions are known. Such exact expressions can be invaluable for validating numerical methods.

\subsubsection*{Example 2: Stationary front.} In the case when $g(u,v)=2u(1-v^2)$, the following stationary front solutions were found by \cite{BHW19}: $\wh{u}=-D\tanh z=-\tanh z$ and $\wh{v}=\tanh z$ (here we set $D=1$ but this solution is valid for $D\ne1$). The second linearised equation in \eqref{eq:linearisedcombo} then becomes
\[
\lambda q = q_{zz} + 6 \sech^2(z) q - 4q\,.
\]
The first equation in \eqref{eq:linearisedcombo} contributes the half line $(-\infty,0]$ to the essential spectrum. Again, there is no point spectrum coming from the first equation. The point spectrum of the second equation is the spectrum $\{ 1 ,4 \}$ shifted by $4$, or $\{ -3, 0 \}$ and so our wave for $D = 1$ is spectrally, linearly and nonlinearly stable.

\subsection{Stationary waves, $v$ does not diffuse, $D=c=0$}

In the case of stationary waves when $v$ does not diffuse, that is $D=0$, we see that the second equation in \eqref{eq:TWsys1} requires that $g(\hat u(z),\hat v(z))=0$. Substituting this into the first equation of \eqref{eq:TWsys1} with $c=0$, we see that $\hat u_{zz}=0$, so that the only solution to satisfy the far field boundary conditions is if  $u(z)$ is a constant function. \black{Applying the left far-field boundary condition} 
we obtain $\hat u(z)\equiv u_-$. We then require that $g(u_-,\hat v(z))=0$, so that any constant $v_*$ that satisfies $g(u_-,v_*)=0$ is a solution. Generically, there will be a finite set of such values of $v_*$. 
We note that the set must be non-empty because $v_*=v_-$ is a solution. If there is more than one value of $v_*$, then a nontrivial solution of \eqref{eq:TWsys1} can  be constructed by setting $\hat u\equiv u_-$ and choosing $\hat v$ to be any piecewise constant solution in which each of the piecewise constants is chosen from the set of $v_*$. Clearly, both pulse and front solutions can be constructed by choosing $v_-=v_+$ and $v_-\ne v_+$ respectively.

Our linearised operator reduces to a second order algebro-differential equation
\begin{equation} \label{eq:second}
\begin{aligned}
\lambda p & = p'' + g_u p + g_v q\,, \\ 
\lambda q & = -g_u p - g_v q\,.
\end{aligned}
\end{equation}
Here, $g_u=g_u(u_-,\hat v(z))$ and $g_v=g_v(u_-,\hat v(z))$ are given by piecewise constant functions since $\hat v(z)$ is a piecewise constant.

We note that if $\lambda=-g_v$ then from \eqref{eq:second}(b) we immediately see that $p\equiv0$ on the corresponding piecewise region. Then from \eqref{eq:second}(a) we also see that $q\equiv0$ on that region as long as $g_v\ne0$. However, $p$ and $p'$ must be continuous at the edges of the piecewise regions and so $p$ must also be identically zero on the neighbouring regions and hence on the entire domain. The case in which $g_v=\lambda=0$ also cannot yield eigenfunctions. Therefore, we conclude that $\lambda\ne-g_v$. We can therefore make the substitution
$$
q = - \frac{g_u}{g_v +\lambda} p, 
$$
and rewrite \cref{eq:second} as 
\begin{equation}\label{eq:herg}
p '' + \left(g_u  - \frac{g_u g_v}{\lambda+g_v} - \lambda \right) p = 0.
 \end{equation}
 We now suppose that $p(z)$ is an eigenfunction and multiply \eqref{eq:herg} by the complex conjugate $p^*$ and integrate from $z=-\infty$ to $z=\infty$ to obtain
\begin{equation}\label{eq:herg_cc}
\int_{-\infty}^\infty p^*p '' dz+ \int_{-\infty}^\infty\left(g_u  - \frac{g_u g_v}{\lambda+g_v} - \lambda \right) p^*p dz = 0.
 \end{equation} 
 Integrating by parts and applying the boundary conditions at infinity we obtain
 \begin{equation}\label{eq:herg_int}
-\int_{-\infty}^\infty |p'|^2 dz+ \int_{-\infty}^\infty\left(g_u  - \frac{g_u g_v}{\lambda+g_v} - \lambda \right) |p|^2 dz = 0.
 \end{equation} 
 Taking the imaginary component and recalling that $g_u$ and $g_v$ are piecewise constants we obtain
\begin{equation}\label{eq:herg_imag}
\Im(\lambda)\sum_k\left(\frac{g_u g_v}{\Im(\lambda)^2+[g_v+\Re(\lambda)]^2} \right) \int_{I_k}|p|^2 dz = 0,
 \end{equation}  
 where the sum is over all of the piecewise regions $I_k$ that compose the stationary wave solution.
 We immediately see that if either $g_u g_v<0$ or $g_u g_v>0$ on every one of the piecewise regions then $\Im(\lambda)=0$ and the eigenvalues must be real.
 We note that a similar result was obtained by \cite{BBR19}.
 
Likewise, the {\em Fredholm borders} of the essential spectrum of the pencil will be the $\lambda$ such that the far-field matrices 
 $$
 A_\pm(\lambda) = \lim_{z \to \pm \infty} \begin{pmatrix} 0 & 1 \\ \left(\lambda - \frac{(-g_u g_v)}{\lambda+g_v} \right) - g_u & 0 \end{pmatrix}
 $$ 
 have a purely imaginary eigenvalue. These will be the borders of the set where the asymptotic matrices will have distinct numbers of eigenvalues with positive real parts, i.e, where the {\em Fredholm index} of the pencil is non-zero. We will again refer to this set as the {\em essential spectrum} of the operator (and of the pencil).

\subsubsection*{Example 3: Piecewise constant stationary pulse.} As described above, piecewise constant solutions exist in the case when $D=c=0$. To illustrate a piecewise constant standing pulse, we choose \cite{BHW19}
\[
g(u,v) = v(v-\gamma)(v-1) - u =: \tilde{g}(v) - u, \quad \textrm{where} \quad 0<\gamma<1.
\]
In order to make our analysis explicit we begin by choosing a $\wh{v}$ which is a `standing pulse' in the sense that it is $0$ outside of some interval $[-L,L]$ and $1$ inside it. That is,
\[
\wh{v} = \begin{cases} 0 &  \text{for } |x| > L \\  1 & \text{for } |x|<L. \end{cases}
\]
Noting that $g_u = -1$ and $g_v = \tilde{g}'(\wh{v})$, our linearised operator is the matrix pencil
\begin{equation}\label{eq:herglotz}
p'' - p + \left(\frac{\tilde{g}'(\wh{v})}{\lambda + \tilde{g}'(\wh{v})} - \lambda \right) p  = 0.
\end{equation}
We can then follow a modified proof of the first part of Theorem 3.1 from \cite{BBR19} to show that if $\lambda$ is such that there is an $\Hbb^1$ solution to \cref{eq:herglotz}, then it must be real. Defining the operator pencil from \cref{eq:herglotz} as $M(\lambda)$, suppose that we had a $\lambda$ with a bounded $p \in \Hbb^1(\R)$  such that $M(\lambda)p=0$, then taking the ($L^2$) inner product with $p$ we have that the function $h(\lambda)$ given by 
\[
h(\lambda) := - \langle p, M(\lambda)p\rangle 
\]
is a so-called Herglotz function in the sense that if $\imag{\lambda} > 0$, $\imag{h(\lambda)}>0$. \black{Thus the roots are necessarily real \cite{BBR19,Simon05}}. To see this note that because $\wh{v}$ is piecewise constant, the function $\tilde{g}'(\wh{v})$ is also a real, piecewise constant function. Thus the inner product will be the sum of three Herglotz functions (and hence Herglotz). 

Next we note that the essential spectrum of the linearised operator is also real. Indeed we have that the Fredholm borders for this specific piecewise constant `standing pulse' $\wh{v}$ will consist of the set
\begin{equation}\label{eq:pulseessspec}
\lambda \in (-\infty, -1-\gamma] \cup [-\gamma,0]. 
\end{equation}
Then, for $\lambda \not \in \sigma_{\text{ess}}$ we can compute the characteristic (Evans) function of the operator by matching the decaying solutions at $\pm \infty$ across the part where $\wh{v} =1$. This gives 
\[
F(\lambda) = \sqrt{\lambda}\left[ \sqrt{1 + \frac{1}{\lambda + \gamma}} +\sqrt{1 + \frac{1}{\lambda + 1- \gamma}} \tanh\left(\sqrt{\lambda} L \sqrt{1 + \frac{1}{\lambda + 1- \gamma}}\right)
\right].
\]
Note that for $\lambda = \omega^2 >0$, $F(\lambda)$ has no roots since the argument of the $\tanh$ function in this expression is positive in this case. Consequently, there are no positive eigenvalues $\lambda$ for this specific $(\wh{u},\wh{v})$ pair. We thus have spectral stability in this case. This is in agreement with the solutions shown in \cite{BHW19} where they numerically demonstrate the stability of a similar standing wave.

\subsubsection*{Example 4: Piecewise constant stationary front.} The reaction function from Example 3 can also be used to produce a standing piecewise front connecting 0 to 1,
\[
\wh{v} = \begin{cases} 0 & \text{if } x <0 \\ 1 & \text{if } x >0 \end{cases}
\]
The essential spectrum will be the union of the set in \cref{eq:pulseessspec} with
\begin{equation}\label{eq:frontspecess}
\lambda \in (-\infty, -2+\gamma] \cup [-1+\gamma,0]. 
\end{equation}
In this case, determining the Evans function requires finding a $C^1$ function that matches the values and the derivatives of the unstable solution at $-\infty$ and the stable one at $+ \infty$.
In this instance we have 
$$
F(\lambda) = \lambda \left( \frac{1}{\lambda + 1-\gamma} + \frac{1}{\lambda + \gamma}\right).
$$
As usual, the roots of this function indicate an eigenvalue, and we find
$$
\lambda = 0,~ -1 \pm \frac{\sqrt{2}}{2} \sqrt{\gamma^2 + (\gamma-1)^2},
$$
two of which are negative for $0<\gamma<1$. Consequently, the linearised operator has no spectrum in the right-half plane and we have spectral stability. Note that in this case, for all $\gamma \neq 1/2$ there will be no eigenvalues in the gaps of continuous spectrum, but rather we have a pair of embedded eigenvalues for all $\gamma$. When $\gamma = 1/2$, there is a pair of embedded eigenvalues at the edges of the continuous spectrum.

\begin{remark}
We remark here that these arguments can be applied to investigate stability for $\wh{v}$ with any number of real piecewise constant parts. The Evans functions necessarily become more complicated the more jumps that occur, and depending on the far field values (i.e. $\tilde{g}'(0)$ and $\tilde{g}'(1)$), the essential spectrum may be slightly different. The general form of the Evans function is not amenable to deduce general stability.
\end{remark}

\section{Numerical results}

In the case of nonstationary waves when $v$ does not diffuse, and in the most general case ($D\ne0,1$), we resort to investigating the stability numerically using the Ricatti-Evans function. 

The Riccati-Evans function exploits the fact that we are solving a linear ODE and as such we have a flow on the Grassmannian \cite{HvHMPW15}. We can thus evolve subspaces themselves via a matrix Riccati equation \cite{KvHMRPW2020}.  Then the appropriately chosen far-field subspaces are evolved and compared, and the Evans function determinant is converted to a function on elements in the Grassmannian and evaluated. The Riccati-Evans function has the advantage that the exponential growth of the spanning sets (which is different for different eigenvectors) is tempered by the solutions themselves, but has the disadvantage that typically solutions will blow up at finite values of the independent variable. This makes the Riccati-Evans function a meromorphic one instead of analytic like the Evans function. We adopt the Ricatti-Evans function rather than the `classical' Evans function in this problem because the latter is not numerically robust enough, see Ex. 7 for a comparison of the two.

This blow-up can be avoided by choosing a `chart' on which to evolve the Riccati flow on the Grassmannian. In particular, if we are interested in real eigenvalues, then a complex chart has a good chance of not leading to a singularity for such $\lambda$ (though there still may be other values of $\lambda$'s which are poles of the Riccati-Evans function). The precise details are in \cite{KvHMRPW2020}.

\subsection{Nonstationary waves, $v$ does not diffuse, $D=0,\,c>0$}

When $D=0$, the second equation in \eqref{eq:eigenvalue1} becomes
\[
q' =  \frac{1}{c} \left(g_u p + (\lambda +g_v) q \right)
\]
so that the third-order stability problem can be written as
\begin{equation}\label{eq:stabode3d}
\begin{pmatrix}
p \\ q \\ r 
\end{pmatrix} '  = \begin{pmatrix} - c & 0 & 1 \\ \frac{g_u}{c} & \frac{\lambda+ g_v}{c} & 0 \\ \lambda-g_u & - g_v & 0  \end{pmatrix} 
\begin{pmatrix}
p \\ q \\ r 
\end{pmatrix}=A(\lambda)\,\begin{pmatrix}
p \\ q \\ r 
\end{pmatrix}
\end{equation}
where $r=p'+cp$.

It is still true that in the case when $D=0$, the linearised operator $\cL$ can be written as $\cL_\infty + \cL_1$, with $\cL_\infty$ a constant coefficient operator representing the linearisation in the far-field. 
In this case, we can still investigate the point spectrum of the operator using the Evans function. 
Here we present numerical findings for the spectrum of the linearised operator about two example solutions presented in \cite{BHW19,WM06}.

\subsubsection*{Example 5: Nonstationary pulse.}

The case where the reaction term is given by
\begin{equation}\label{eq:d0trav1pot}
g(u,v) = 8u^3 - 6u^2 + c^2(u+v),
\end{equation}
has nonstationary pulse solutions given by \cite{BHW19}
\begin{equation}\label{eq:d0trav1sols}
\wh{u} = \frac{1}{1+z^2} \quad \text{and } \quad \wh{v} = \frac{-c z^2 + 2 z - c}{c(z^2 +1)^2}.
\end{equation}

We want to investigate the possibility of a solution to \eqref{eq:stabode3d} with \eqref{eq:d0trav1pot} and \eqref{eq:d0trav1sols} in the region where $\lambda \in \C$ is such that $A_\pm(\lambda)$ is hyperbolic. In this case the limiting matrices of \eqref{eq:stabode3d} as $z\to\pm\infty$, $A_\pm(\lambda)$, are equal, so this condition is sufficient to begin our search for an eigenvalue. The region where $A_\pm(\lambda)$ is hyperbolic in particular contains the positive real axis. Direct calculations show that the characteristic polynomial of $A_\pm(\lambda)$ is 
$$
\det (A_\pm(\lambda) - \mu I) = -\mu^3 + \frac{\lambda}{c} \mu^2 + 2 \lambda \mu - \frac{\lambda^2}{c}.
$$
The discriminant of this cubic in $\mu$ is 
$$
\lambda^3\left(32 + \frac{\lambda(13c^2 + 4\lambda)}{c^4}\right)
$$
which is greater than zero when $\lambda >0$ and only vanishes for imaginary $\lambda$. We thus have three distinct real roots. Applying Descartes rule of signs we see that there is one negative root and two positive roots for $\mu$.

\begin{figure}
\includegraphics[scale=0.25]{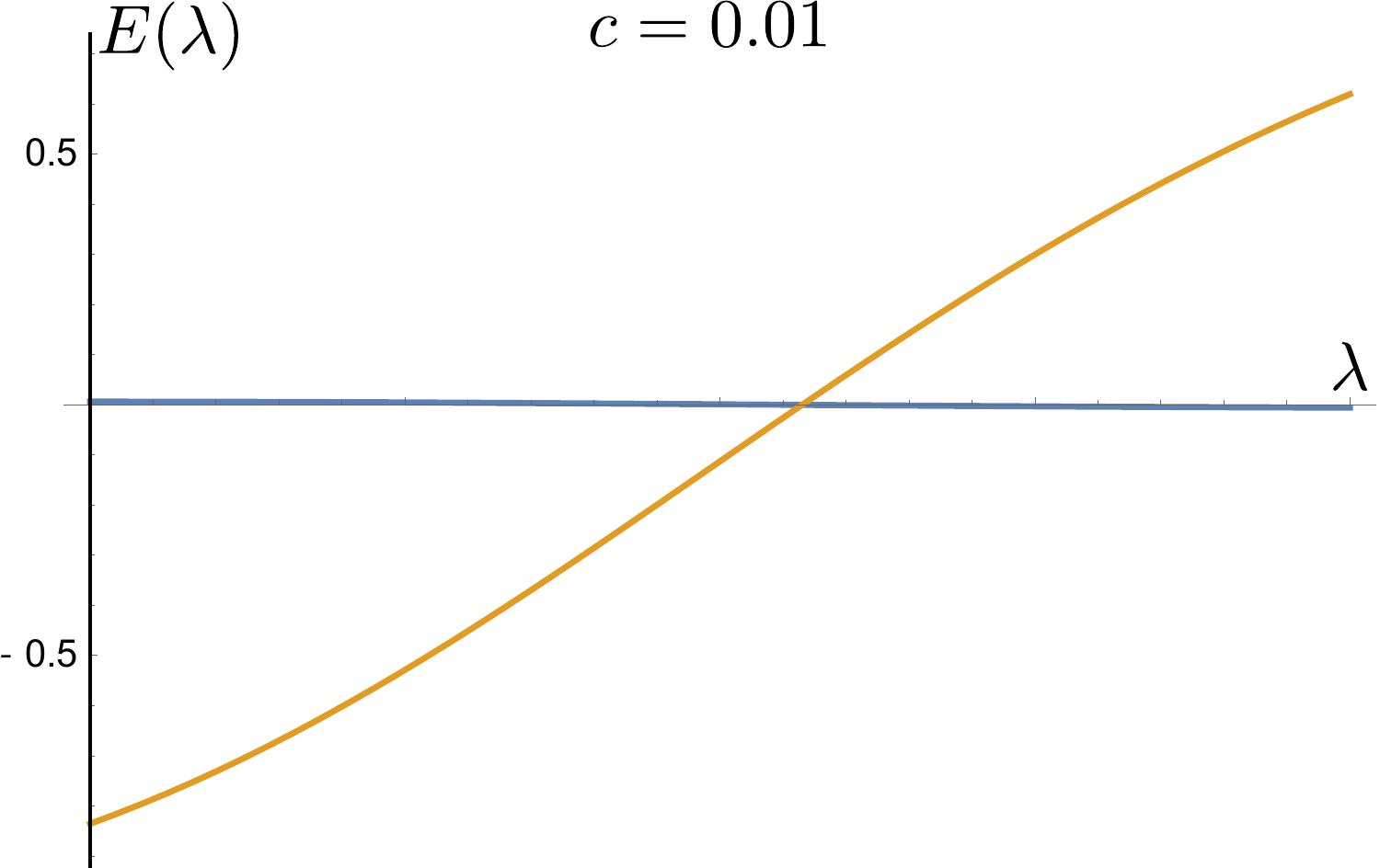}
\includegraphics[scale=0.25]{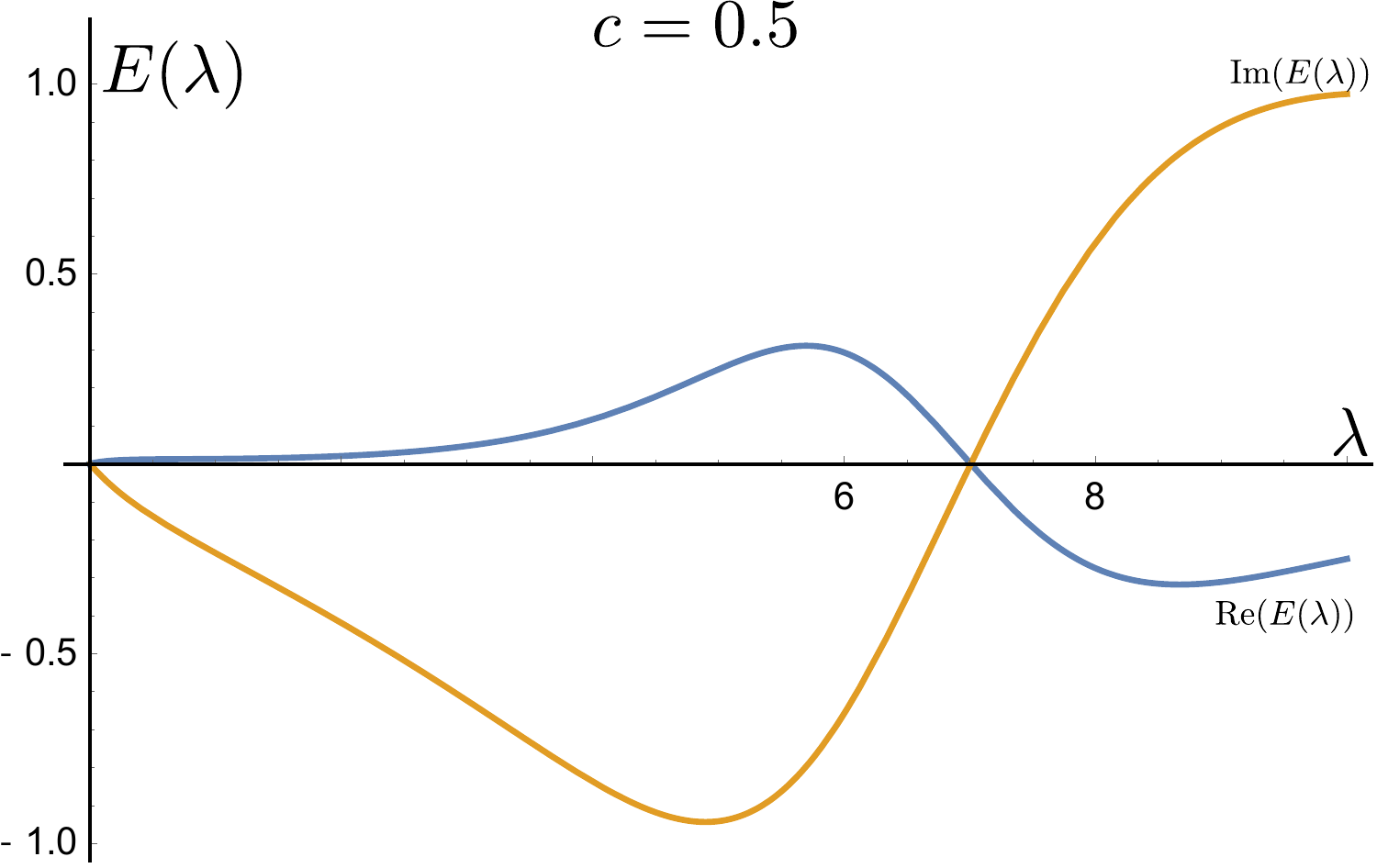} \\ 
\includegraphics[scale=0.25]{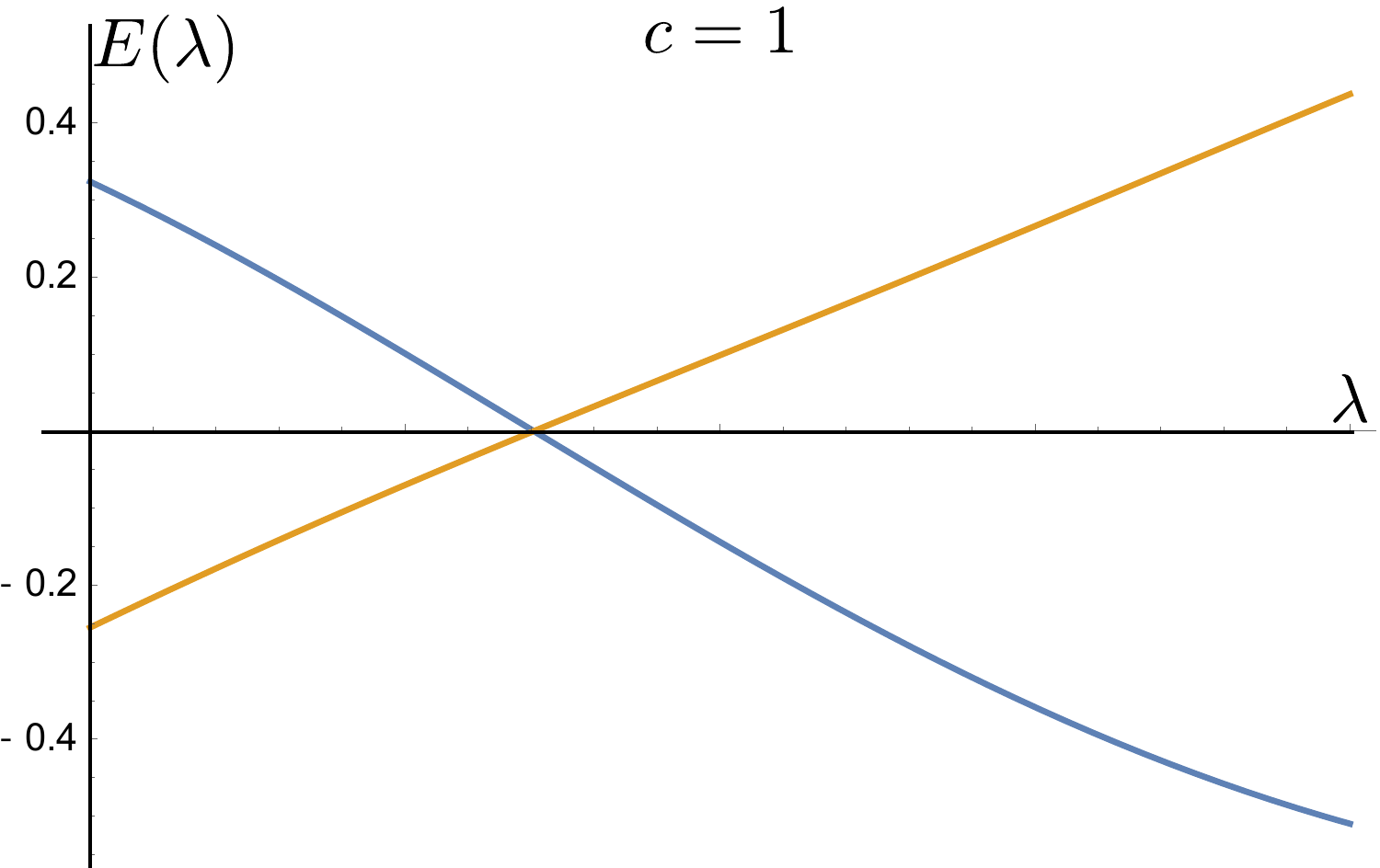}
\includegraphics[scale=0.25]{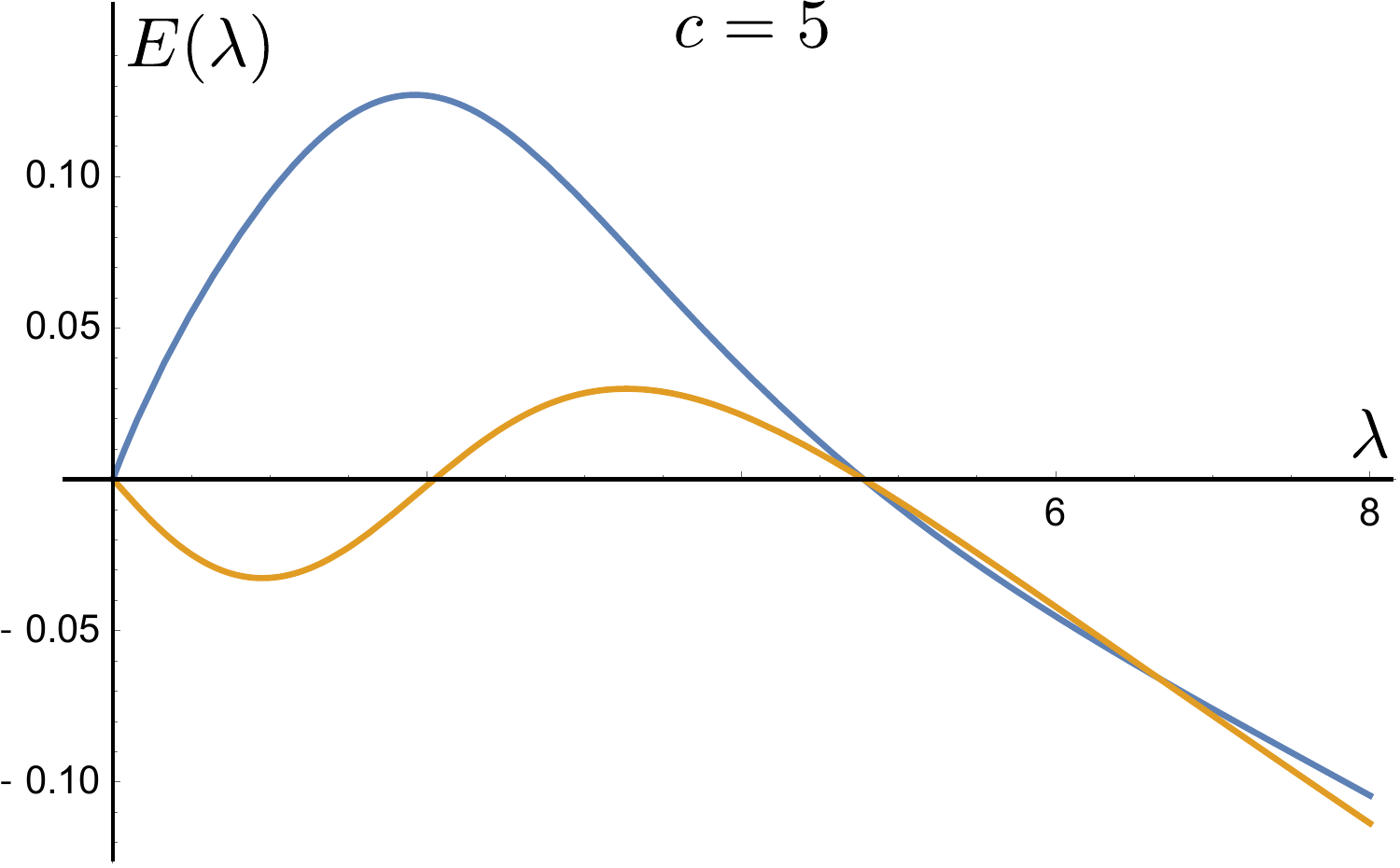}
\caption{Example 5: Nonstationary pulse, $D=0$. Real (blue online) and imaginary (gold) parts of the Riccati-Evans function, \eqref{eq:ricevans3d} for the waves given in \cref{eq:d0trav1sols}  for various values of $c$. In all cases for $c>0$ we have a positive eigenvalue, and hence point spectrum in the right-half plane.}\label{fig:algwave2} 
\end{figure}

The essential spectrum for this operator is not confined to a left-opening cone in the complex plane, so even though one can readily show that the essential spectrum is stable for all parameter values, the standard results concerning nonlinear stability do not apply \cite{KP13}. We now turn our attention to
the point spectrum for this solution. Following \cite{HvHMPW15,KvHMRPW2020} we choose a chart so that the Riccati flow for $\alpha := \frac{q}{p}$ and $\beta := \frac{r}{p}$ will only have imaginary poles.   
We note that this will also hold for the Ricatti flows on $Gr(2,3)$. Such a chart was found to be 
\begin{equation}\label{eq:t3dchart}
T := \begin{pmatrix} i & 0 & 1 \\ 0 & -i & 0 \\ 0 & 0 & 1 \end{pmatrix}.
\end{equation} 
We note that $\det (T) = 1$ so that there is no scaling of the Evans function. We then solve the Riccati equation on $Gr(1,3)$, namely the Grassmannian of lines in $\C^3$. We want to solve the Riccati equation for $y' = A_T y$, where $A_T= TAT^{-1}$.  
Writing $(A_T)_{i,j} = a_{ij}$ we have that $\alpha$ and $\beta$ solve
\begin{equation}\label{eq:ricstab}
\begin{aligned}
\alpha' & = a_{2,1} + a_{2,2} \alpha + a_{2,3}\beta - \alpha\left( a_{1,1} + a_{1,2} \alpha + a_{1,3} \beta \right) \\ 
\beta' & = a_{3,1} + a_{3,2} \alpha + a_{3,3} \beta - \beta \left( a_{1,1} + a_{1,2} \alpha + a_{1,3} \beta \right).
\end{aligned}
\end{equation}
We then shoot backwards from the stable eigenvector of $A_\pm(\lambda)$ (written as $(1,\alpha,\beta)$ ), and compare it to the evolution of the unstable subspace. This is evolved according to the following prescription. We write a pair of three vectors as a $2\times2$ matrix over a $1 \times 2$ matrix. Then splitting the linear ODE appropriately we want to solve 
\begin{equation}
\begin{pmatrix} X_{2,2} \\ Y_{1,2} \end{pmatrix}' = \begin{pmatrix} A_{2,2} & B_{2,1} \\ C_{1,2} & D_{1,1} \end{pmatrix} \begin{pmatrix} X_{2,2} \\ Y_{1,2} \end{pmatrix}.
\end{equation}
(where the subscripts indicate the size of the matrix). Inverting $X_{2,2}$ we have that $W_{1,2}:= Y_{1,2} X_{2,2}^{-1}$ solves
\begin{equation}
W_{1,2} = C_{1,2} + D_{1,1} W_{1,2}  - W_{1,2}A_{2,2} - W_{1,2}B_{2,1}W_{1,2}. 
\end{equation}
We then shoot forwards from the unstable subspace at $- \infty$ and compare with our evolution (in the Riccati coordinates) of the stable subspace. If we write $W_{1,2} = (w_1,w_2)$, then the Riccati-Evans function in this case is 
\begin{equation}
\label{eq:ricevans3d}
E(\lambda) = - w_1(z_0;\lambda) - w_2(z_0;\lambda) \alpha(z_0,\lambda) + \beta(z_0,\lambda)
\end{equation}
for an appropriately chosen $z_0$ (this will typically be $z_0 = 0$). These were numerically evaluated for a series of wave speeds and a positive real root was found in each case. We thus have a real positive eigenvalue for \cref{eq:stabode3d}. Figure \ref{fig:algwave2} shows the real (blue online) and imaginary (gold online) parts of the Riccati-Evans function for this nonstationary, for $\lambda>0$. %

For an additional example of the Riccati-Evans function in the case of a third order ODE, see \cite{LM23}.

\subsubsection*{Example 6: Nonstationary front}

Wylie and Miura \cite{WM06} found a nonstationary front solution in the case when the reaction term is given by
\begin{equation}\label{eq:3dsech}
g(u,v) = (u+v) - \delta u(u+\gamma)(u+1).
\end{equation}
with $\delta>0$ and $0<\gamma<1$. In this case, solutions to \eqref{eq:TWsys2} are given by
\begin{equation}\label{eq:3dsechsol}
u(z) = -\frac{1}{2} + \frac{1}{2} \tanh\left(\sqrt{\delta/8} z \right) \quad v(z) = - \frac{1}{c} u'(z) - u(z), 
\end{equation}
with 
\begin{equation}
c = \frac{\sqrt{\delta} (1-2 \gamma
   )+\sqrt{\delta (1-2 \gamma
   )^2+8}}{2 \sqrt{2}}.
\end{equation}

\begin{figure}
\includegraphics[scale=0.35]{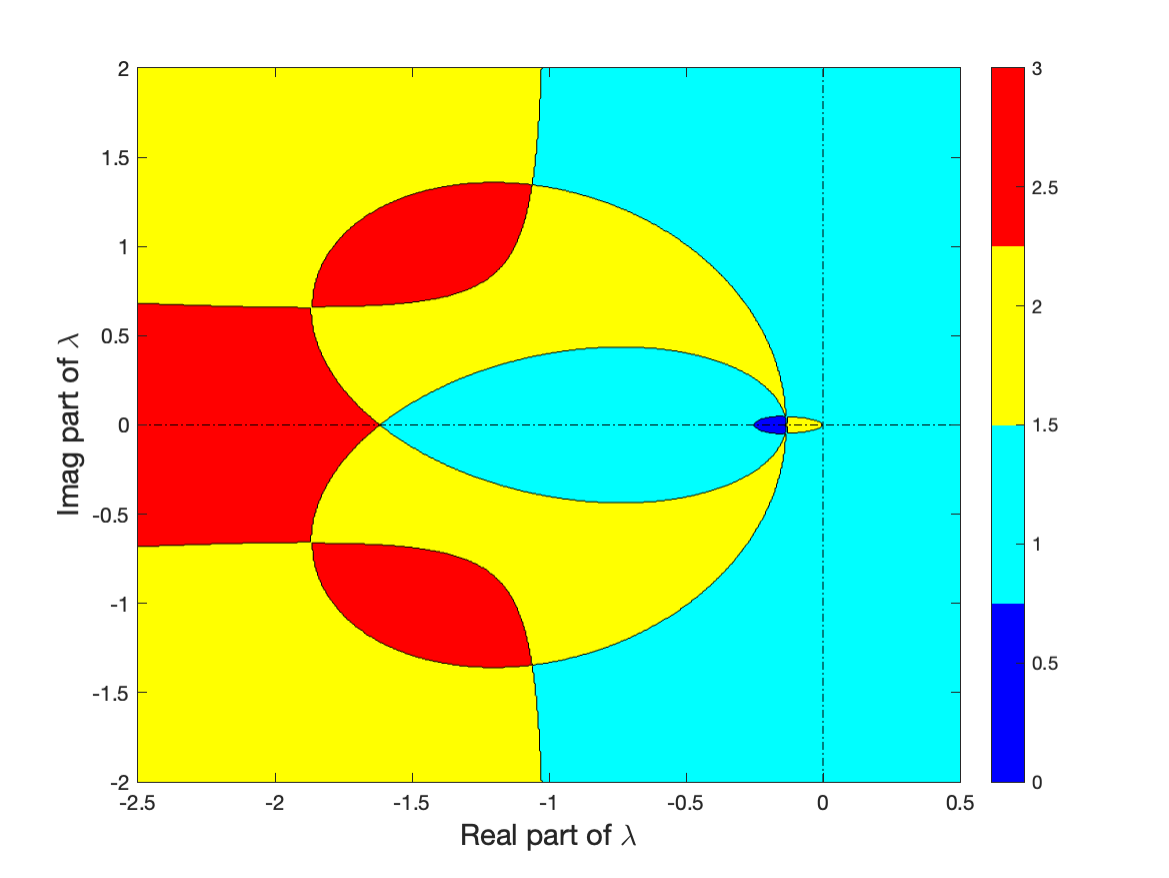}
\includegraphics[scale=0.35]{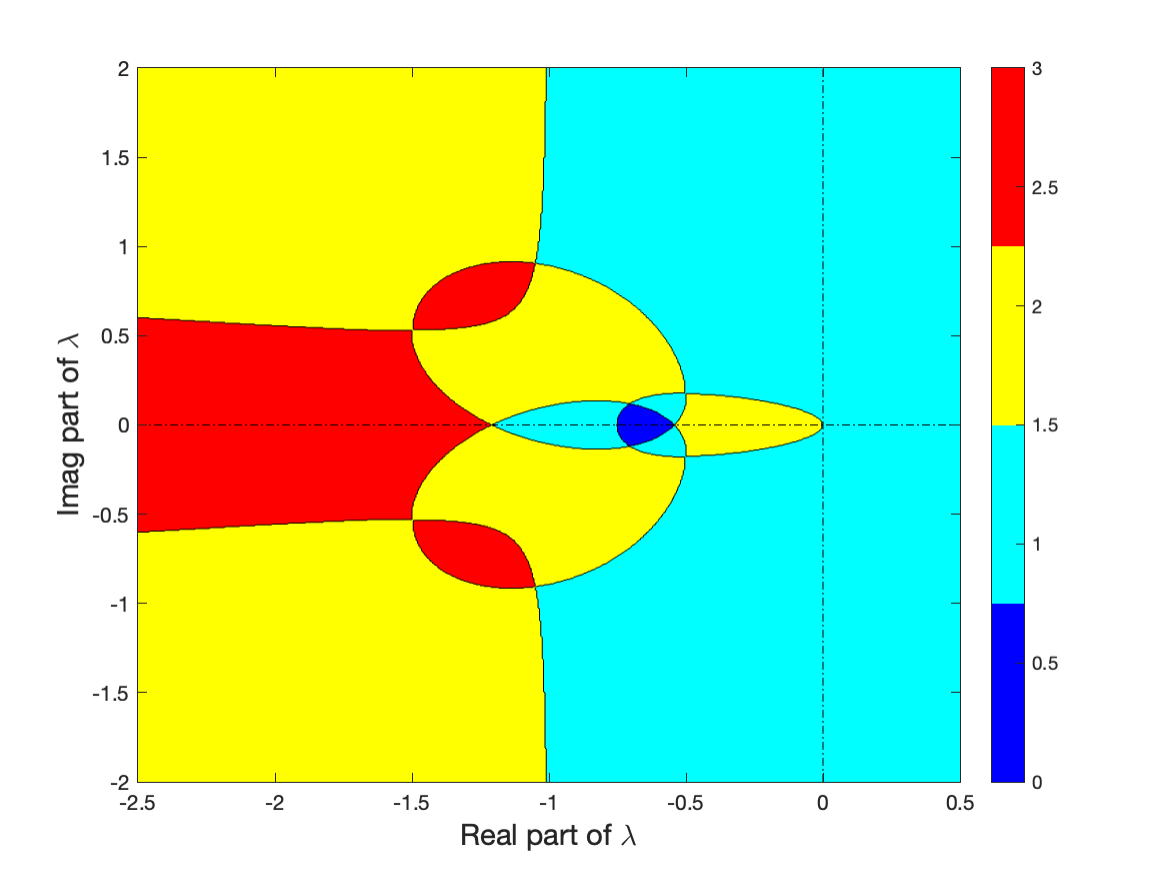}
\includegraphics[scale=0.4]{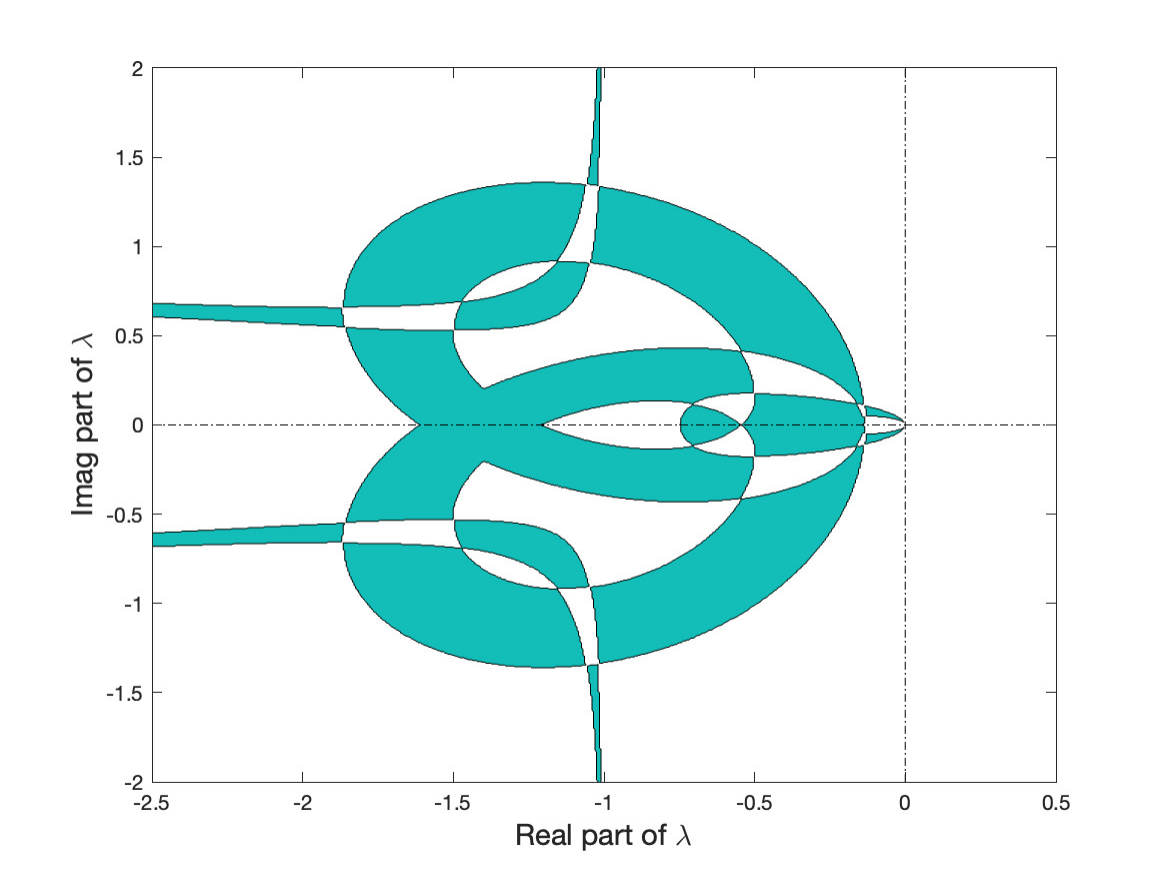}
\caption{Example 6: Nonstationary front, $D=0$. The top two figures show the number of unstable of eigenvalues of the far-field matrices $A_\pm(\lambda)$ for the nonstationary front given in equation \cref{eq:3dsechsol} ($A_-(\lambda)$ on the top left, $A_+(\lambda)$ top right). A value of $\lambda$ for which there is a mismatch indicates a non-zero Fredholm index and hence $\lambda$ is an element of the essential spectrum of $\cL$ which is shown in the bottom figure. The operator is not sectorial (as one expects because the ODE is $3\times 3$), but the real part of the essential spectrum is bounded above. Moreover, we see that the essential spectrum is contained in the left half plane. Here, $\gamma = \frac{3}{4}$ and $B = 1$.}  %from $A_+(\lambda)$.
 \label{fig:sechfredholm}
\end{figure}

The essential spectrum can be obtained by setting $D=0$ in \eqref{FredholmBorder} which yields 
\begin{equation}\label{eq:disprel6}
\begin{aligned}
\lambda_{1,2}&=\frac{1}{2}\left(-k^2-\delta\gamma \pm\sqrt{(k^2+2k+\delta\gamma)(k^2-2k+\delta\gamma)}+2ick\right),\\
\lambda_{3,4}&=\frac{1}{2}\left(-k^2-\delta(1-\gamma)\pm\sqrt{(k^2+2k+\delta(1-\gamma))(k^2-2k+\delta(1-\gamma))}+2ick\right),
\end{aligned}
\end{equation}
so that the continuous spectrum is contained in the left-half plane. An example plot of the continuous spectrum is shown in Figure \ref{fig:sechfredholm}. 

\begin{figure}
\includegraphics[scale=0.25]{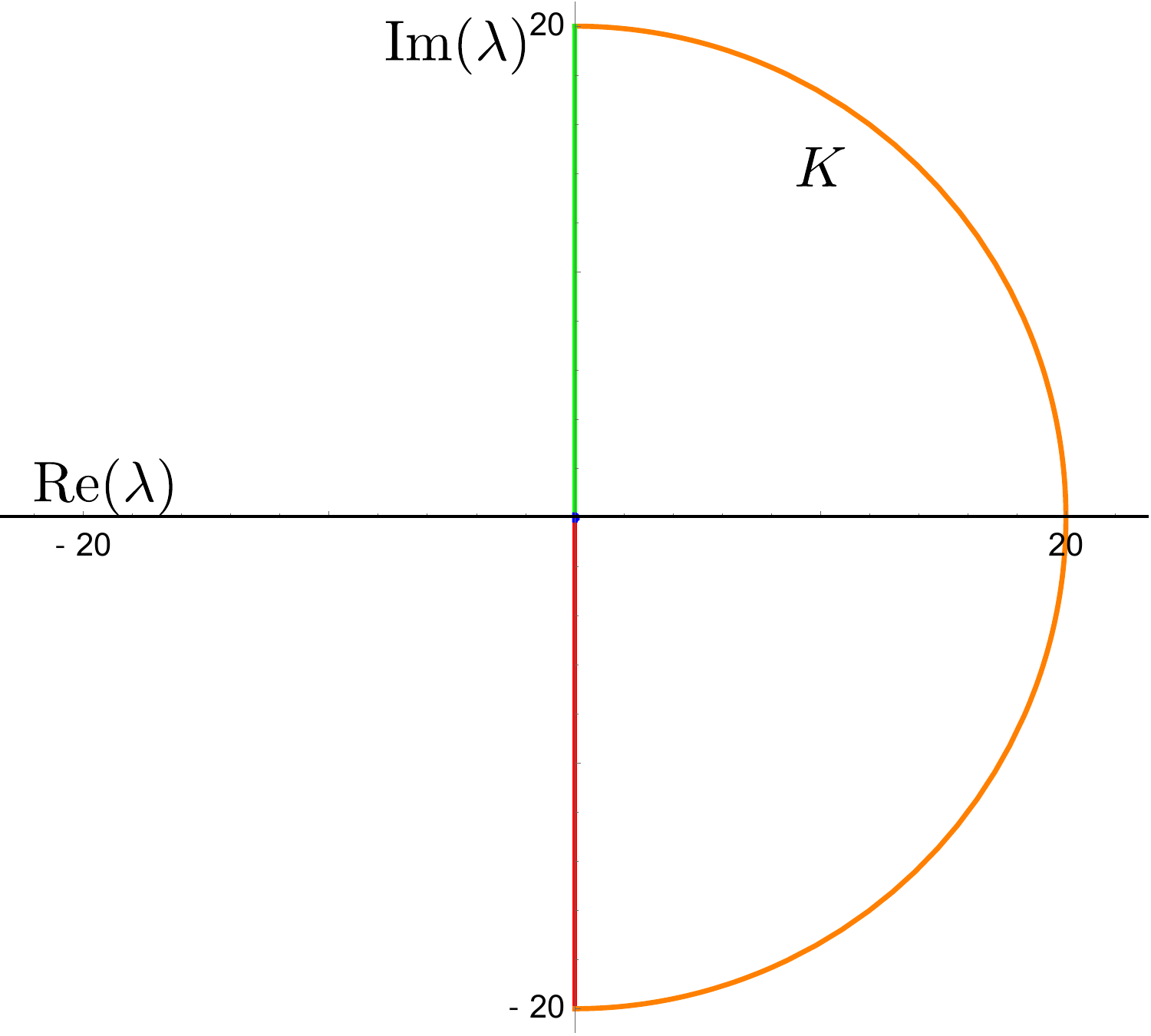}
\includegraphics[scale=0.25]{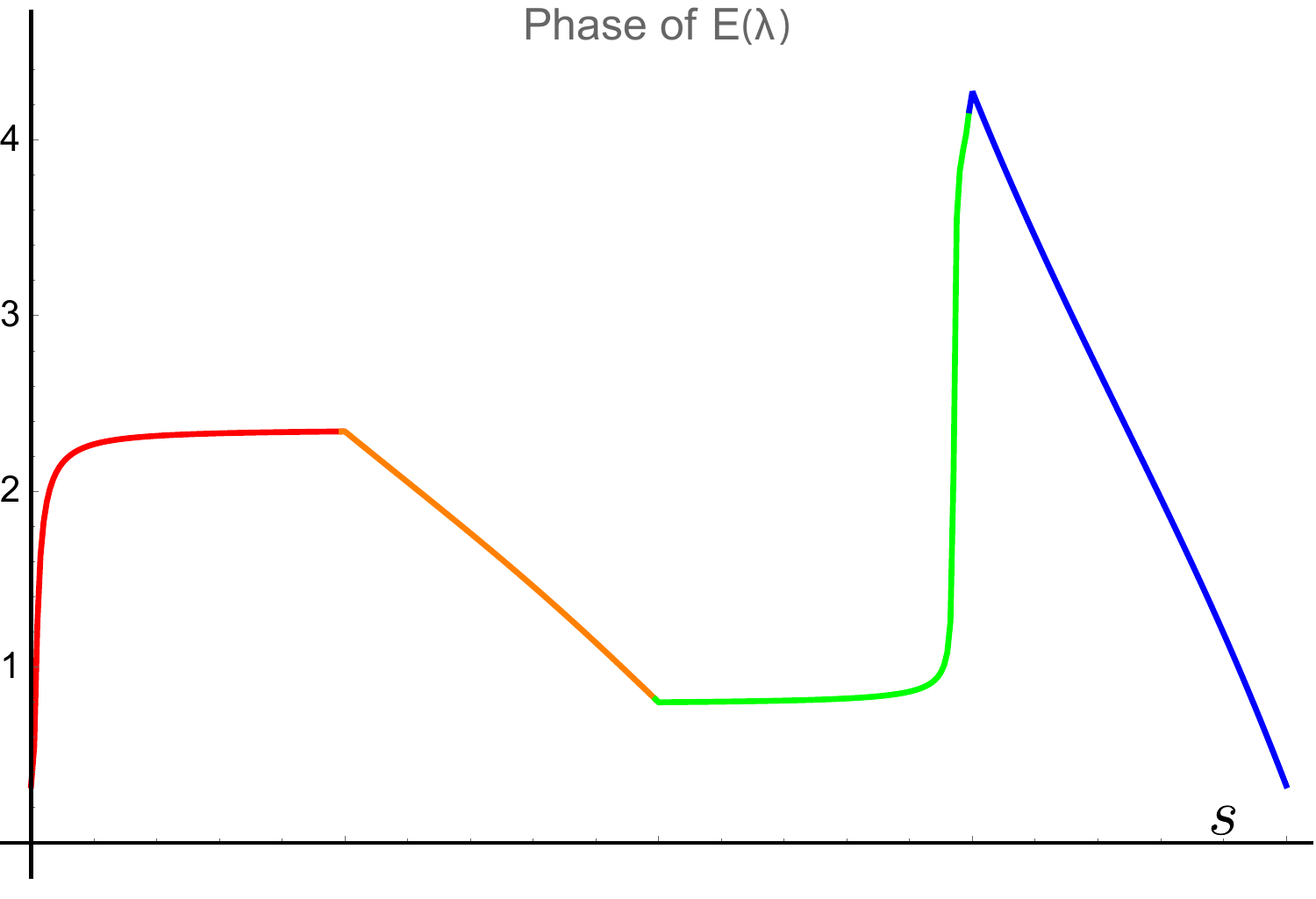}\\
\includegraphics[scale=0.25]{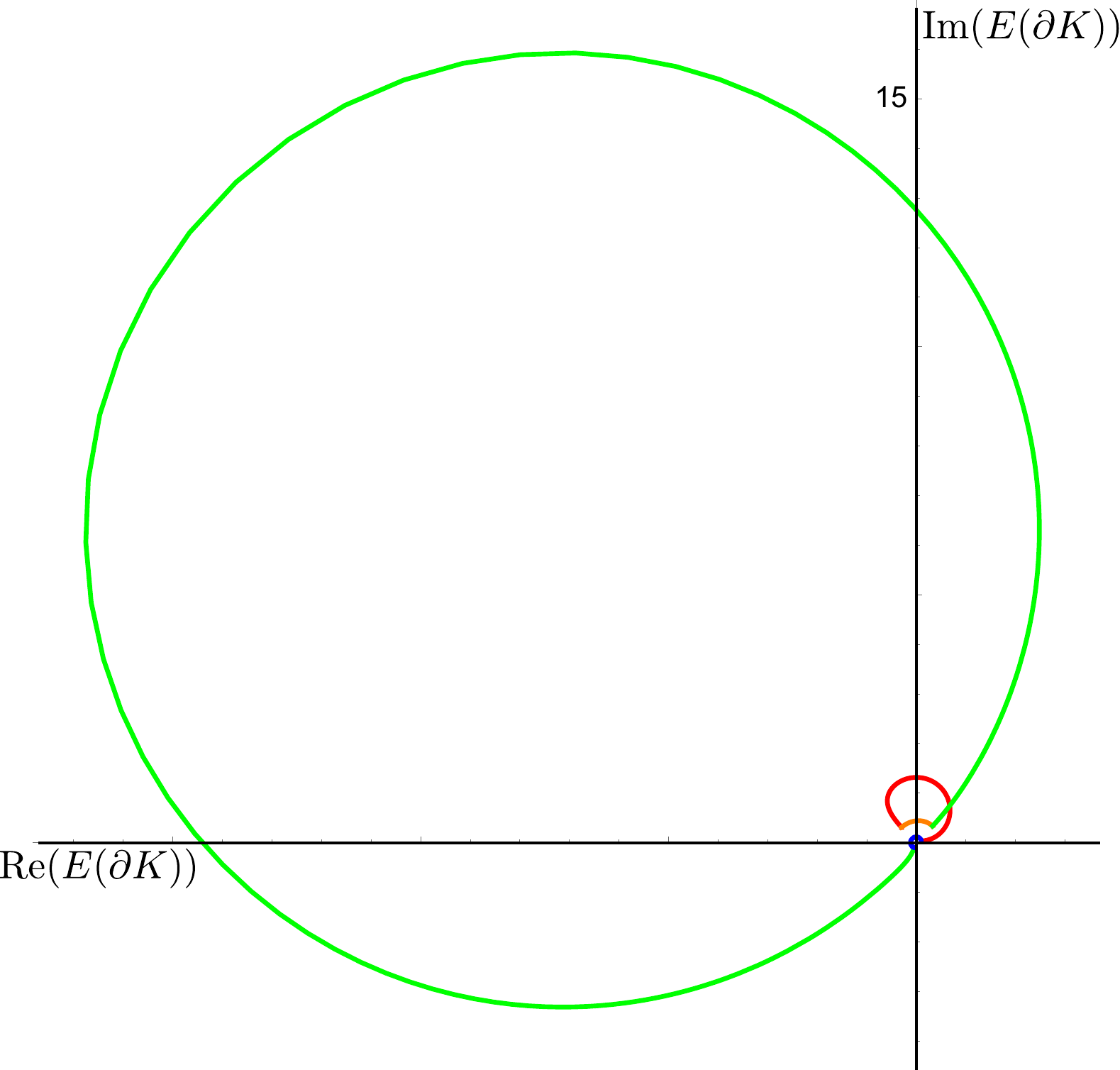}
\includegraphics[scale=0.25]{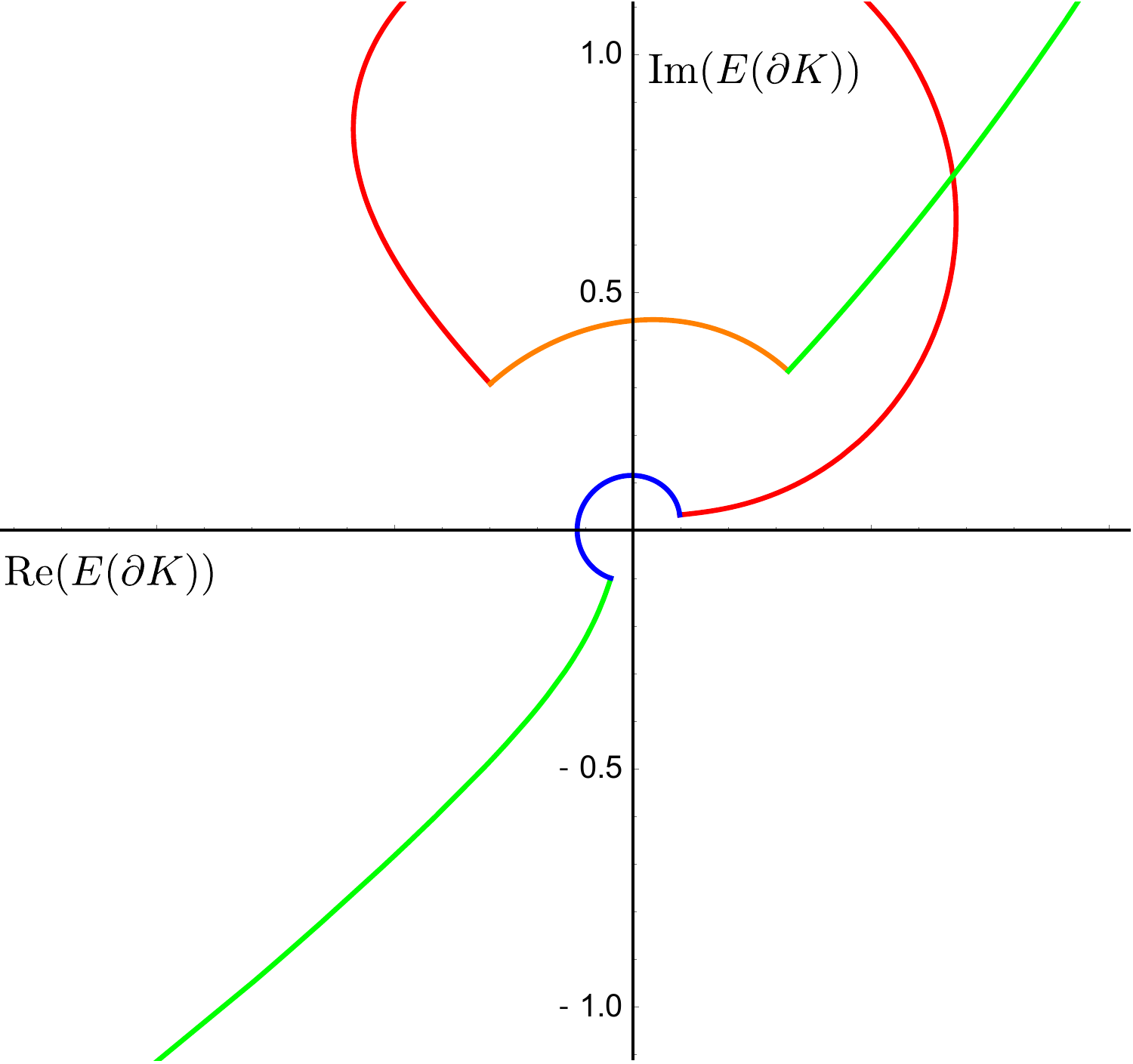}
\caption{Example 6: Nonstationary front, $D=0$. Determination of the point spectrum. Top right panel shows the phase of $E(\lambda)$ \eqref{eq:ricevans3d} as $\lambda$ moves along the path shown in the top left panel (for solutions \eqref{eq:3dsechsol}). Bottom panels show the image of $E(\lambda)$, with the bottom right zooming in at the origin. Here, $\gamma = \frac{3}{4}$ and $B=1$. The graph shows no winding and thus we infer spectral stability of the fronts in \eqref{eq:3dsechsol}.}\label{fig:sechricstabimag}
\end{figure}

Employing the same methods as described in Example 5, we search for point spectrum in the right-half plane. For solutions \eqref{eq:3dsechsol}, we tested a wide range of parameters and found that the Riccati-Evans function \eqref{eq:ricevans3d} has no zeros for real $\lambda$. 
In order to determine if there are complex eigenvalues with positive real part we introduce the closed curve in the complex plane as shown in the upper left panel of Figure~\ref{fig:sechricstabimag}. This curve consists of a large semi-circular arc (orange), a small semi-circular arc (blue) going around the origin and two straight lines connecting the two arcs (green and red). We then integrate along this curve to obtain the winding number.
The top right panel shows the phase of the Riccati-Evans function, $E(\lambda)$ as we wind from red to orange to green to blue along the path shown in the top left panel. The bottom panels show the image of $E(\lambda)$, illustrating that there is no winding, corresponding to there being no change in phase in the top right panel. Similar results were observed across a broad range of parameters, leading us to conclude that the solutions \eqref{eq:3dsechsol} exhibit spectral stability.

\subsection{Stationary waves, $v$ diffuses, $D>0$, $c=0$}

Here we use the Riccati-Evans function to investigate the stability of stationary pulse and front solutions when $v$ is diffusing. 

\subsubsection*{Example 7: Stationary pulse}

The stationary pulse solution described in Example 1 \eqref{eq:Eg1sol} is also valid in the case of unequal diffusivities. In this case, $D\ne1$ and the stability must be determined numerically. The stability problem \eqref{eq:eigenvalue1} can be recast as a system of four first order equations,
\begin{equation} \label{eq:fourth}
\begin{pmatrix}
p \\ q \\ r \\ s  
\end{pmatrix}' = \begin{pmatrix} -c & 0 & 1 & 0 \\ 0 & -\frac{c}{D} & 0 & \frac{1}{D} \\ \lambda -  g_u  & - g_v & 0 & 0 \\ g_u & \lambda + g_v& 0 & 0 \end{pmatrix} \begin{pmatrix} p \\ q \\ r \\ s \end{pmatrix} =A(\lambda)\,\begin{pmatrix}
p \\ q \\ r \\ s  
\end{pmatrix}
\end{equation}
where $r = p' + c p $ and $Dq' + c q = s$. Both far field matrices are the same (since pulse solutions have the same behaviour in both directions), and are given by
\begin{equation*}
A_{\pm}(\lambda)=\lim_{z\to\pm\infty}A(\lambda)=\begin{pmatrix} 0 & 0 & 1 & 0 \\ 0 & 0 & 0 & \frac{1}{D} \\ \lambda +4  & 0 & 0 & 0 \\ -4 & \lambda & 0 & 0 \end{pmatrix} \begin{pmatrix} p \\ q \\ r \\ s \end{pmatrix}
\end{equation*}
The characteristic polynomial of $A_{\pm}(\lambda)$ is
\[
\det(A_{\pm}(\lambda)-\mu)=\frac{(\lambda+\mu^2+4)(\lambda+D\mu^2)}{D}.
\]
The continuous spectrum can be found by setting $\mu=ik$ (equivalent to setting $p(z)\sim\e^{ikz}$ and $q(z)\sim\e^{ikz}$ in equations \eqref{eq:eigenvalue1}) whence we see that the continuous spectrum is the negative real axis, $(-\infty,0]$. The spatial eigenvalues are
\[
\mu_{1,2}=\pm\sqrt{\frac{\lambda}{D}}\,,\qquad \mu_{3,4}=\pm\sqrt{\lambda+4}\,,
\]
so that, away from the continuous spectrum, $A(\lambda)$ is hyperbolic with two eigenvalues with positive real part, and two with negative real part.

\begin{figure}
\includegraphics[scale=0.25]{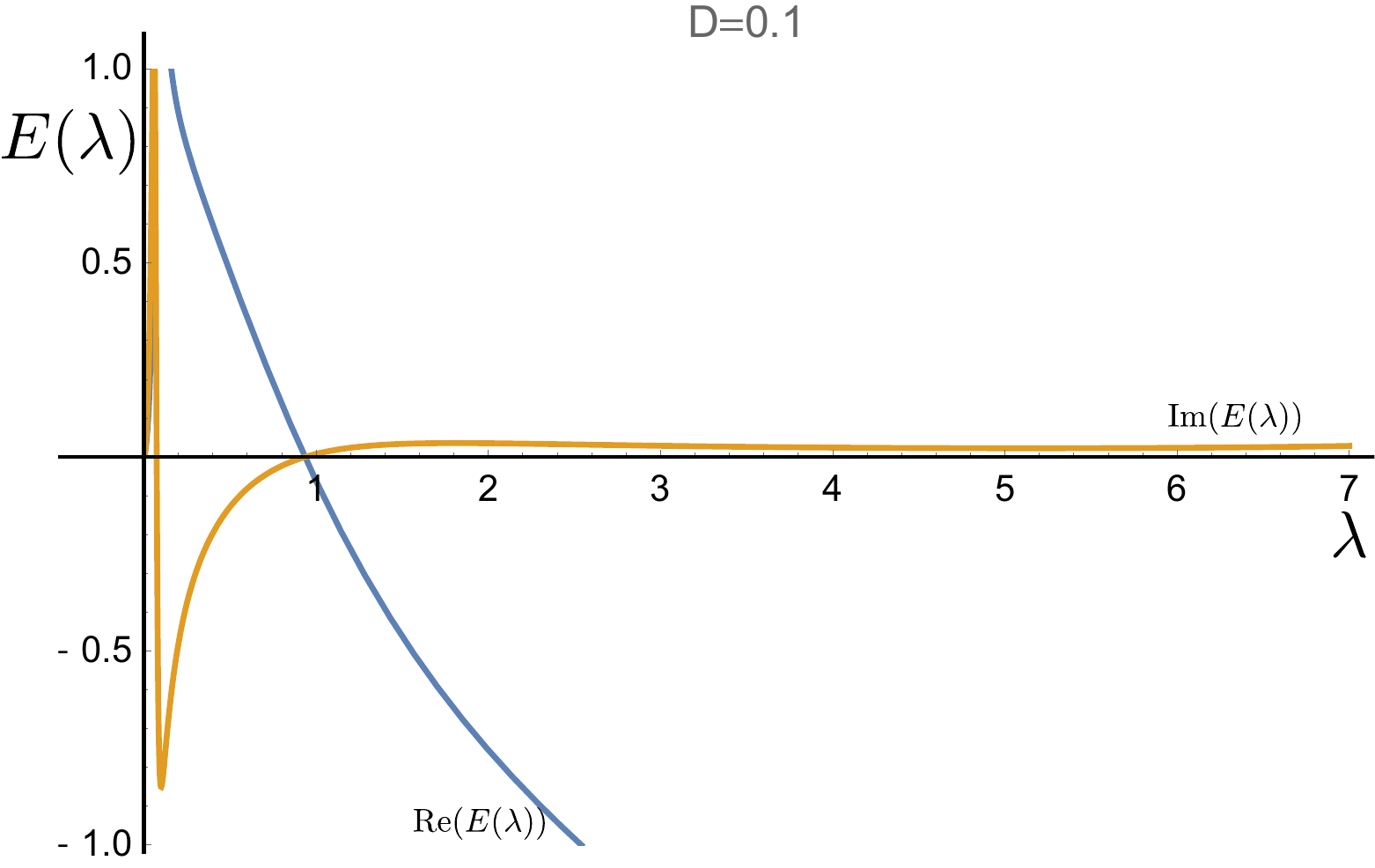}
\includegraphics[scale=0.25]{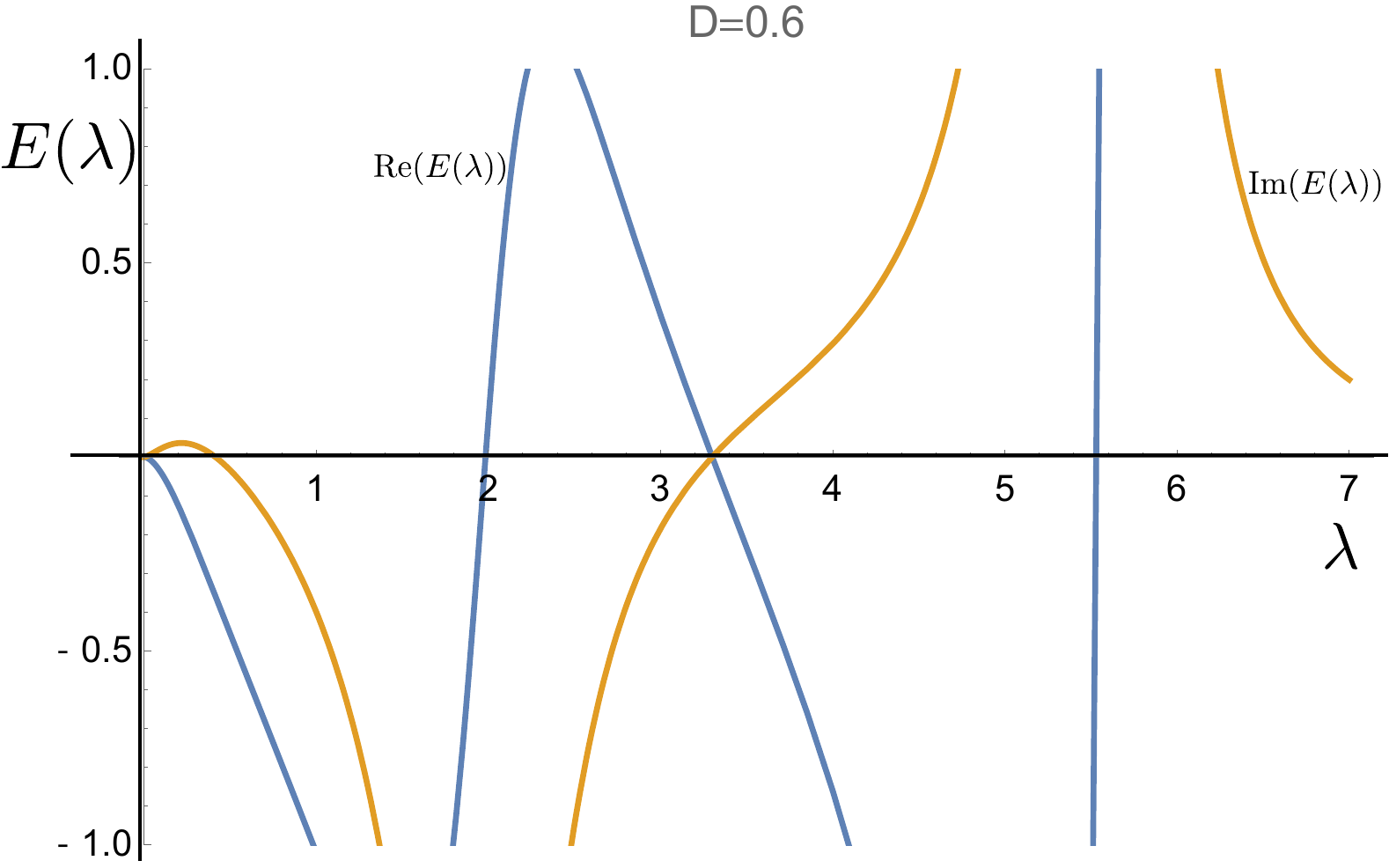}
\includegraphics[scale=0.25]{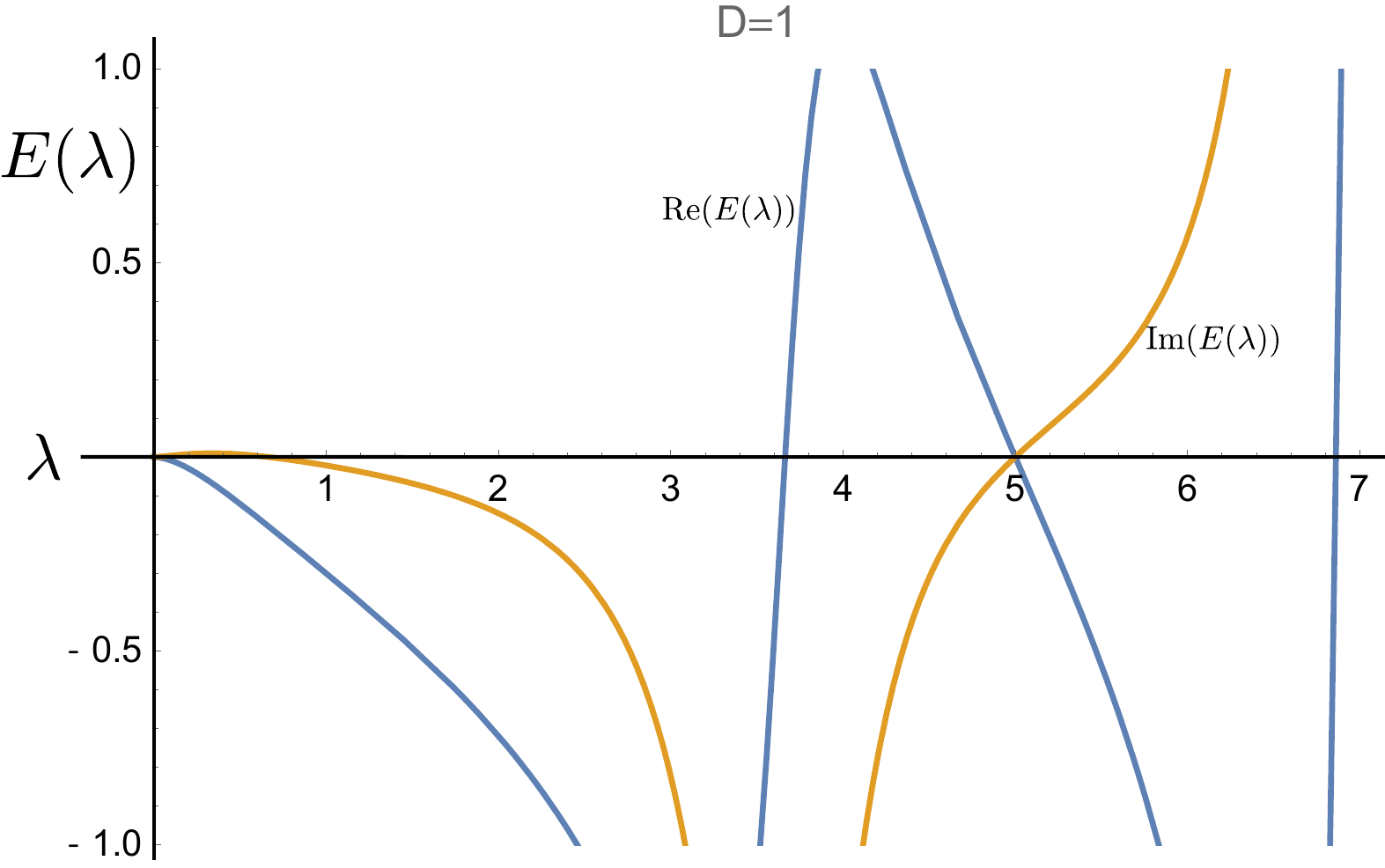}
\includegraphics[scale=0.25]{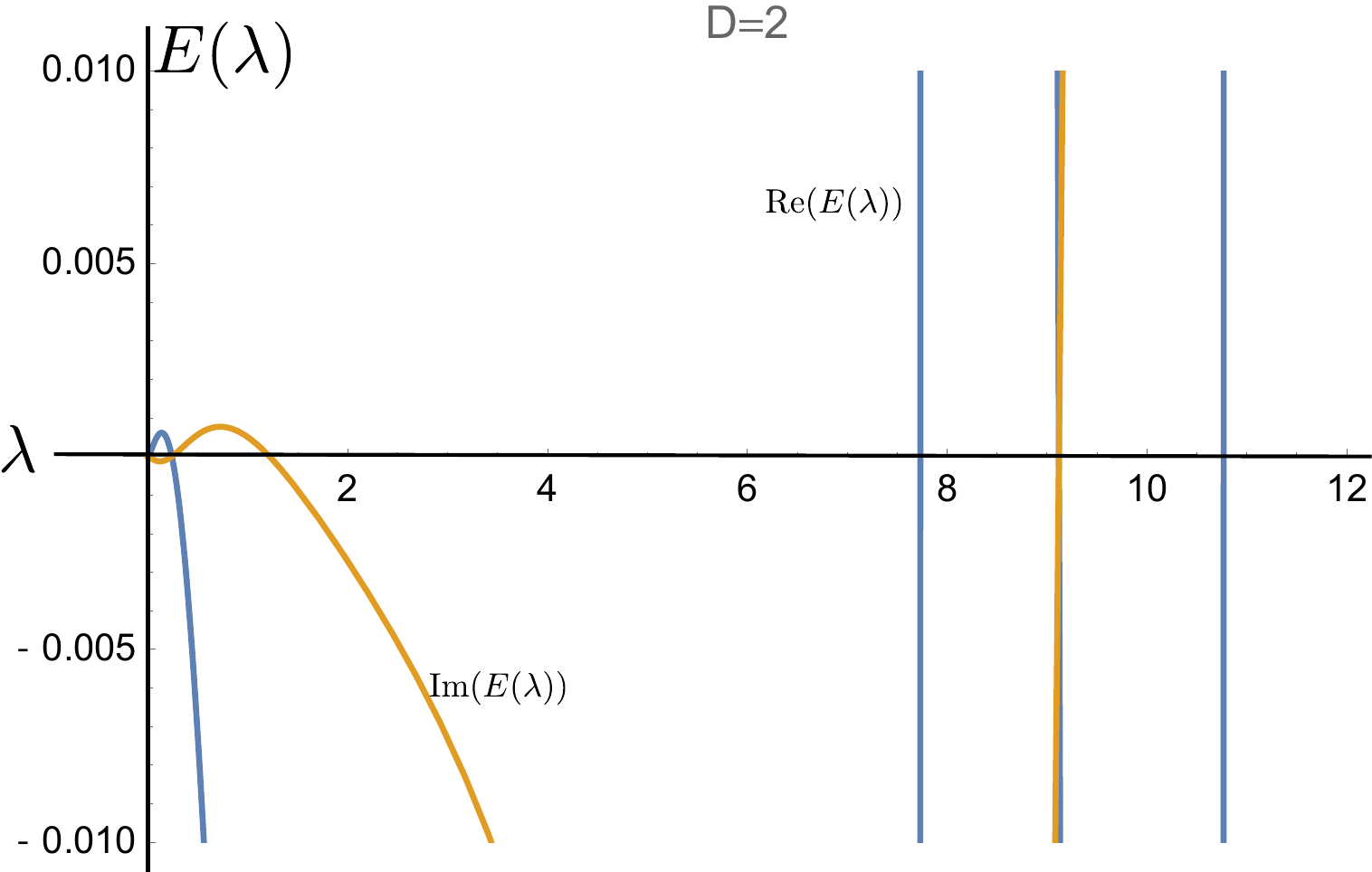}
\caption{Example 7: Stationary pulse, $D>0$. Real (blue online) and imaginary (yellow) parts of the Riccati-Evans function for  \eqref{eq:Eg1sol}. When $D = 1$ there is a simple eigenvalue at $\lambda = 5$, which moves to the right as $D$ is increased and to the left as $D \to 0$. For $D>1$, a second eigenvalue occurs, bifurcating out of the essential spectrum.}
\label{fig:ricevstanding}
\end{figure}

To complete the investigation of the stability of this solution, we must also find the point spectrum. Using the chart given by
\begin{equation}\label{eq:Tchart}
T=\begin{pmatrix} i & 0 & 1 & 0 \\ 0 & 1 & 0 & 1 \\ 0 & 0 & -i & 0 \\ 0 & 0 & 0 & 1 \end{pmatrix}
\end{equation}
and following the technique described in Example 5,
we computed the Evan’s function for a wide array of parameters and in all cases we found an eigenvalue on the positive real axis. Some examples of the Evan’s function are shown in Figure \ref{fig:ricevstanding}.
Since, in all cases the point spectrum lies on the positive real axis, we conclude that solutions of this type are always unstable.

\subsubsection*{Example 8: Stationary front}

The stationary pulse solution described in Example 2 is also valid in the case of unequal diffusivities, $D\ne1$. When $D\ne1$, we can determine the stability numerically using the Riccati-Evans functions as in the last examples.

\begin{figure}
\includegraphics[scale=0.25]{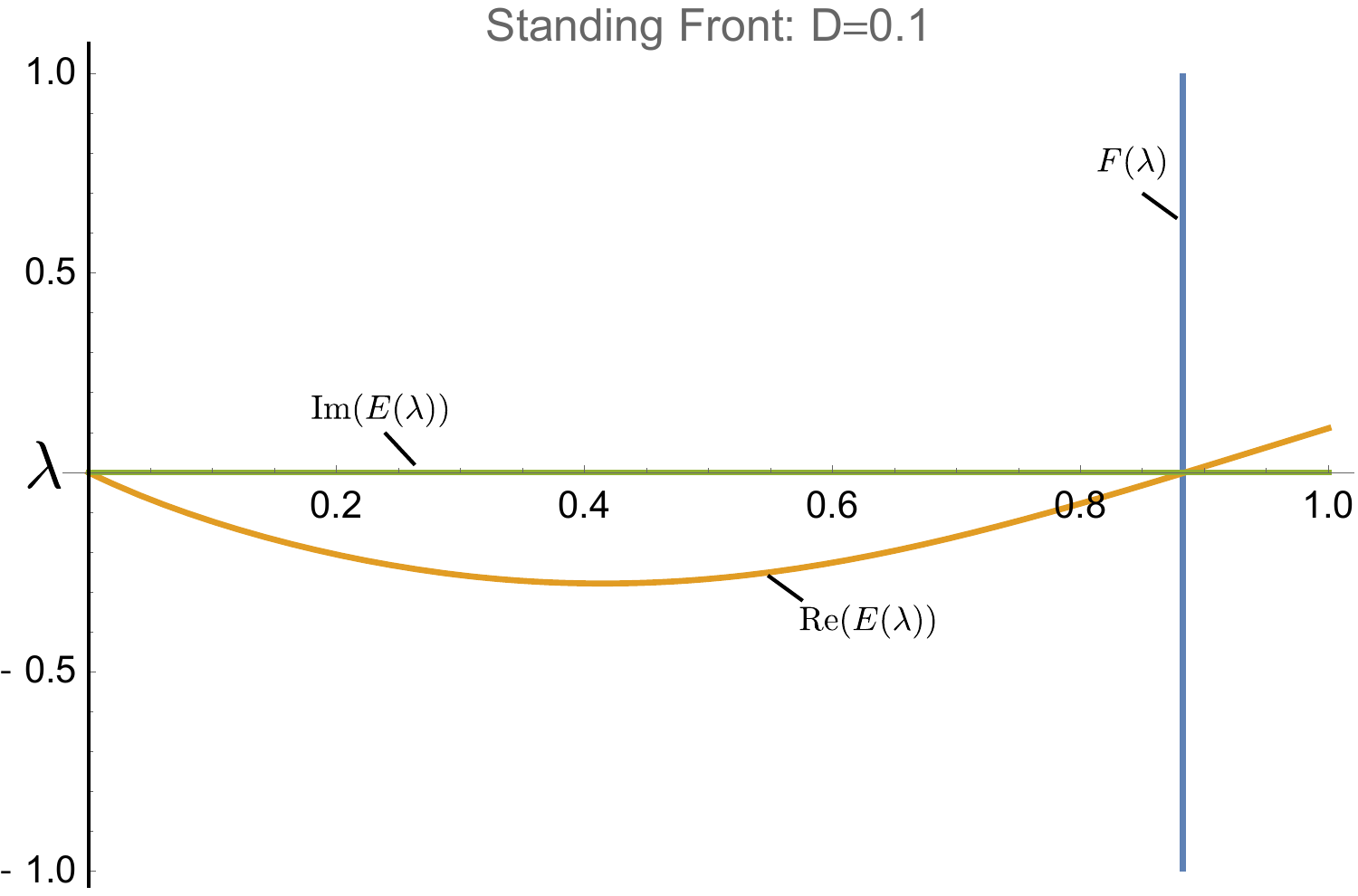}
\includegraphics[scale=0.25]{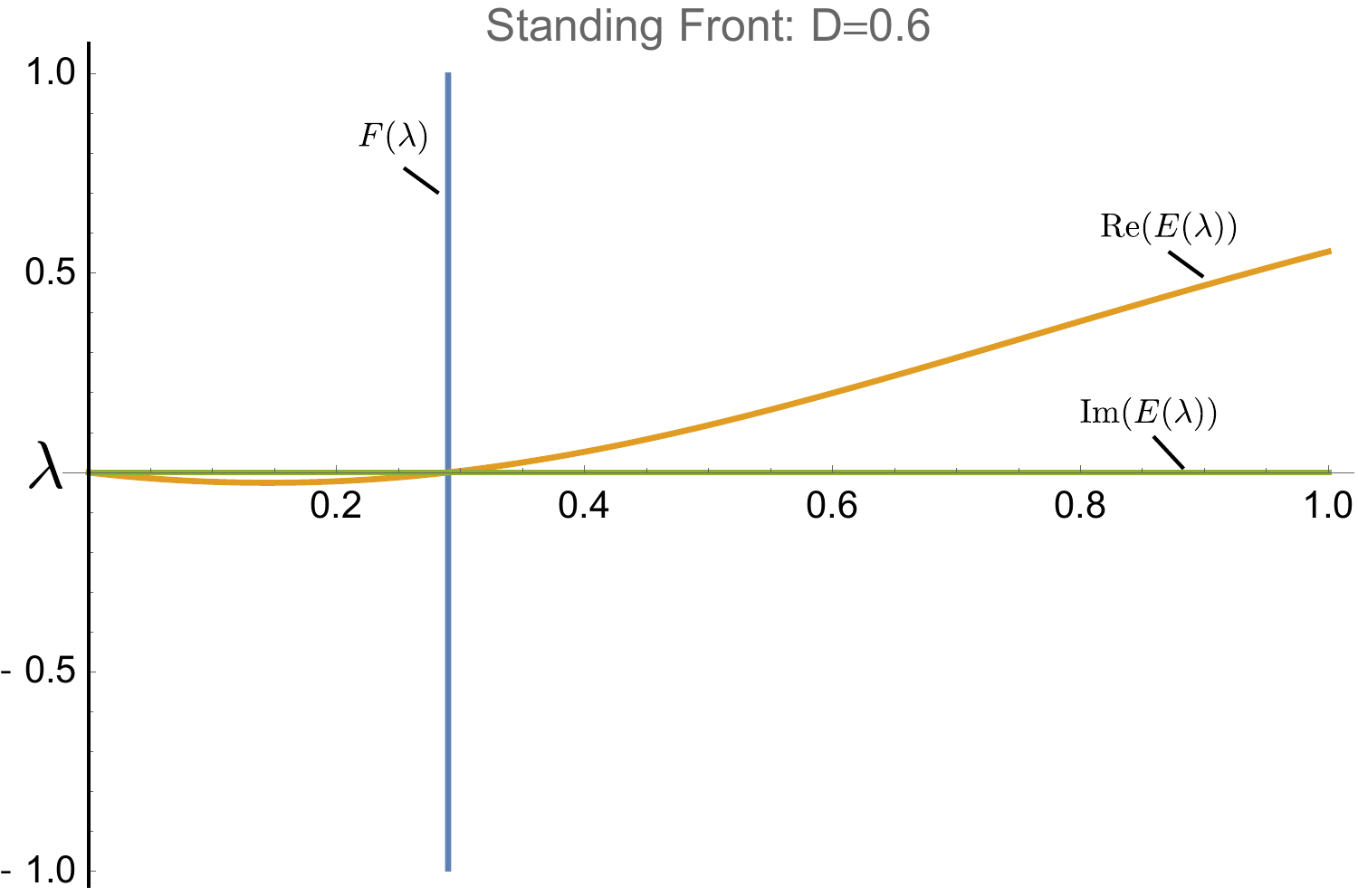}
\includegraphics[scale=0.25]{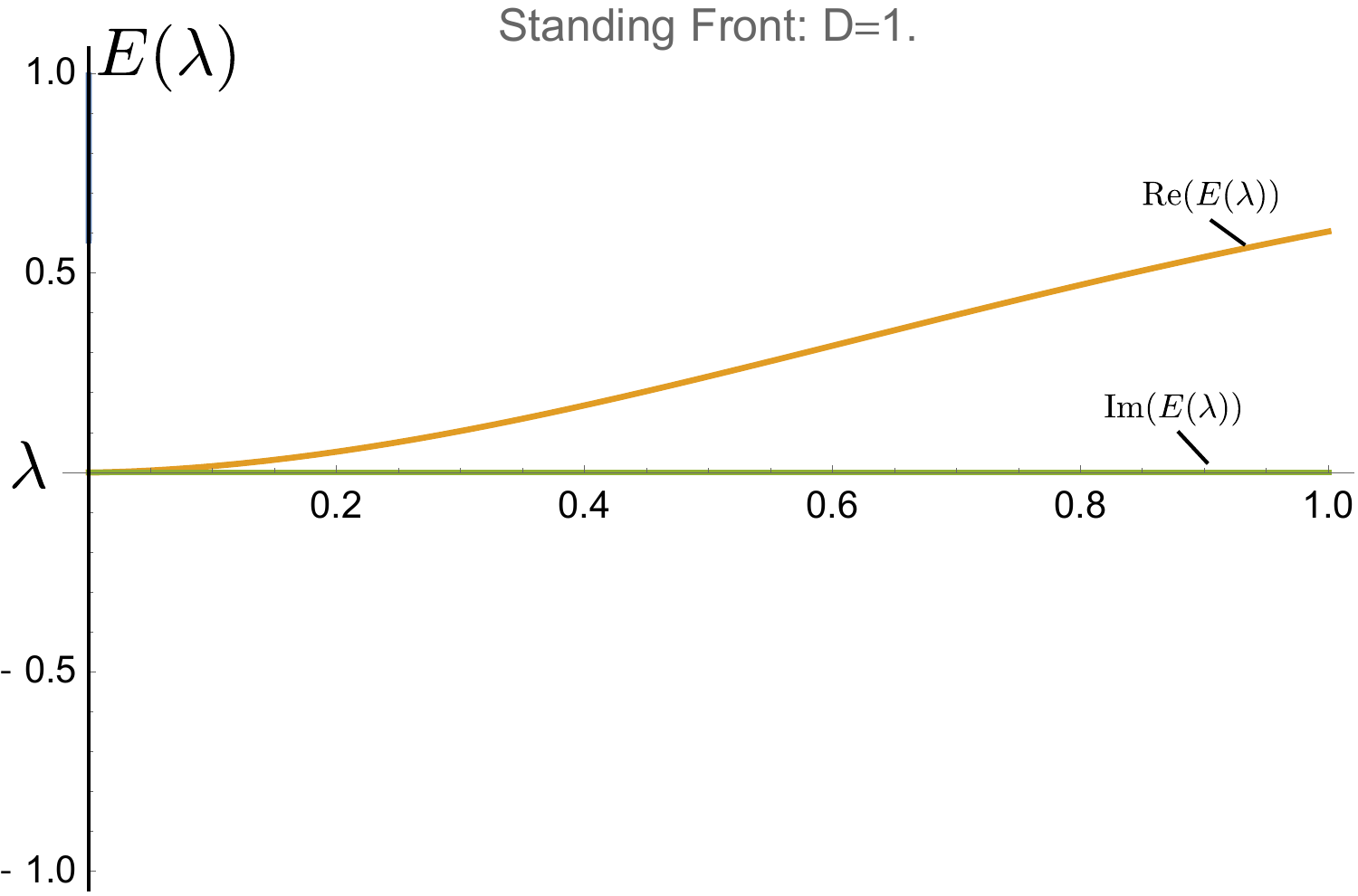}
\includegraphics[scale=0.25]{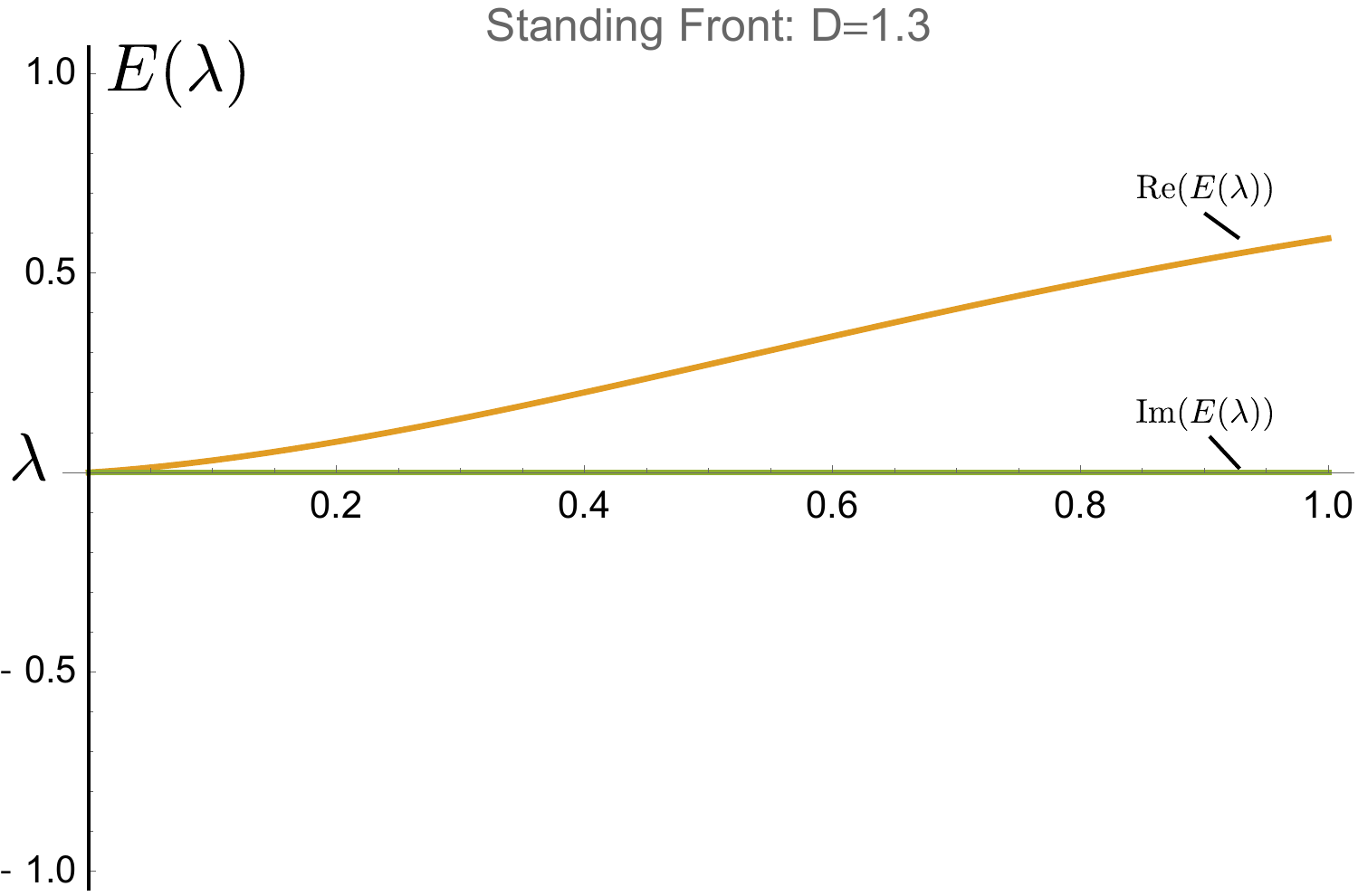}
\caption{Example 8: stationary front, $D>0,~c=0$. A plot of the real (yellow online) and imaginary (green) parts of the Riccati-Evans function with the chart given in \cref{eq:Tchart2}. When $D <1$ there is a simple eigenvalue, which moves to the left as $D$ is increased and to the right as $D \to 0$. For $D\geq1$, no real eigenvalues appear. The blue (online) is the standard Evans function in this case -- note that it is numerically very unstable and does not produce reliable results for $\lambda$ much larger than $5$ in this case.}
\label{fig:ricevstandfront}
\end{figure}

Despite the fact the far field behaviour is not the same in both directions, the far field matrices are the same,
\begin{equation}
A_\pm(\lambda) = 
\left(
\begin{array}{cccc}
 0 & 0 & 1 & 0 \\
 0 & 0 & 0 & \frac{1}{D} \\
 \lambda  & -4 D & 0 & 0 \\
 0 & 4 D+\lambda  & 0 & 0 \\
\end{array}
\right).
\end{equation}
The characteristic equation in the case where $\mu=ik$ is
$$
\det(A(\lambda) - i k) = \frac{(k^2 + \lambda)(\lambda + D(4 + k^2))}{D},
$$
and so the essential spectrum in this case is $(-\infty,0]$ again. The Riccati-Evans function can be used to find any point spectrum. Using the chart
\begin{equation}\label{eq:Tchart2} 
T := \left(
\begin{array}{cccc}
 i & 0 & 1 & 0 \\
 0 & -i & 0 & 1 \\
 0 & 0 & 1 & 0 \\
 0 & 0 & 0 & 1 \\
\end{array}
\right)
\end{equation}
we find that there is a root at $\lambda = 0$ for all values of $D$. We considered real eigenvalues for a wide range of parameters and found that for $D<1$ there is always a single positive eigenvalue. On the other hand, for $D \geq 1$, we never observed real positive eigenvalues. Thus, the stationary front in Example 8 is unstable when $D<1$. 
We show a representative sample of the real and imaginary parts of the Riccati-Evans function in Figure \ref{fig:ricevstandfront}.

\begin{remark}
As with the pulse solution from Example 7, the front solution and stability analysis in Example 8 is also valid in the case when $D=1$ (see Section \ref{sec:D1pulse}). Example 2 showed analytically that the continuous spectrum is the negative half-line, $(-\infty,0]$, and that there was no point spectrum in the right-half plane, in agreement with the numerical stability calculations here. This can be seen for Example 7 in the bottom left panel of Figure \ref{fig:ricevstanding} which shows a simple eigenvalue at $\lambda=5$. For Example 8, the bottom left panel of Figure \ref{fig:ricevstandfront} shows no real eigenvalues in the right-half plane. 
\end{remark}

\begin{figure}
\includegraphics[scale=0.25]{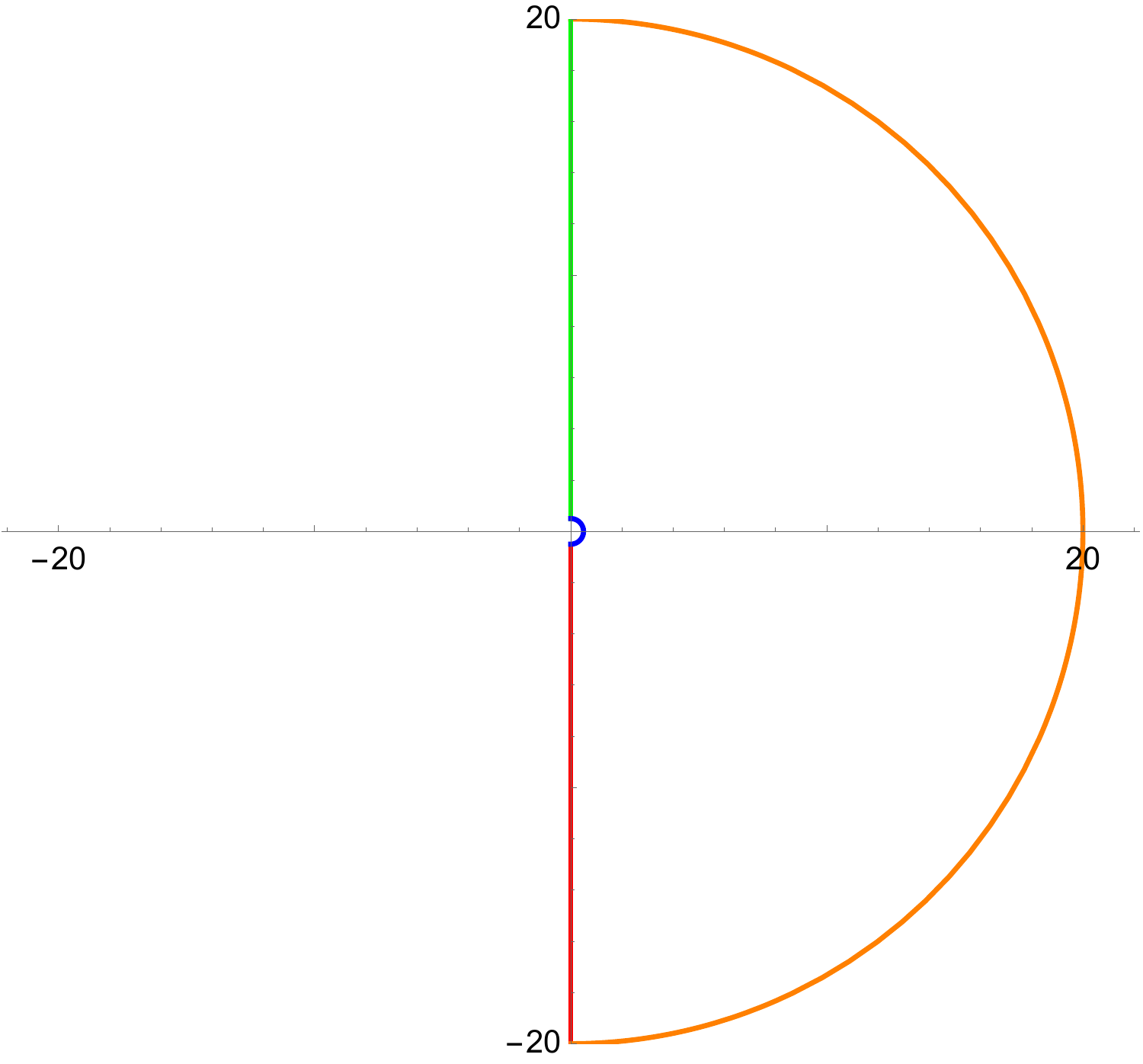}
\includegraphics[scale=0.25]{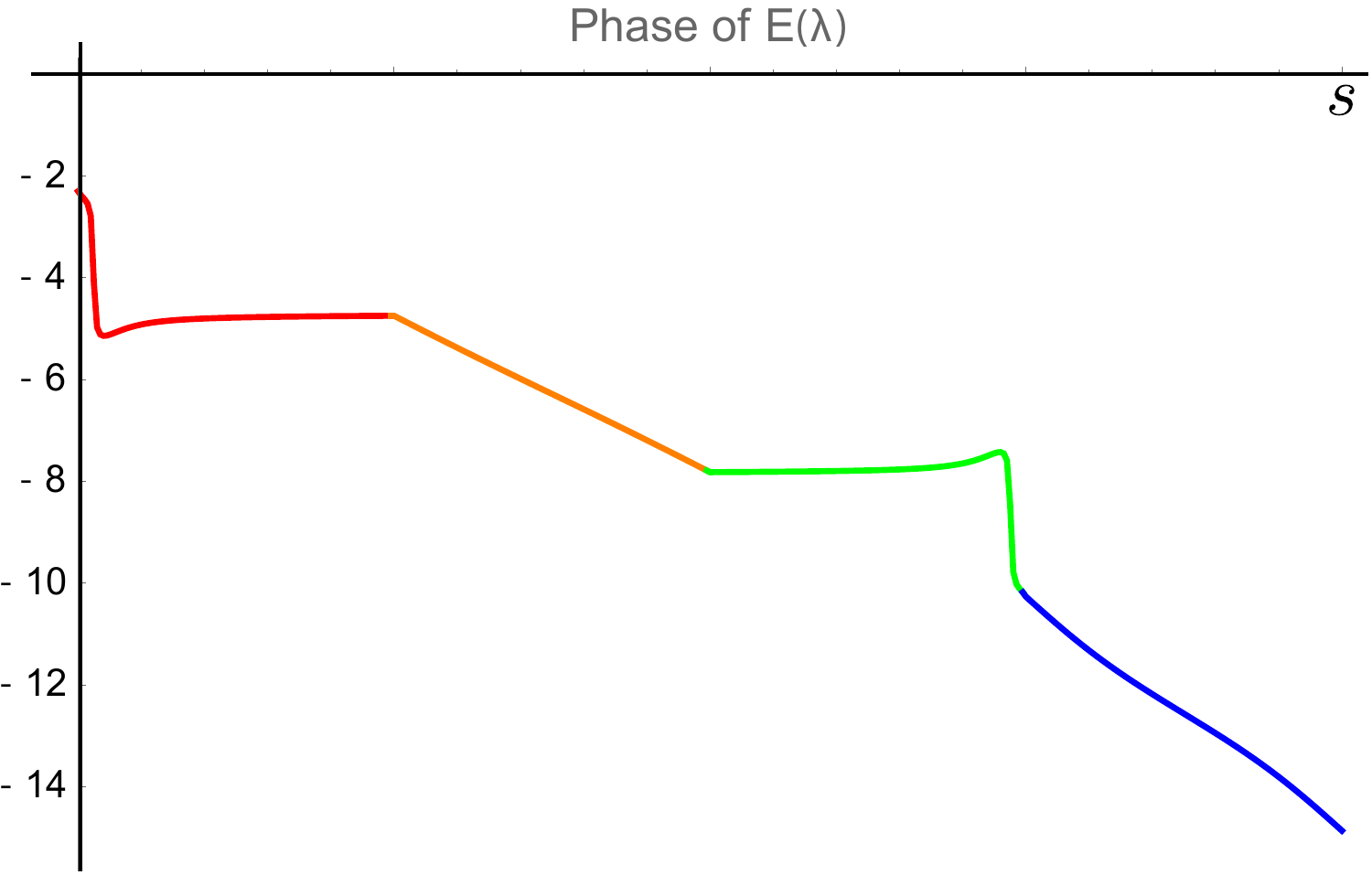}
\includegraphics[scale=0.25]{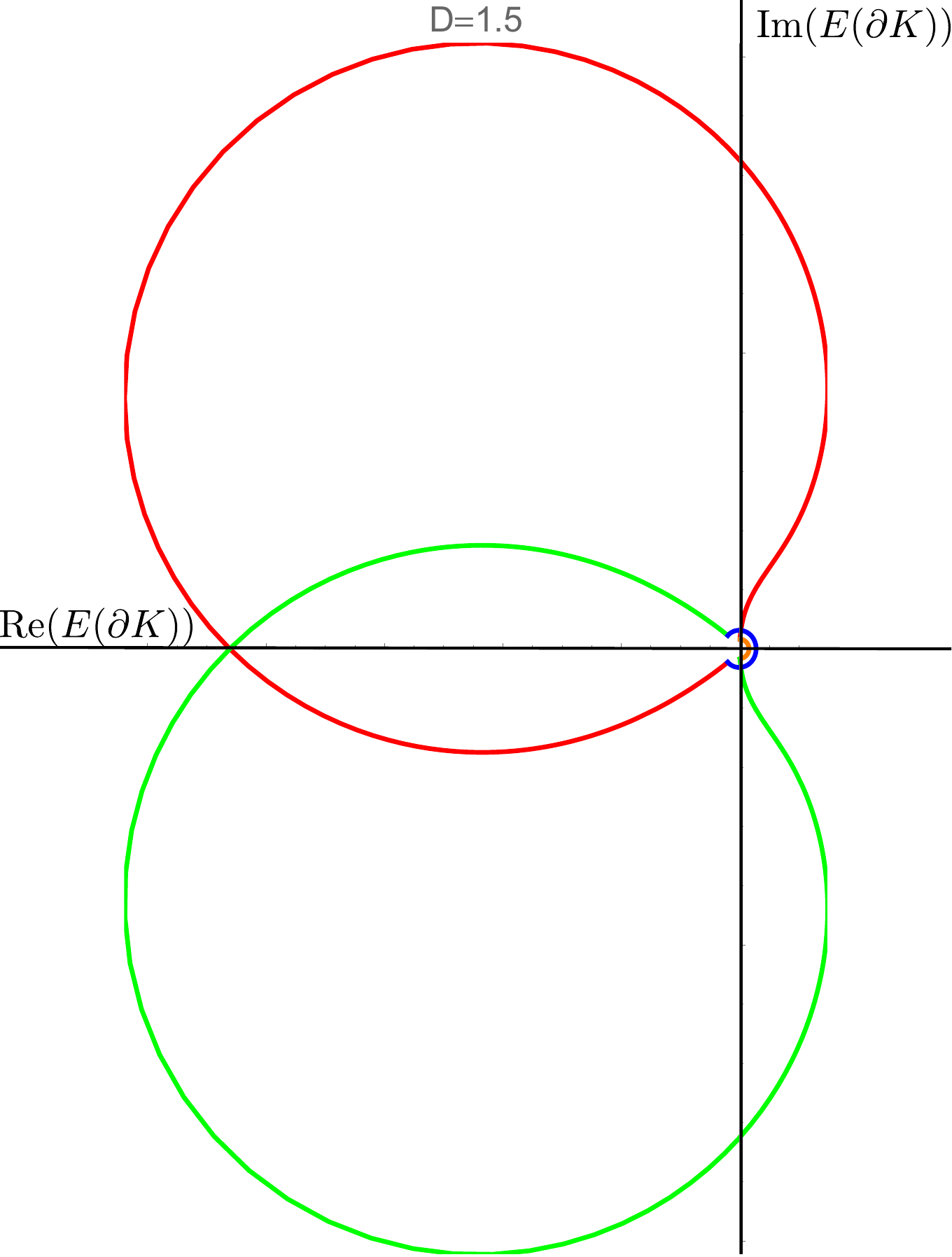}
\includegraphics[scale=0.25]{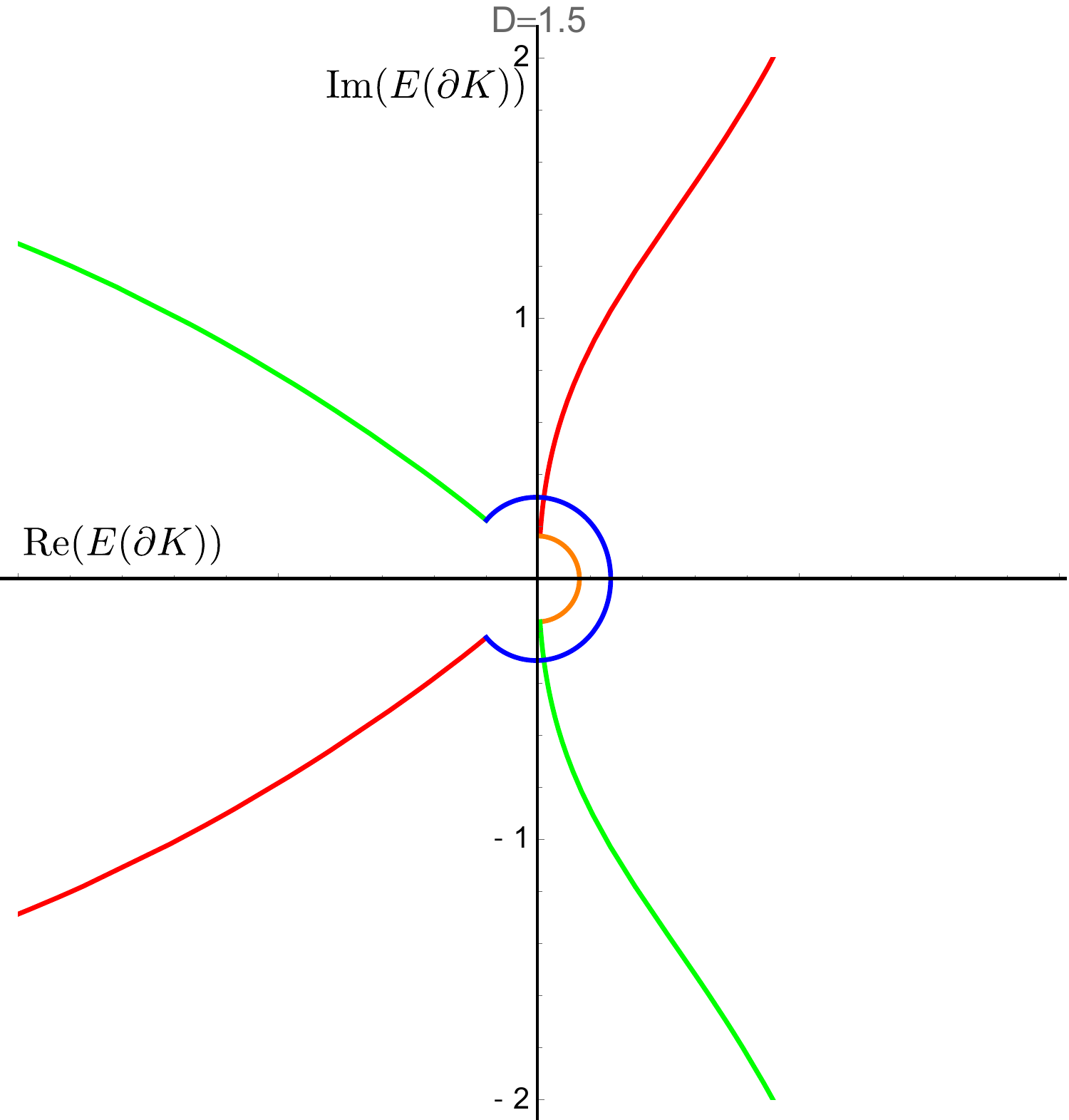}
\caption{Example 8: stationary front, $D>0,~c=0$. Top left shows the region $K$ in the complex plane over which we evaluate $E(\lambda)$. Top right shows the change in phase of $E(\lambda)$ as the boundary of $K$ is traversed. Bottom panels show the image of the boundary of $K$ under $E(\lambda)$ for $D = 1.5$, colours correspond to the first panel.}
\label{fig:dwind}
\end{figure}

In order to determine the spectral stability in the case when $D>1$ we need to investigate the possibility of complex eigenvalues with positive real part. This can be determined by investigating the phase change of the Evans function and the winding number of the Riccati-Evans function. Consider a large region $K$ in the complex plane bounded by a large outer half disc, the imaginary axis and a small inner half disc, as shown in the top left panel of Figure \ref{fig:dwind}. In this instance we can compute the total phase change of $E(\lambda)$ as $\lambda$ is varied around the boundary of $K$ (top right panel). The phase change is seen to be $-4\pi$ so that the winding number is $-2$, suggesting (along with previous computations for $D =1$) that there are two poles and no roots inside $K$. The bottom panels of Figure \ref{fig:dwind} provide visual confirmation that for $D > 1$, the winding number is $-2$. This means that there are two more poles in $K$ than there are roots. Since when $D=1$ this is also the case, and we did not find evidence of a root and pole pair entering the region $K$ as we varied $D$, we conclude numerical stability of the standing fronts when $D>1$. We numerically confirmed the existence of a pair of poles at complex conjugate points $(x + i y) \approx (0.092  \pm1.448i)$.

\subsection{The general case: nonstationary waves, $v$ diffuses, $D>0$, $c\ne0$}

The stability of solutions in this most general case can be investigated using the same techniques as used in earlier parts of this section, except that the calculations become more complicated numerically. For a representative sample of parameters, the Evans function was found to be numerically unreliable, so the conclusions in this subsection rely on calculations of the Riccati-Evans function.

\subsubsection*{Example 9: Nonstationary pulse}

When the reaction terms in equations \eqref{eq:TWsys2} are given by
$$
g(u,v) = 6(u^2 - D^2 v^2) + \frac{(c^2-4)u + (c^2 - 4 D^2)v}{1-D}
$$
with $D\ne0,1$ and $c>0$, solutions were found to be \cite{BHW19}
\begin{equation}\label{eq:travpulse}
\begin{aligned}
\wh{u}(z) & = -\frac{\sech^2(z)(2D \tanh(z) -c)}{c(1-D)}\,, \qquad
\wh{v}(z) & = \frac{\sech^2(z)(2 \tanh(z) -c)}{c(1-D)}. 
\end{aligned}
\end{equation}
The essential spectrum is stable for all parameters values and is given by a pair of left-opening parabolas centred on the real axis, qualitatively like those shown in Figure~\ref{fig:FredholmBoundariesCZeroReal}. The Fredholm index is zero everywhere and the Fredholm borders are the essential spectrum. We determined the point spectrum for a wide range of parameter values using the Riccati-Evans function and always found an eigenvalue in the right-half plane indicating that nonstationary pulses of this type are always unstable.

\subsubsection*{Example 10: Nonstationary front}

When the reaction terms in equations \eqref{eq:TWsys2} are given by
$$
g = \frac{1}{1-D}(u+v)(4Du^2 + 4D^3v^2+8D^2u v - 2 c u - 2 c D^2 v + c^2)\,,
$$
nonstationary front solutions are \cite{BHW19}
\begin{equation}\label{eq:NSFgen}
\begin{aligned}
\wh{u}(z)=\frac{D\tanh^2z-c\tanh z-D}{c(D-1)}\,,\qquad
\wh{v}(z)=-\frac{\tanh^2z-c\tanh z-1}{c(D-1)}.
\end{aligned}
\end{equation}
One can readily show that stability of the essential spectrum requires $(c-2)(C-2D)\ge0$ and $(c+2)(C+2D)\ge0$. However, even when these conditions are satisfied we have numerically determined the point spectra for a wide range of values of $c$ and $D$ and in all cases a root with positive real part was found indicating that these solutions are always unstable.

\begin{figure}
    \includegraphics[scale=0.25]{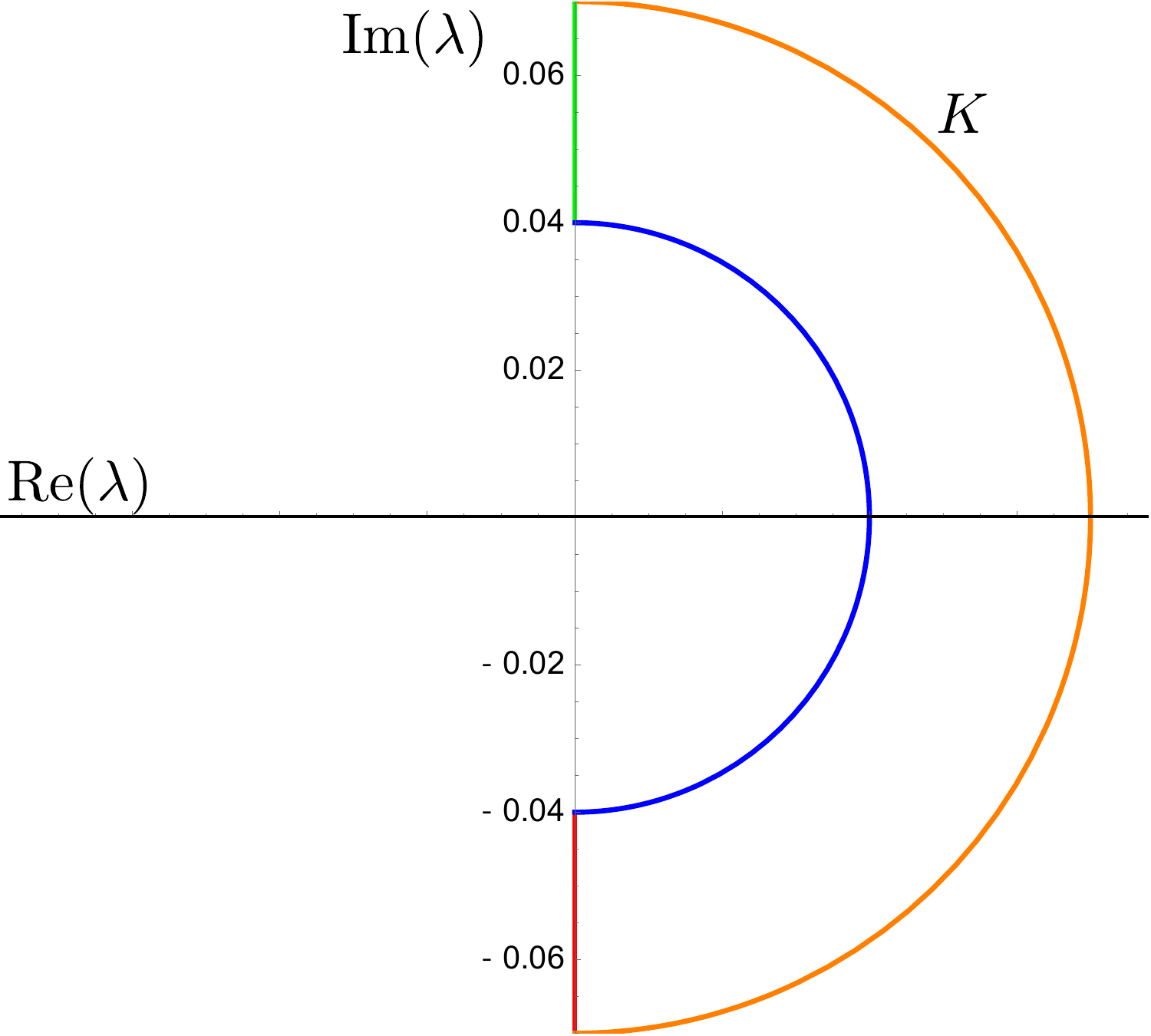} \quad
        \includegraphics[scale=0.25]{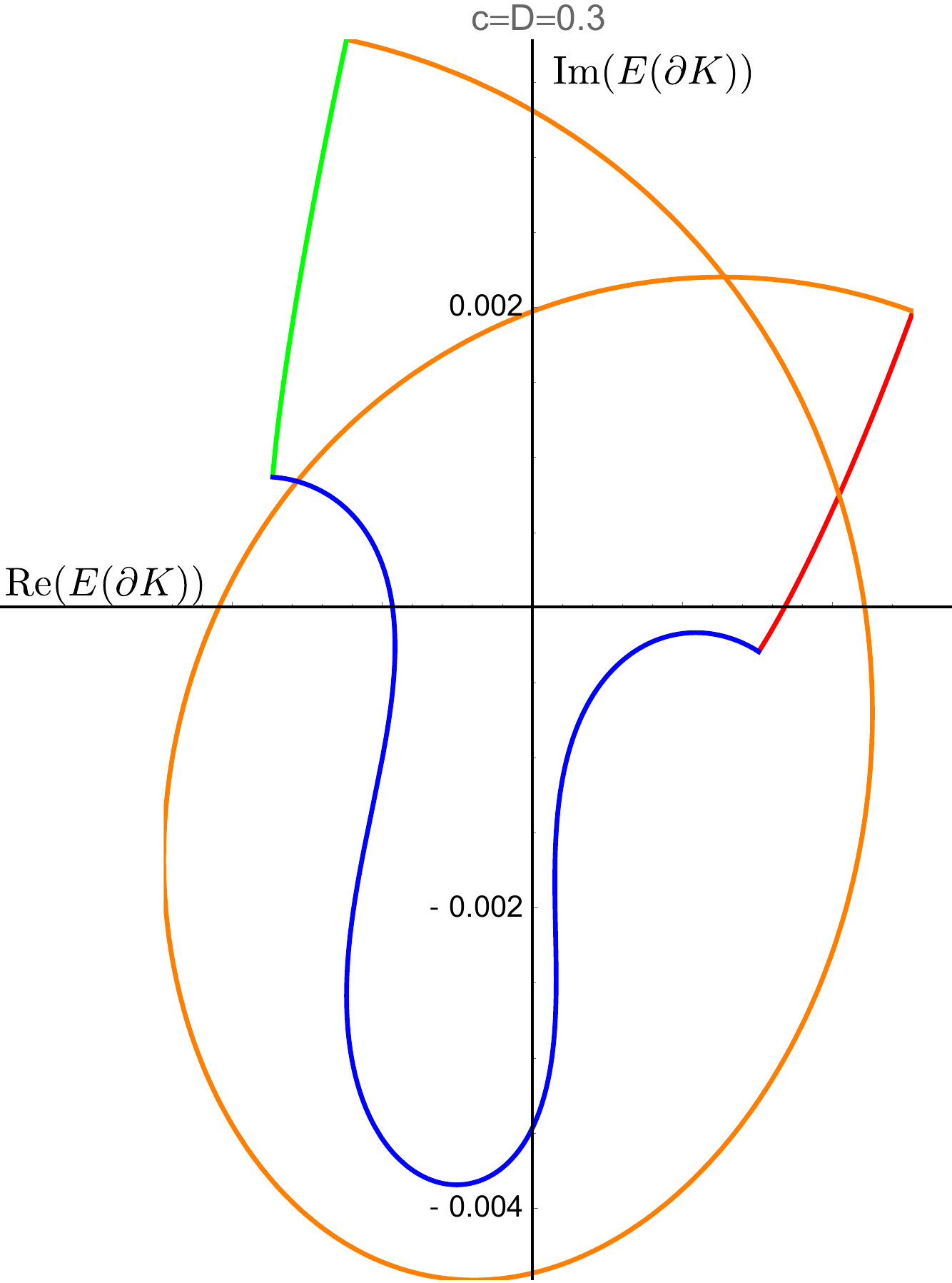}\\
    \includegraphics[scale=0.25]{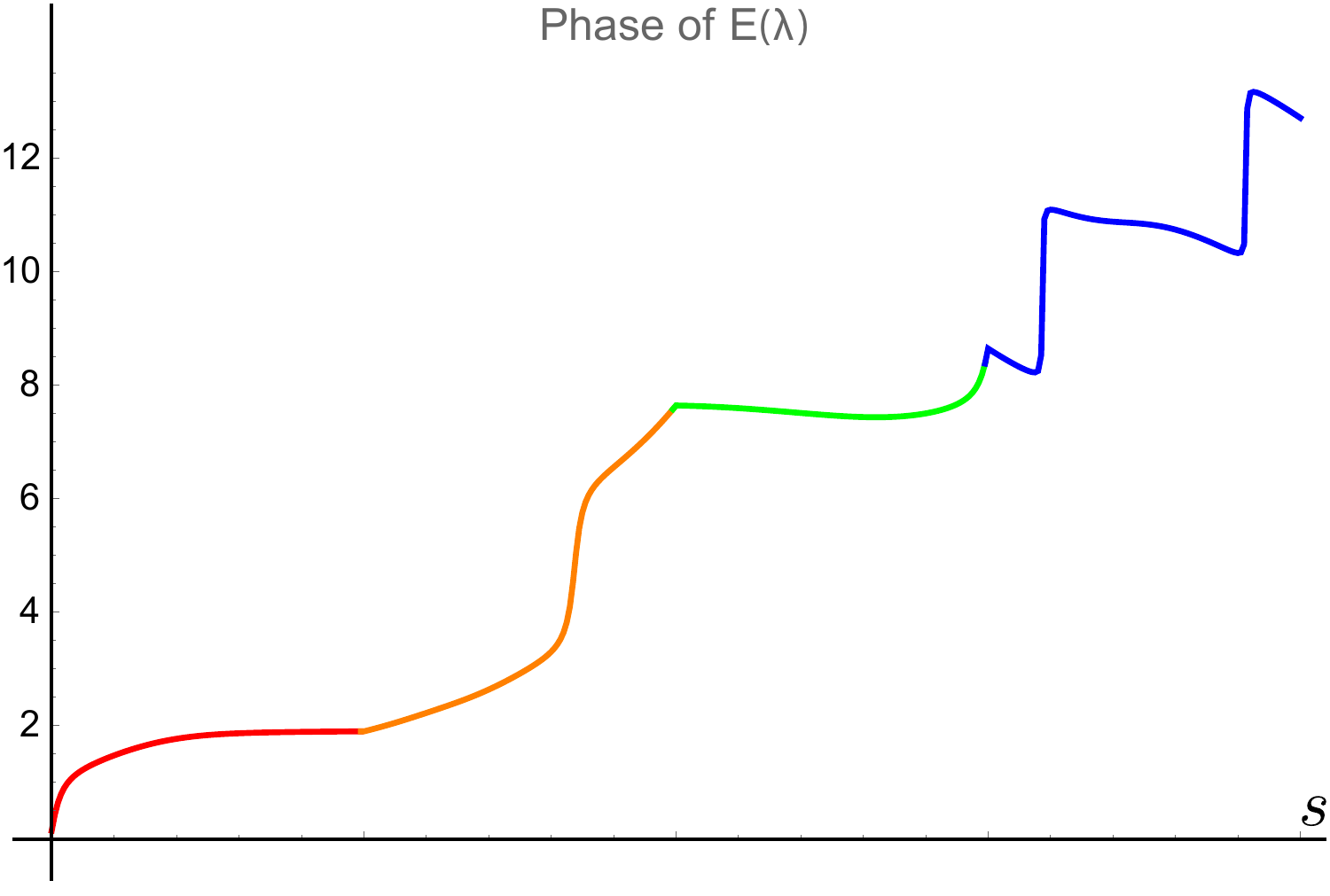}
    \caption{Top right is the plot of Evans function for the travelling front for the region $K$ in the top left. These parameter values have stable essential spectrum ($c=0.3=D$). As seen in the bottom figure, there is a root in the right half plane showing that the wave is nevertheless unstable. The middle figure is the winding of the phase of $E(\lambda)$}
    \label{fig:unsttravfront}
\end{figure}

\section{Summary of results and Discussion}
In this section, we summarize all of the various theoretical and numerical results that we have obtained and discuss their importance in the context of applications.

We have determined the equations that govern the stability of travelling and stationary pulses and fronts for a system of reaction diffusion equations with degenerate source terms given by \eqref{eq:sys1}. 
For $D\ne0$ the stability is determined by a fourth-order eigenvalue problem. However, for $D=0$ the highest derivative term is absent and for $D=0$ and $c\ne0$ the system reduces to a third-order eigenvalue problem indicating that $D=0$ is a singular limit. Moreover, if $D=c=0$ then the two highest derivative terms are absent and the system is doubly singular. Consequently, these two cases must be treated separately. In fact, as we have shown, these singular cases have completely different behaviour from the case $D\ne0$. The case $D=1$ is also special in the sense that it can be integrated twice and this reduces to a second-order eigenvalue problem. This case also has very different stability results from the generic case. 

For the generic problem, we have determined the Fredholm borders and obtained the large $\lambda$ asymptotics of the Fredholm index, and this has allowed us to determine the form of the essential spectrum. In particular, we have obtained the necessary and sufficient conditions for the essential spectrum to be stable. Moreover, the most unstable mode of the essential spectrum can be either  oscillatory or non-oscillatory, and we have determined conditions that delineate these two cases.
We have also shown that for stationary waves ($c=0$) the essential spectrum consists only of the Fredholm borders.  

We have also found a number of surprising theoretical results if $D=0$ or $D=1$. For $D=1$ we have shown that pulses are very different from fronts. Firstly, stationary pulses are necessarily unstable while non-stationary pulses cannot exist. Secondly, stationary and non-stationary fronts are necessarily stable if they are monotonic and unstable otherwise. Thirdly, we have constructed examples of stationary pulses, stationary fronts, and non-stationary fronts for which we can explicitly determine the Evan's function and the entire point spectrum. 

In the case of $D=0$ we have shown that stationary pulses and fronts can only be made up of piecewise constant solutions. Moreover, if any nontrivial piecewise solution exists then there an infinite family of solutions must also exist. For any solution we have derived a necessary condition for the point spectrum to be real. We have also shown that both pulses and fronts can be stable or unstable. We have constructed examples of pulses and fronts for which the Evan's function and the point spectrum can be computed explicitly. In fact, it is always possible to obtain the Evan's function explicitly but, in general, obtaining explicit expressions for the spectrum will not be possible for solutions with more than two jumps. The main theoretical results are summarized in Table 1.

In the remaining cases, one must make use of numerical methods to compute the point spectrum. We used a numerical method that is based on the Riccati-Evan's function. This method proved to be much more robust and reliable than methods that use the traditional Evan's function approach and allowed us to examine a very broad range of parameter values that would otherwise have been problematic. We used our numerical method to determine the stability of the point spectrum for all known explicit pulse and front solutions of \eqref{eq:sys1}. Our numerical results illustrated that all known pulse solutions are unstable. On the other hand, stable non-stationary front solutions were found in the case $D=0$ while stable stationary fronts were found in the case $D\ne0,1$. In addition, all known nonstationary front solutions with $D\ne0,1$ are unstable. These results are summarized in Table 2.

The results for $D=0$ have a number of important consequences for the modelling of ion transport in cellular matter using equations of the type first developed by \cite{Tuck1978}. There has been an ongoing discussion in the literature regarding the minimal model required to generate nonstationary ionic waves in which the ionic concentrations return to their initial values after the wave had passed. Waves of this type are characteristic of many biological processes, including cortical spreading depression \cite{Tuck1978}. Miura et al. \cite{MHW07} presented an intuitive argument that suggested that to obtain non-stationary pulse solutions, the minimal model with degenerate source terms must have at least two ionic species (corresponding to four equations) and that pulse solutions could not exist for systems with a single ionic species (corresponding to two equations) of the form \eqref{eq:sys1}. Moreover, despite trying multiple different functional forms for the reaction terms, direct numerical simulations of \eqref{eq:sys1} never generated pulse solutions. In order to address this question, \cite{BHW19} derived an explicit condition on the source term that guarantees the existence of pulse solutions and determined a parametric family of exact solutions. This seemed to suggest that travelling pulses could be generated by equations of the form \eqref{eq:sys1} after all. However, in this paper, we have shown that even though stationary pulse solutions and travelling front solutions of \eqref{eq:sys1} can be stable, all of the non-stationary pulse solutions obtained by \cite{BHW19} are unstable. Therefore, this work suggests that the TM model may, after all, be the minimal model that can give {\it stable} travelling pulse solutions.

\begin{table}[ht]
    \centering
    \caption{Theoretical Results}
    \begin{tabular}{|l||l|l|}
        \toprule
          & \textbf{Pulse} & \textbf{Front}\\
        \midrule\midrule
        $D=1$, $c=0$  & Unstable & Stable if front is monotone\\ 
                      &          &\quad and unstable otherwise\\
                      & & \\
                      & Example with explicit Evan's & Example with explicit Evan's  \\
                      &\quad function and spectrum & \quad function and spectrum\\
                      & & \\
                      &(Example 1)&(Example 2)\\
        \midrule
        $D=0$, $c=0$  & Only piecewise solutions exist & Only piecewise solutions exist\\ 
                    & & \\
                      &Can be stable or unstable &Can be stable or unstable \\
                      & & \\
                      & Example with explicit  & Example with explicit   \\
                      &\quad Evan's function& \quad Evan's function\\
                      & & \\
                      &(Example 3)&(Example 4)\\
        \midrule
        $D=1$, $c\ne0$ & No solutions exist & Stable if front is monotone\\ 
                      &          &\quad and unstable otherwise\\
                      & & \\
                      &(Equation \eqref{D1PulseStationary})&(Equation \eqref{DvSatisfiesLambda0} and \cite{KP13} )\\
        \bottomrule
    \end{tabular}
\end{table}

\begin{table}[ht]
    \centering
    \caption{Numerical Results}
    \begin{tabular}{|l||l|l|}
        \toprule
          & \textbf{Pulse} & \textbf{Front}\\
        \midrule
        \midrule
        $D=0$, $c\ne 0$  & All examples unstable & Stable examples exist\\ 
        & &\\
        &(Example 5, Figure \ref{fig:algwave2})&(Example 6, Figures \ref{fig:sechfredholm} \& \ref{fig:sechricstabimag})\\
        \midrule
        $D\ne 0, 1$, $c=0$  & All examples unstable & Stable examples exist\\
        & &\\
        &(Example 7, Figure \ref{fig:ricevstanding})&(Example 8, Figures \ref{fig:ricevstandfront} \& \ref{fig:dwind})\\
        \midrule
        $D\ne 0, 1$, $c\ne 0$  & All examples unstable & All examples unstable\\
        & &\\
        &(Example 9)&(Example 10)\\ %, Figures \ref{fig:tfessential} )\\
        \bottomrule
    \end{tabular}
\end{table}

\section*{Acknowledgments}

This work was supported by Australian Research Council Discovery Projects DP200102130 to RM and BHBH, DP210101102 to RM, and DP210102246 to BHBH. JJW acknowledges support from the Sydney Mathematical Research Institute International Visitor Program. JJW was supported by the Research Grants Council of Hong Kong Special Administrative Region, China (CityU 11309422).

\newpage
%%%%%%%%%%%%%%%%%%%%%%%%%%%%%%%%%%%%%%%%%%%%%%%%%%%%%%%%%%%%
\bibliographystyle{alpha}
\bibliography{stab} 
%%%%%%%%%%%%%%%%%%%%%%%%%%%%%%%%%%%%%%%%%%%%%%%%%%%%%%%%%%%%

\end{document}